\documentclass[11pt]{article}
\usepackage{amsmath, amssymb, epsfig, eepic}

\numberwithin{equation}{section}

\newcommand{\vep}{\varepsilon}
\newcommand{\eps}{\varepsilon}

\newcommand{\qed}{\hfill $\Box$}
\newcommand{\proof}{\noindent {\bf Proof}. \hspace{2mm}}
\newcommand{\prob}{\mathbb P}
\newcommand{\Q}{\mathbb Q}
\newcommand{\expec}{\mathbb E}
\newcommand{\smallsup}[1] {{\scriptscriptstyle{({#1}})}}
\newtheorem{theorem}{Theorem}[section]
\newtheorem{lemma}[theorem]{Lemma}
\newtheorem{prop}[theorem]{Proposition}
\newtheorem{corr}[theorem]{Corollary}
\newtheorem{remark}[theorem]{Remark}

\newcommand{\eq}{\begin{equation}}
\newcommand{\en}{\end{equation}}

\newcommand{\sN} {{\scriptscriptstyle{N}}}
\newcommand{\sss} {\scriptscriptstyle }

\title{Distances in random graphs with finite variance degrees}
\author{Remco van der Hofstad\footnote{Department of Mathematics and
Computer Science, Eindhoven University of Technology, P.O.\ Box
513, 5600 MB Eindhoven, The Netherlands. E-mail: {\tt
rhofstad@win.tue.nl}}\\ Gerard Hooghiemstra and Piet Van
Mieghem\footnote{Delft University of Technology, Electrical Engineering,
Mathematics and Computer Science P.O. Box 5031, 2600 GA Delft, The
Netherlands. E-mail: {\tt G.Hooghiemstra@ewi.tudelft.nl,
P.VanMieghem@ewi.tudelft.nl}}}

\begin{document}
\maketitle

\begin{abstract}
In this paper we study a random graph with $N$ nodes, where
node $j$ has degree $D_j$ and $\{D_j\}_{j=1}^N$ are i.i.d.\
with $\prob(D_j\leq x)=F(x)$. We assume that
$1-F(x)\leq c x^{-\tau+1}$ for some $\tau>3$ and
some constant $c>0$. This graph model is a variant of the
so-called configuration model, and includes heavy tail degrees
with finite variance.

The minimal number of edges between two arbitrary connected nodes, also known as the
graph distance or the hopcount, is investigated when $N\rightarrow \infty$.
We prove that the graph distance grows like $\log_{\nu}N$, when the base
of the logarithm equals $\nu=\expec[D_j(D_j -1)]/\expec[D_j]>1$. This confirms the
heuristic argument of Newman, Strogatz and Watts \cite{NSW00}. In
addition, the random fluctuations around this asymptotic mean $\log_{\nu}{N}$
are characterized and shown to be uniformly bounded. In particular, we
show convergence in distribution of the centered graph distance along
exponentially growing subsequences.
\end{abstract}


\section{Introduction}
The study of complex networks plays an increasingly important role
in science. Examples of such networks are electrical power
grids and telephony networks, social relations, the World-Wide Web
and Internet, co-authorship and citation networks of scientists,
etc. The structure of these networks affects their performance.
For instance, the topology of social
networks affects the spread of information and disease (see e.g.,
\cite{Strogatz_Nature01}). The rapid evolution in, and the success
of, the Internet have incited fundamental research on the topology
of networks.

Different scientific disciplines report their
own viewpoints and new insights in the broad area of networking.
In computer science and electrical engineering, massive Internet
measurements have lead to fundamental questions in the modelling
and characterization of the Internet topology
\cite{FFF99,Tangmunarunkit_sigcom02}. These modelling questions drive
the understanding of the Internet's complex behavior and allow to
plan and to control end-to-end communication. The pioneering work of
Strogatz and Watts (see e.g.\ \cite{Strogatz_Nature01,Watts}
and the references therein) have triggered an immense number of
research papers in the field of theoretical physics. Strogatz and
Watts proposed `small world networks' and illustrated how such small
worlds can arise due to underlying mechanisms in different practical
networks such as social networks, growing structures in nature, the
Web, etc.

Albert and Barab\'asi in \cite{AB99}
showed that preferential attachment of nodes gives rise to a
class of graphs often called `scale free networks'. See also
\cite{AB02, Barabasi_linked} and the references therein.
Scale free networks seem to explain the structure of the World-Wide Web,
the autonomous domain structure of Internet, citation graphs and many other
complex networks (see e.g., \cite{AB02, Newm03}). The essence of scale
free networks is that the nodal degree is a power law, or, alternatively,
heavy-tailed, meaning    that the number of nodes with degree
equal to $k$ is proportional to
$k^{-\tau}$ for some power exponent $\tau>1$.
On the World-Wide Web, it has indeed been shown that there are power law
degree sequences, both for the in- and out degrees
(see \cite{BKMRRSTW, KRRT99}). The work of Albert and Barab\'asi
have inspired substantial work on scale-free graphs and can be
seen as a way to understand the emergence of power law degree sequences.
In the model by
Albert and Barab\'asi \cite{AB99}, this power exponent is restricted
to $\tau=3$ \cite{BRST01}, but in refinements of the model, different
values of $\tau$ can be obtained. See, e.g., \cite{ACL01b, BBCR03, CF03, KRRSTU00}
and the references therein. We will comment on the relations between our work
and preferential attachment models in Section \ref{sec-RW} below.
For an overview of the extensive field of random
graphs, we refer to the books of Bollob\'as \cite{bol01} and
Janson {\em et al.}  \cite{jans}.

The current paper presents a rigorous mathematical derivation
for the random fluctuations of the graph distance between two
arbitrary nodes in a graph with finite variance degrees.
These finite variance degrees include power laws with power
exponent $\tau >3$. We consider the configuration model
with power law degree sequences, a variation on
a model originally proposed by Newman, Strogatz and Watts
\cite{NSW00}, prove their conjecture and proceed beyond their
results by combining coupling theory, branching processes and
shortest path graphs.


\subsection{Model definition}
Fix an integer $N$. Consider an i.i.d.\ sequence $D_1,D_2,\ldots,D_{\sN}$.
We will construct an undirected graph with $N$ nodes where node $j$ has degree $D_j$.
We will assume that $L_{\sN}=\sum_{j=1}^N D_j$ is even. If $L_{\sN}$ is odd, then we add a
stub to the $N^{\rm th}$ node, so that $D_{\sN}$ is increased by 1.
This single stub will make hardly any difference in what follows,
and we will ignore this effect. We will later specify the distribution of
$D_1$.

To construct the graph, we have $N$
separate nodes and incident to node $j$, we have $D_j$ stubs. All stubs
need to be connected to build the graph. The stubs are
numbered in a given order from $1$ to $L_{\sN}$. We start by
connecting at random the first stub with one of the $L_{\sN}-1$
remaining stubs. Once paired, two stubs form a single edge of the
graph. Hence, a stub can be seen as the left or the right half of
an edge. We continue the procedure of randomly choosing and
pairing the stubs until all stubs are connected.
Unfortunately, nodes having self-loops may occur. However,
self-loops are scarce when $N \to \infty$.

We now specify the degree distribution we will investigate in
this paper. The probability mass
function and the distribution function of the nodal degree $D$ are
denoted by
    \begin{equation}
    \label{kansen}
    \prob(D=j)=f_j,\quad j=0,1,2,\ldots, \quad \mbox{and} \quad
    F(x)=\sum_{j=0}^{\lfloor x \rfloor} f_j,
    \end{equation}
where $\lfloor x \rfloor$ is the largest integer smaller than or
equal to $x$. Our main assumption is that for some $\tau >3$  and
some positive constant $c$,
    \begin{equation}
    \label{distribution}
    1-F(x)\leq c x^{-\tau+1}, \qquad (x>0).
    \end{equation}
This condition implies that the second moment of $D$ is finite.
The often used condition that $1-F(x)= x^{-\gamma+1}L(x), \, \gamma>3,$ with $L$ a
slowly varying function is covered by (\ref{distribution}),
because by Potter's Theorem \cite[Lemma 2, p.~277]{fellerb}, any slowly varying
function $L(x)$ can be bounded above and below by an arbitrary small
power of $x$, so that (\ref{distribution}) holds for any $\tau<\gamma$.

The above model is closely related to the so-called {\it configuration model},
in which the degrees of the nodes are often assumed to be fixed (rather than i.i.d.). See
\cite[Section 4.2.1]{Newm03} and the references therein. We will review some
results proved for the configuration model in Section \ref{sec-RW} below.


\subsection{Main results}

We denote
\begin{equation}
    \label{nu}
    \mu=\expec[D], \qquad \nu=\frac{\expec[D(D-1)]}{\expec[D]},
\end{equation}
and we define the distance or hopcount $H_{\sN}$
between the nodes $1$ and $2$ as the minimum number of edges that
form a path from $1$ to $2$ where, by convention, the distance equals $\infty$ if nodes $1$
and $2$ are not connected. Since the nodes are exchangeable, the
distance between two randomly chosen nodes is equal in distribution to
$H_{\sN}$. Our main result is the following theorem:

\begin{theorem}[Limit law for the typical nodal distance]
\label{thm-tau>3}
Assume that $\tau>3$ in (\ref{distribution}) and that
$\nu>1$. For $k\geq 1$, let $a_k=\lfloor \log_{\nu}k \rfloor-\log_{\nu}k \in (-1,0]$.
There
exist random variables $(R_{a})_{a\in (-1,0]}$ such that as $N\to \infty$,
\begin{equation}
    \label{hop}
    \prob\big(H_{\sN}- \lfloor \log_{\nu} N\rfloor = k\big|H_{\sN}<\infty\big)
    = \prob(R_{a_{\sN}}=k)+o(1), \qquad k\in {\mathbb Z}.
\end{equation}
\end{theorem}

In words, Theorem \ref{thm-tau>3} states that for $\tau>3$, the
graph distance $H_{\sN}$ between two randomly chosen connected nodes
grows like the $\log_{\nu}{N}$, where $N$ is the size of the graph, and that
the fluctuations around this mean remain uniformly bounded in $N$.
Theorem \ref{thm-tau>3} proves the conjecture in Newman, Strogatz
and Watts \cite[Section II.F, (54)]{NSW00}, where
a heuristic is given that the number of edges between
arbitrary nodes grows like $\log_{\nu} N$.
In addition, Theorem \ref{thm-tau>3} improves upon that conjecture
by specifying the fluctuations around the value $\log_{\nu} N$.

We will identify the laws of $(R_{a})_{a\in (-1,0]}$ in Theorem \ref{thm-ll} below. Before doing so, we state
two consequences of the above theorem:

\begin{corr}[Convergence in distribution along subsequences]
\label{cor-weak}
Fix an integer $N_1$. Under the assumptions in
Theorem \ref{thm-tau>3}, and conditionally on $H_{\sN}<\infty$,
along the subsequence $N_k=\lfloor N_1 \nu^{k-1}\rfloor,$ the sequence of random variables
$H_{\sss N_k}- \lfloor \log_{\nu} {N_k}\rfloor$ converges in distribution to $R_{a_{\sss N_1}}$
as $k\rightarrow \infty$.
\end{corr}
Simulations illustrating the
convergence in Corollary \ref{cor-weak}
are discussed in Section \ref{sec-sim}.

\begin{corr}[Concentration of the hopcount]
Under the assumptions in Theorem \ref{thm-tau>3},
\begin{itemize}
\item[{\rm (i)}] with probability $1-o(1)$ and conditionally
on $H_{\sN}<\infty$, the random variable $H_{\sN}$
is in between $(1\pm \vep)\log_{\nu} N$ for any $\vep>0$;
\item[{\rm (ii)}] conditionally on $H_{\sN}<\infty$, the random variables
$H_{\sN}- \log_{\nu} {N}$ form a tight sequence, i.e.,
    \begin{equation}
    \label{hoptight}
    \lim_{K\rightarrow \infty}
    \limsup_{N\rightarrow \infty} \prob\big(|H_{\sN}- \log_{\nu} N|\leq K\big|H_{\sN}<\infty\big)
    =1.
    \end{equation}
\end{itemize}
\end{corr}

We need a limit result from branching process theory before we can
identify the limiting random variables $(R_{a})_{a\in (-1,0]}$.
In Section \ref{sec-BP} below, we introduce a
{\it delayed} branching process  $\{ {\cal Z}_k\}$, where in the
first generation, the offspring distribution is chosen according to
(\ref{kansen}) and in the second and further generations, the
offspring is chosen in accordance to $g$ given by
    \begin{equation}
    \label{outgoing degree}
        g_j=\frac{(j+1) f_{j+1}}{\mu},\quad j=0,1,\ldots.
    \end{equation}
The process $\{{\cal Z}_k/\mu\nu^{k-1}\}$ is a martingale with
uniformly bounded expectation and consequently
converges almost surely to a limit:
    \begin{equation}
    \label{martin}
    \lim_{n\to \infty}\frac{{\cal Z}_n}{\mu\nu^{n-1}}
    = {\cal W}\qquad a.s.
    \end{equation}
In the theorem
below we need two independent copies ${\cal W}^\smallsup{1}$ and
${\cal W}^\smallsup{2}$ of ${\cal W} $.

\begin{theorem}[The limit laws]
\label{thm-ll}
Under the assumptions in Theorem \ref{thm-tau>3}, and for $a\in (-1,0]$,
\begin{equation}
    \label{Zdistr2}
    \prob(R_{a}>k) = \expec\big[\exp\{-\kappa\nu^{a+k}
    {\cal W}^\smallsup{1}{\cal W}^\smallsup{2}\}
    \big|{\cal W}^\smallsup{1}{\cal W}^\smallsup{2}>0
    \big],
\end{equation}
where ${\cal W}^\smallsup{1}$ and $ {\cal W}^\smallsup{2}$ are
independent limit copies of ${\cal W}$ in (\ref{martin}) and where
$
\kappa=\mu(\nu-1)^{-1}.
$
\end{theorem}
We will also provide an error bound of the convergence stated in Theorem \ref{thm-tau>3}.
Indeed, we show that for any $\alpha>0$, and for all $k\leq \eta \log_\nu N$ for some
$\eta>0$ sufficiently small,
    \begin{equation}
    \label{uitspraakrep}
    \prob(H_{\sN}>\lfloor \log_{\nu} N\rfloor+k)
    =\expec\big(\exp\{-\kappa \nu^{a_{\sN}+k}{\cal W}^\smallsup{1}
    {\cal W}^\smallsup{2}\}
    \big)+O((\log N)^{-\alpha}).
    \end{equation}
Unfortunately, due to the conditioning in Theorem \ref{thm-tau>3},
it is hard to obtain an explicit error bound in (\ref{hop}).

The law of $R_{a}$ is involved, and can in most cases not be computed
exactly. The reason for this is the fact that the random variables
${\cal W}$ that appear in its statement are hard to compute explicitly.
For example, for the power-law degree graph with $\tau>3$, we do not
know what the law of ${\cal W}$ is. See also Section \ref{sec-BP}.
There are two examples where the law of ${\cal W}$ is known.
The first is when all degrees in the graph are equal to some $r>2$,
and we obtain the $r$-regular graph (see also \cite{BdlV82}, where the diameter
of this graph is studied). In this case, we have that $\mu=r, \nu=r-1$, and ${\cal W}=1$
a.s. In particular, $\prob(H_{\sN}<\infty)=1+o(1).$ Therefore, we obtain that
\begin{equation}
    \label{Zdistrex1}
    \prob(R_{a}>k) = \exp\{-\frac{r}{r-2}(r-1)^{a+k}\},
\end{equation}
and $H_{\sN}$ is asymptotically equal to $\log_{r-1} N$. The second example is when the law $g$ is
geometric, in which case the branching process with offspring $g$ conditioned to be positive
converges to an exponential random variable with parameter 1. This example corresponds to
    \begin{equation}
        g_j=p(1-p)^{j-1}, \quad \text{so that }\quad f_j = \frac{1}{j c_p} p(1-p)^{j-2},\quad \forall j\geq 1,
    \end{equation}
and $c_p$ is the normalizing constant. For $p>\frac 12$, the law of ${\cal W}$ has the same law as the sum of
$D_1$ copies of a random variable ${\cal Y}$, where ${\cal Y}=0$ with probability $\frac{1-p}{p}$
and equal to an exponential random variable with parameter 1 with probability $\frac{2p-1}p$. Even in
this simple case, the computation of the exact law of $R_a$ is non-trivial.
Although the laws $R_{a}$ are hard to compute exactly, Theorems \ref{thm-tau>3} and \ref{thm-ll}
make it possible to simulate the hopcount in random graphs of arbitrary size since the law of
${\cal W}$ is simple to approximate numerically, for example using Fast Fourier Transforms.

In \cite{HHVexpec}, the expected value of the random variable $R_a$ is computed numerically,
by comparing it to $\expec[\log {\cal W}|{\cal W}>0]$.
One would expect that for some $\beta$ with $0<\beta <\alpha$,
    \eq
    \label{convexp}
    \expec[H_{\sN}|H_{\sN}<\infty]= \lfloor \log_{\nu} N\rfloor +\expec[R_a] +O((\log{N})^{-\beta}).
    \en
If so, an accurate computation of $\expec[R_a]$ would yield the fine asymptotics
of the expected hopcount, and this would yield an extension of the conjectured
results in \cite[(54)]{NSW00}. Our methods stop short of proving (\ref{convexp}),
and this remains an interesting
question.

Our final result describes the size of the largest connected component and the maximal
size of all other connected components. In its statement, we write $G$ for the random graph
with degree distribution given by (\ref{kansen}), and we write $q$ for the survival
probability of the delayed branching process $\{{\cal Z}_k\}$ described above.
Thus, $1-q$ is the extinction probability of the branching process.

\begin{theorem}[The sizes of the connected components]
    \label{thm-cluster}

    With probability $1-o(1)$, the largest connected component in $G$ has
    $qN(1+o(1))$ nodes, and there exists $\gamma<\infty$
    such that all other connected components have at most $\gamma \log N$ nodes.
\end{theorem}


\subsection{Methodology and heuristics}
\label{sec-MH}
One can understand Theorems \ref{thm-tau>3} and \ref{thm-ll}
intuitively as follows. Denote by $Z_k^{\smallsup{1}}$, respectively,
$Z_k^{\smallsup{2}}$ the number of stubs of nodes at distance $k-1$ from node
$1$, respectively, node 2 (see Section \ref{sec-coupling} for the precise
definitions). Then for $N\to\infty$, the random process
$Z_1^{\smallsup{i}},Z_2^{\smallsup{i}},
\ldots,Z_k^{\smallsup{i}}$, which will be called shortest path
graphs (SPG's), behave as a delayed branching process as long as
$Z_k^{\smallsup{i}}$ is of small order compared to $N$. Thus,
the local neighborhood of the node $i$ is close in distribution
to a branching process.

We sample the stubs uniformly from all stubs and thus, for
large $N$, we attach the stubs to the SPG
proportionally to $jf_j$. Moreover, when a new
stub is attached to the SPG, the chosen stub is used to
attach the new node and forms an edge together with the present stub.
Therefore, the number of stubs of the freshly chosen
node decreases by one and is equal to $j$ if the number of stubs
of the chosen node was originally equal to $j+1$. This motivates
(\ref{outgoing degree}).

The offspring of the node 1 is distributed as $D_1$,
whereas the offspring distribution of
$Z_2^{\smallsup{1}},Z_3^{\smallsup{1}},\ldots$ has (for $N\to
\infty$) probability mass function (\ref{outgoing degree}).
Consequently, as noted in \cite[(51)]{NSW00}, the mean number of
free stubs at distance $k$ is close to $\mu \nu^{k-1}$,
where $\nu=\sum_{j=1}^\infty jg_j$ is defined in (\ref{nu}).
Moreover, a stub in $Z_k^{\smallsup{1}}$ is attached with a positive
probability to a stub in $Z_k^{\smallsup{2}}$ whenever
$Z_{k}^{\smallsup{1}}Z_{k}^{\smallsup{2}}$ is of order $L_{\sN}$.
The total degree $L_{\sN}$ is proportional to $N$
by the law of large numbers, because $\mu=\expec[D_1]<\infty$. Since
both sets grow at the same rate, each has to be of
order $\sqrt{N}$. Therefore, $k$ is typically $\frac 12
\log_{\nu}N$, and the typical distance between 1 and 2 is of order
$2k=\log_{\nu}N$. This can be made precise by coupling
$Z_1^{\smallsup{1}},Z_2^{\smallsup{1}},\ldots$
to a branching process ${\hat Z}_1^{\smallsup{1}},{\hat
Z}_2^{\smallsup{1}},\ldots$ having offspring distribution
$g_j^{\smallsup N}$ given by
        \begin{equation}
        g_j^{\smallsup{N}}=\sum_{i=1}^N I[D_i=j+1]\frac{D_i}{L_{\sN}}=
        \frac{j+1}{L_{\sN}}\sum_{i=1}^N I[D_i=j+1],
        \end{equation}
where $I[E]$ is the indicator of the event $E$.
This coupling will be described in Section \ref{subsec-BP1}.
In turn, the branching process ${\hat
Z}_1^{\smallsup{1}},{\hat Z}_2^{\smallsup{1}},\ldots$ will be
coupled, in a conventional way, to a branching process ${\cal
Z}_1^{\smallsup{1}},{\cal Z}_2^{\smallsup{1}},\ldots$ with
offspring distribution $\{g_j\}$ defined in (\ref{outgoing degree}).
The limit result of Theorem \ref{thm-tau>3} and Theorem \ref{thm-ll}  depends on the
martingale limit for super-critical
branching processes with finite mean.

The proof of Theorems \ref{thm-tau>3} and \ref{thm-ll} are based upon
a comparison of the local neighborhoods
of nodes to branching processes. Such techniques are used extensively
in random graph theory. An early example is in \cite{BdlV82},
where the diameter of a random regular graph was investigated.
See also \cite[Chapter 10]{AS00}, where comparisons to branching processes are used to describe the phase transition and the birth of
the giant component for the random graph $G(p,N)$.

The proof of Theorem \ref{thm-cluster} makes essential use of
results by Molloy and Reed \cite{MR95,MR98} for the usual configuration model.
We will now describe their result. When the number of
nodes with degree $i$ in the graph of size $N$ equals $d_i(N)$
where $\lim_{N\rightarrow \infty} d_i(N)/N=Q(i)$, Molloy and Reed \cite{MR95,
MR98} identify the condition $\sum_{i=1}^{\infty} i(i-2)Q(i)>0$ as the
necessary and sufficient condition to ensure that a `giant
component' proportional to the size of the graph exists. By
rewriting the condition $\nu>1$ in Theorem \ref{thm-tau>3} as $\expec[D^2]-2\expec[D]>0$,
we see that a similar condition as in the model of Molloy and Reed is needed
here. To prove Theorem \ref{thm-cluster}, we need to check
that the technical conditions in \cite{MR95,MR98} are satisfied
in our model. In fact, we need to alter the graph $G$ a little bit in order to
apply their results, since in \cite{MR95} it is assumed that no nodes of degree larger than
$N^{\frac14-\epsilon}$ exist for some $\epsilon>0$.


The novelty of our results is that we investigate {\it typical distances}
in random graphs. In random graph theory, it is more customary to investigate
the {\it diameter} in the graph, and in fact, this would also be an interesting
problem. The research question investigated in this paper is inspired by
the Internet. In a seminal paper \cite{FFF99}, Faloutsos {\em et al.} have shown
that the degree distribution of autonomous systems
in Internet follows a power law with power exponent
$\tau\approx 2.2$. Thus, the power law random graph with this
value of $\tau$ can possibly lead to a good Internet model
on the autonomous systems (AS) level (see \cite{FFF99, Tangmunarunkit_sigcom02}).
For the Internet on the more detailed router level, extensive
measurements exist for the hopcount, which is the
number of routers traversed between two typical routers, as well as for
the AS-count, which is the number of autonomous systems traversed
between two typical routers.  To validate the configuration model
with i.i.d.\ degrees,
we intend to compare the distribution of the distance
between pairs of nodes to these measurements in Internet.
For this, a good understanding of the typical distances
between nodes in the degree random graph are necessary, which formed the main
motivation for our work. The hopcount in Internet seems to be close to a
Poisson random variable with a fairly large parameter. In turn, a Poisson random variable
with large parameter can be approximated by a normal
random variable with equal expectation and variance. See e.g.
\cite{Paxs97, vanmieghem} for data
of the hopcount in Internet.

From a practical point of view, there are good reasons to study the
typical distances in random graphs rather than the diameter. For one,
typical distances are simpler to measure, and thus allow for a simpler
validation of the model. Also, the diameter is a number, while the {\it
distribution} of the typical distances contains substantially more information.
Finally, the diameter is rather sensitive to small changes to a graph.
For instance, when adding a string of a few nodes, one can dramatically alter
the diameter, while the typical distances in the graph hardly change. Thus,
typical distances in the graph are more robust to modelling discrepancies.


\subsection{Related work}
\label{sec-RW}
There is a wealth of related work which we will now summarize. The
model investigated here was also studied in \cite{norros}, with
$1-F(x)=x^{-\tau+1}L(x),$ where $\tau\in (2,3)$ and $L$ denotes a slowly varying function. It was shown in
\cite{norros} that the average distance is bounded from above by
$2\frac{\log\log N}{|\log(\tau-2)|}(1+o(1))$.
We plan to return to the question of average distances and
connected component sizes when $\tau<3$ in three future publications \cite{HHZ04a, HHZ04b, HHZ04c}.

There is substantial work on random graphs that are, although different
from ours, still similar in spirit. In \cite{ACL01}, random graphs
were considered with a degree sequence that is {\it precisely}
equal to a power law, meaning that the number of nodes with degree
$k$ is precisely proportional to $k^{-\tau}$. Aiello {\em et al.} \cite{ACL01}
show that the largest connected component is of the order of the size of the graph
when $\tau<\tau_0=3.47875\ldots$, where $\tau_0$ is the solution of
$\zeta(\tau-2)-2\zeta(\tau-1)=0$, and where $\zeta$ is the Riemann Zeta
function. When $\tau>\tau_0$, the largest connected component is
of smaller order than the size of the graph and more precise
bounds are given for the largest connected component. When
$\tau\in (1,2)$, the graph is with high probability connected. The proofs
of these facts use couplings with branching processes and strengthen
previous results due to Molloy and Reed \cite{MR95,MR98} described above.
For this same model, Dorogovtsev {\em et al.}
\cite{DGM01, DGM02} investigate the leading asymptotics and
the fluctuations around the mean of the distance between arbitrary
nodes in the graph from a theoretical physics point of view, using
mainly generating functions.

A second related model can be found in \cite{CL02a} and \cite{CL02b}, where edges
between nodes $i$ and $j$ are present with probability
equal to $w_iw_j/\sum_l w_l$ for some `expected degree vector'
$w=(w_1, \ldots, w_{\sN})$. Chung and Lu \cite{CL02a} show that when $w_i$ is
proportional to $i^{-{\frac 1{\tau-1}}}$ the average distance
between pairs of nodes is $\log_{\nu}N(1+o(1))$ when $\tau>3$, and
$2\frac{\log\log N}{|\log(\tau-2)|}(1+o(1))$ when $\tau\in (2,3)$.
The difference between this model and ours is that the nodes are
not exchangeable in \cite{CL02a}, but the observed phenomena
are similar. This result can be heuristically understood as follows. Firstly,
the actual degree vector in \cite{CL02a} should be close to the
expected degree vector. Secondly, for the expected degree vector, we can compute that
the number of nodes for which the degree is less than or equal to $k$ equals
    $$
    |\{i: w_i\leq k\}|\propto |\{i: i^{-\frac{1}{\tau-1}}\leq k\}|\approx k^{-\tau+1}.
    $$
Thus, one expects that the number of nodes with degree at most $k$ decreases as
$k^{-\tau+1}$, similarly as in our model. In \cite{CL02b}, Chung and Lu study
the sizes of  the connected components in the above model. The advantage of this
model is that the edges are {\it independently} present, which makes the resulting
graph closer to a traditional random graph.

All the models described above are {\it static}, i.e., the size of the graph
is {\it fixed}, and we have not modeled the {\it growth} of the graph.
As described in the introduction, there is a large body of work investigating
{\it dynamical} models for complex networks, often in the context of the World-Wide
Web. In various forms, preferential attachment has been shown to lead to
power law degree sequences. Therefore, such models intend to {\it
explain} the occurrence of power law degree sequences in random graphs.
See \cite{ACL01b, AB99, AB02, BBCR03, BR02, BR03a, BR03b, BRST01,
CF03, KRRSTU00} and the references therein. In the preferential attachment model,
nodes with a fixed degree $m$ are added sequentially. Their stubs are attached
to a receiving node with a probability proportionally to the degree of the receiving
node, thus favoring nodes with large degrees. For this model, it is shown that the
number of nodes with degree $k$ decays proportionally to $k^{-3}$ \cite{BRST01}, the diameter
is of order $\frac{\log{N}}{\log\log{N}}$ when $m\geq 2$ \cite{BR02}, and couplings
to a classical random graph $G(N,p)$ are given for an appropriately chosen $p$ in \cite{BR03b}.
See also \cite{BR03a} for a survey.

It can be expected that our model is a snapshot of the above models,
i.e., a realization of the graph growth processes at the time instant that the graph
has a certain prescribed size. Thus, rather than to describe the growth of the model,
we investigate the properties of the model at a given time instant.
This is suggested in \cite[Section VII.D]{AB02}, and it would be very interesting
indeed to investigate this further mathematically, i.e., to investigate the relation
between the configuration and the preferential attachment models.

The reason why we study the random graphs
at a given time instant is that we are interested in the topology of the random graph.
In \cite{Tangmunarunkit_sigcom02}, and inspired by the observed power law degree
sequence in \cite{FFF99}, the configuration model with
i.i.d.\ degrees is proposed as a model for the AS-graph
in Internet, and it is argued on a qualitative basis that this simple model serves
as a better model for the Internet topology than currently used topology generators.
Our results can be seen as a step towards the quantitative understanding of whether the
hopcount in Internet is described well by the average graph distance in the
configuration model.

In \cite[Table II]{Newm03}, many more examples are given of real
networks that have power law degree sequences. Interestingly, there
are also many examples where power laws are {\it not} observed,
and often the degree law falls off faster than a power law. These
observed degrees can be described by a degree distribution as
in (\ref{kansen}) with $1-F(x)$ smaller than any power, and
the results in this paper thus apply. Such examples are described
in more detail in \cite[Section II]{AB02}. Examples where the
tails of the degree distribution are lighter than power laws are
power and neural networks \cite[Section II.K]{AB02}, where the
tails are observed to be exponential, and protein folding
\cite[Section II.L]{AB02}, where the tails are observed to
be Gaussian. In other examples, a degree distribution is found
that for small values is a power law, but has an exponential
cut off. An example of such a degree distribution is
    \eq
    f_k=C k^{-\gamma} e^{-k/\kappa},
    \en
for some $\kappa>0$ and $\gamma \in {\mathbb R}$. The size of $\kappa$ indicates up to what
degree the power law still holds, and where the exponential
cut off starts to set in. For this example, our results apply
since the exponential tail ensures that (\ref{distribution})
holds for {\it any} $\tau>3$ by picking $c>0$ large enough.
Thus, we prove the conjectures on the expected path lengths
in \cite[(55), (56)]{NSW00} and \cite[Section V.C, (63) and (64)]{AB02}
for this particular model.

\subsection{Simulation for illustration of the main results}
\label{sec-sim} To illustrate Theorem \ref{thm-tau>3}, we have
simulated the random graph with degree distribution $D=\lceil
U^{-\frac{1}{\tau-1}} \rceil$, where $U$ is uniformly distributed over
$(0,1)$ and where for $x\in \mathbb R$, $\lceil x \rceil$ is the
smallest integer greater than or equal to $x$. Thus,
    \begin{eqnarray*}
    1-F(k) =\prob(U^{-\frac{1}{\tau-1}}>k)
    =k^{1-\tau},\quad k=1,2,3,\ldots,
    \end{eqnarray*}
for which $\mu=1+\zeta(\tau-1)$ and $\nu=2\zeta(\tau-2)/\mu$.

\renewcommand{\epsfsize}[2]{0.6#1}
\begin{figure}[t]
\begin{center}
\epsfbox[20 40 576 469]{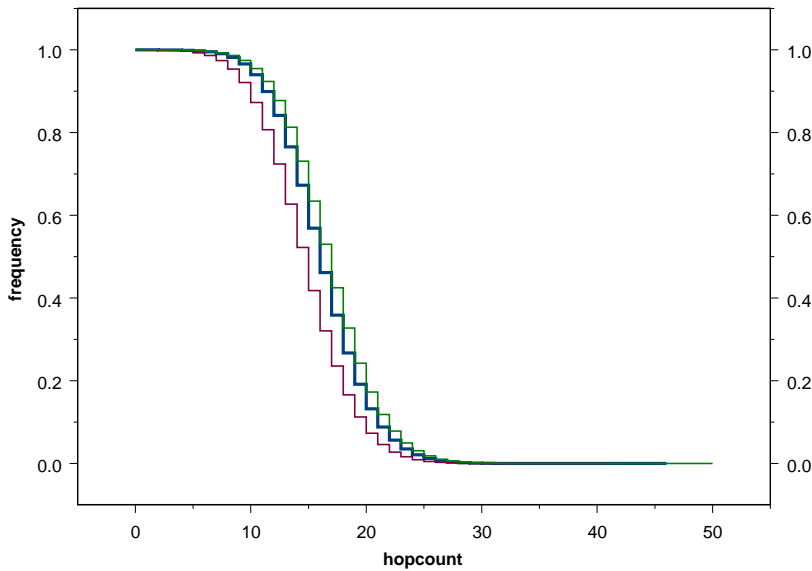} \caption{Empirical
survival functions of the hopcount for $\tau=3.5$ and the
values $N=25,000$, $N=75,000$ (bold) and $N=125,000$, based on
samples of size $1,000$.} \label{survival plots 1}
\end{center}
\end{figure}

We observe that for $\tau=3.5$ and $N=25,000$ and $N=125,000$, the values
$a_{\sN}=-0.62\ldots$ are identical up to two decimals.
We hence expect, on the basis of our main theorem, that
the survival functions $\prob(H_{\sN}>k)$ for these two cases
are similar. Because  $\lfloor \log_{\nu}25,000\rfloor=12$
and $\lfloor \log_{\nu}125,000\rfloor=14$, we expect that
the empirical survival function for $N=125,000$ is a shift
of the empirical survival function for $N=25,000$,
over the horizontal distance $14-12=2$. Figure \ref{survival plots 1}
supports this claim, given the statistical inaccuracy.
In Figure \ref{survival plots 1} we have also
included the empirical survival function for $N=75,000$,
for which $a_{\sN}=-0.99\ldots$, as the bold line. This empirical survival
function clearly has a different shape. Thus, the empirical survival
function for $N=75,000$ is not a shift of the empirical survival
function for $N=25,000$ or $N=125,000$.

\begin{figure}[t]
\begin{center}
\epsfbox[20 40 576 469]{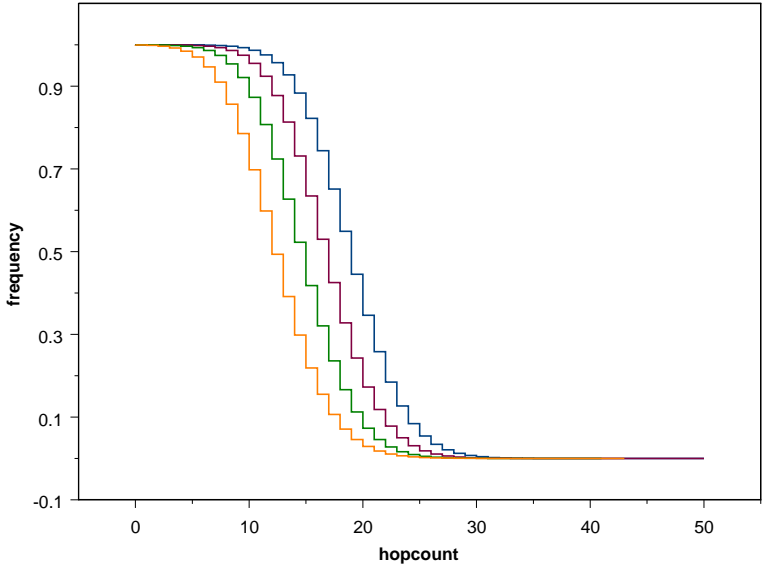} \caption{ Empirical survival
functions of the hopcount for $\tau=3.5$ and the four
 values
$N_k=5,000 \nu^{2k},\, k=0,1,2,3$, based on $1,000$ runs.}
\label{survival plots 2}
\end{center}
\end{figure}

We finally demonstrate Corollary \ref{cor-weak} for $\tau=3.5$
in Figure \ref{survival plots 2}. In this case $\nu^2\approx 5$ and $N_k=N_1\nu^{2k}, \, k=0,1,2,3$.
We take $N_1=5,000$, and so
    $
    N_2=25,000,\, N_3=125,000,\, N_4=625,000.
    $
For these values of $N_1,\ldots,N_4$, we have
simulated the hopcount with $1,000$ replications and we expect from
Corollary \ref{cor-weak} that the survival functions run parallel at mutual distance $2$.

\renewcommand{\epsfsize}[2]{0.6#1}

\subsection{Organization of the paper}
\label{sec-org}
We will first review the relevant literature on branching processes in Section
\ref{sec-BP}. We will then explain how we can couple our degree
model to independent branching processes in Section
\ref{sec-coupling}. This section is also valuable for our coming
paper \cite{HHZ04a}, where we study the case $\tau \in (2,3)$. In particular,
in \cite{HHZ04a}, we will use Lemmas \ref{lem-firststubs} and \ref{lem-stochbds}
and Proposition \ref{prop-caftrep}. The bounds for the coupling are formulated
in Sections \ref{subsec-BP1}, \ref{sec-bdscoupling} and \ref{subsec-BP3}.
In these sections, we will state the results on the coupling that
are needed in the proof of the main results, Theorems \ref{thm-tau>3} and
\ref{thm-ll}. Parts of this section apply more generally, i.e., to
$\tau \in (2,3)$. We prove Theorems \ref{thm-tau>3} and  \ref{thm-ll}
in Section \ref{sec-pftau>3} and Theorem \ref{thm-cluster}
in Section \ref{sec-clusters}.
The technical details of the coupling of $\{{\hat Z^{\smallsup{i}}}_k\}$
to $\{{\cal Z}^{\smallsup{i}}_k\}$ for $i=1,2$ are contained in Section \ref{appendix},
while the details of the coupling of $\{Z_k^{\smallsup{i}}\}$ to
$\{{\hat Z}_k^{\smallsup{i}}\}$ for $i=1,2$ are in Section \ref{proof Prop4.1}.
Finally, we prove that at any fixed time $m$, with probability converging to
1, $Z_m^{\smallsup{i}}={\cal Z}_m^{\smallsup{i}}$ for $i=1,2$ in Section \ref{sec-pf33}.

\section{Review of branching process theory with finite mean}
\label{sec-BP}

Since we rely heavily on the theory of branching processes, we
will briefly review this theory in the case where the expected
value of the offspring distribution is finite. The theory of
branching processes is well understood (see e.g. \cite{atthreya}).

For the formal definition of the delayed branching process (BP) that we consider here, we define a
double sequence $\{X_{n,i}\}_{n\geq 1, i\geq 1}$ of i.i.d.\ random
variables each with distribution equal to the offspring
distribution $\{g_j\}_{j=0}^\infty$, where we recall
        \begin{equation}
        g_j=\frac{(j+1) f_{j+1}}{\mu},\quad j=0,1,\ldots.
        \end{equation}
We further let $X_{0,1}$ have probability mass function
$f$ in (\ref{kansen}), independently from $\{X_{n,i}\}_{n\geq 1,
i\geq 1}$. The BP $\{{\cal Z}_n\}$ is now defined
by ${\cal Z}_0=1$ and
    $$
    {\cal Z}_{n+1}=\sum_{i=1}^{{\cal Z}_n} X_{n,i},\quad n\ge 0.
    $$

Because $\tau>3$, we have that both $\expec[{\cal
Z}_1]=\expec[X_{0,1}]=\mu<\infty$ and $\nu=\expec[X_{1,1}]<\infty$. We
further assume that $\nu=\expec[X_{1,1}]>1$, so that the BP is super-critical.
Given that the $(n-1)^{\rm st}$ generation
consists of $m$ individuals, the conditional expectation of ${\cal
Z}_n$ equals $m\nu,$ independently of the size of the
preceding generations, so that for $n\geq 1$, we have $
\expec[{\cal Z}_n|{\cal Z}_{n-1}]={\cal Z}_{n-1}\nu.$
Hence, ${\cal W}_n=\frac{{\cal Z}_n}{\mu\nu^{n-1}},$
is a martingale. Since $\expec[|{\cal W}_n|]=\expec[{\cal
W}_n]=1$, the sequence $\expec[|{\cal W}_n|]$ is uniformly bounded
by $1$ and so by Doob's martingale convergence theorem \cite[p.
58]{williams} the sequence ${\cal W}_n$ converges almost surely.
If we denote the a.s.\ limit by a proper random variable ${\cal
W}$, we obtain (\ref{martin}).

There are only few examples where the limit random variable ${\cal
W}$ is known. It is known that ${\cal W}$ has an atom at $0$ of size
$p\ge 0$, equal to the extinction probability of the
(delayed-)BP ($q=1-p$). Conditioned on non-extinction the
limit ${\cal W}$ has an absolute continuous density on
$(0,\infty)$.

We need a result that follows from \cite{soren}
concerning the speed of convergence of ${\cal W}_n$ to ${\cal W}$. Define
    $$
    {\cal R}_n=\frac{{\cal W}_n}{\nu} \int_{\nu^n/n^{\alpha}}^\infty
    x\,dG(x), \quad \alpha>0,
    $$
where $G$ is the distribution function of the offspring with
probabilities $\{g_j\}$. Since
    $$
    \mu_{\alpha}=\int_0^{\infty } x [\log^{+} x]^{\alpha}\,
    dG(x)<\infty, \qquad (\log^{+} x=\max(0,\log x)),
    $$
for each $\alpha>0$, it follows from
(\cite[page~8, line~4]{soren}) that with probability 1,
    \begin{equation}
    \label{small-o} {\cal W}-{\cal W}_{k}+\sum_{n=k}^\infty {\cal
    R}_{n}=o(k^{-\alpha}).
    \end{equation}
An immediate consequence of (\ref{small-o}) is that if
$|{\cal W}-{\cal W}_{k}|>k^{-\alpha}$, then $\sum_{n=k}^\infty
{\cal R}_{n}>k^{-\alpha}.$
Hence, using $\expec[{\cal W}_n]=1$ and partial integration,
    \begin{eqnarray*}
    \prob(|{\cal W}-{\cal W}_{k}|>k^{-\alpha})&&\leq
    \prob\left(\sum_{n=k}^\infty {\cal R}_{n}>k^{-\alpha}\right)\leq
    k^{\alpha}\sum_{n=k}^\infty \expec[{\cal R}_{n}]
    =-\sum_{n=k}^\infty \frac{k^{\alpha}}{\nu}\int_{\nu^n/n^{\alpha}}^\infty x\,d\,[1-G(x)]\\
    &&=\sum_{n=k}^\infty
    \frac{k^{\alpha}}{\nu}[1-G(\nu^n/n^{\alpha})] +\sum_{n=k}^\infty
    \frac{k^{\alpha}}{\nu}\int_{\nu^n/n^{\alpha}}^\infty [1-G(x)]\,dx.
    \end{eqnarray*}
Since $1-F(x)\le c\cdot x^{1-\tau}$ (see (\ref{distribution})), we
find $1-G(x)\le c' \cdot x^{2-\tau}$ so that for each
$\alpha>0$, and with $k=\lfloor\frac12 \log_{\nu} N \rfloor$,
\begin{equation}
\label{big-O} \prob\left(|{\cal W}-{\cal W}_{k}|>(\log
N)^{-\alpha}\right)\leq O((\log N)^{\alpha})\sum_{n=k}^\infty
(\nu^n/n^{\alpha})^{3-\tau}=O(e^{-\beta \log N})=O(N^{-\beta}),
\end{equation}
for some positive $\beta$, because $\tau>3$ and $\nu>1$.

\section{Graph construction and coupling with a BP}
\label{sec-coupling} In this section, we will describe how the
shortest path graph (SPG) from node 1 can be obtained, and we will
couple it to a BP. This coupling works for any
degree distribution. In Sections \ref{sec-bdscoupling}
and \ref{subsec-BP3} below, we
will obtain bounds on the coupling.

The SPG from node 1 is the random graph as observed
from node 1, and consists of the shortest
paths between node 1 and all other nodes $\{2, \ldots, N\}$.
As will be shown below, it is not necessarily a tree because
cycles may occur. Recall that two stubs together form an edge.
We define $Z^{\smallsup 1}_1=D_1$, and for $k\ge 2$, we denote by $Z^{\smallsup 1}_k$ the
number of stubs attached to nodes at distance $k-1$ from node 1,
but are not part of an edge connected to a node at distance $k-2$. We will
refer to such stubs as `free stubs'.  Thus, $Z^{\smallsup 1}_k$ is the number
of outgoing stubs from nodes at distance $k-1$.

In Section \ref{subsec-BP1} we will describe a coupling that,
conditionally on $D_1, \ldots, D_{\sN}$, couples $\{Z^{\smallsup 1}_k\}$ to
a BP $\{{\hat Z}^{\smallsup 1}_k\}$ with the {\it
random} offspring distribution
    \begin{eqnarray}
    g_j^{\smallsup{N}}&=& \sum_{i=1}^N I[D_i=j+1] \prob(\mbox{a stub from node
    $i$ is sampled}|D_1, \ldots, D_{\sN}) \nonumber \\
    &=& \sum_{i=1}^N I[D_i=j+1]\frac{D_i}{L_{\sN}}=
    \frac{j+1}{L_{\sN}}\sum_{i=1}^N I[D_i=j+1] \label{gnN},
    \end{eqnarray}
where as before $L_{\sN}=D_1+D_2+\ldots+D_{\sN}$.
By the strong law of large numbers, for $N\to \infty$,
    $$
    \frac{L_{\sN}}{N}\to \expec[D],\quad \mbox{and}\quad \frac1N \sum_{i=1}^N
    I[D_i=j+1]\to \prob(D=j+1),\qquad a.s.
    $$
so that a.s.,
    \begin{equation}
    \label{convergenceoffspring}
    g_j^{\smallsup{N}} \rightarrow (j+1) \prob(D=j+1)/\expec[D]=g_j, \quad N\to \infty.
    \end{equation}
Therefore, the BP $\{{\hat Z}_k^{\smallsup{1}}\}$ with
offspring distribution $\{g_j^{\smallsup{N}}\}$ is expected to be close to a
BP with offspring distribution $\{g_j\}$ given in
(\ref{outgoing degree}). Consequently, in Section \ref{subsec-BP3}, we will couple the BP $\{{\hat
Z}^{\smallsup 1}_k\}$ to a BP $\{{\cal Z}^{\smallsup
1}_k\}$ with offspring distribution $\{g_j\}$. This will allow us to prove Theorems \ref{thm-tau>3}
and \ref{thm-ll} in Section \ref{sec-pftau>3}.

Throughout the paper we use the following lemma. It shows that $L_{\sN}$ is close to $\expec[L_{\sN}]=\mu N$.

\begin{lemma}[Concentration of $L_{\sN}$]
    \label{lem-mom}
    For each $0<a<\frac12$, $b=1-2a$ and some constant $c>0$,
        \begin{equation}
        \prob\left(\left|\frac{L_{\sN}}{{\mathbb E}[L_{\sN}]}-1\right| \geq N^{-a}\right)
        \leq  cN^{-b}.
        \end{equation}
\end{lemma}

\proof The proof is immediate from the Chebychev inequality, since
    $$
    \prob\left(\left(\frac{L_{\sN}}{{\mathbb E}[L_{\sN}]}-1\right)^2 \geq N^{-2a}\right)
            \leq \frac{N^{2a}}{(N\mu)^2}\mbox{Var}(L_{\sN})=\frac{\mbox{Var}(D)}{\mu^2}N^{2a-1},
    $$
    so that $b=1-2a>0$ and $c=\frac{\mbox{Var}(D)}{\mu^2}<\infty$.
\qed
\subsection{Coupling with a branching process with offspring $g^{\smallsup N}$}
\label{subsec-BP1}

We will construct the SPG in such a way that we simultaneously construct a
BP  with offspring distribution
$\{g_j^{\smallsup{N}}\}$ in (\ref{gnN}). This BP is of course purely imaginary.
The BP is coupled with the SPG such that it enables us to control their
difference.

As above, we will use the notation $Z_k^{\smallsup{1}}$ and
$Z_k^{\smallsup{2}}$ to denote the number of stubs attached
to nodes at distance $k-1$ from node 1, respectively, node 2,
but not part of an edge connected to a node at distance $k-2$.
For $k=1$, $Z_k^{\smallsup{i}}=D_i$.  We start with a
description of the coupling of the SPG with root 1, and
a BP with offspring distribution
$g^{\smallsup N}$ given in (\ref{gnN}). The first stages of the
generation of the SPG are drawn in Figure \ref{fig-1}. We will
explain the meaning of the labels 1, 2 and 3 below.

\begin{figure}[t]
\begin{center}
\setlength{\unitlength}{0.0004in}
{
\begin{picture}(6000,6500)(-1000,-5300)

\put(-2000,800) {\makebox(0,0)[lb]{SPG}}

\put(2500,800) {\makebox(0,0)[lb]{stubs with their labels}}

\path(0,0)(300,300)(600,0) \path(300,0)(300,300)

\path(1000,0)(1200,300)(1400,0)

\path(2000,0)(2400,300)(2800,0) \path(2200,0)(2400,300)(2600,0)

\path(3200,0)(3400,300)(3600,0) \path(3400,0)(3400,300)

\path(4000,0)(4400,300)(4800,0) \path(4200,0)(4400,300)(4600,0)
\path(4400,0)(4400,300)

\path(5200,0)(5400,300)(5600,0)

\path(6300,0)(6300,300)

\path(7000,0)(7200,300)(7400,0)

\path(8000,0)(8200,300)(8400,0)


\path(-2000,-1000)(-1700,-700)(-1400,-1000)
\path(-1700,-1000)(-1700,-700)

\path(0,-1000)(300,-700)(600,-1000) \path(300,-1000)(300,-700)
\put(-50,-1150){\makebox(0,0)[lb]{${\scriptscriptstyle 2}$}}
\put(250,-1150){\makebox(0,0)[lb]{${\scriptscriptstyle 2}$}}
\put(550,-1150){\makebox(0,0)[lb]{${\scriptscriptstyle 2}$}}

\path(1000,-1000)(1200,-700)(1400,-1000)

\path(2000,-1000)(2400,-700)(2800,-1000)
\path(2200,-1000)(2400,-700)(2600,-1000)

\path(3200,-1000)(3400,-700)(3600,-1000)
\path(3400,-1000)(3400,-700)

\path(4000,-1000)(4400,-700)(4800,-1000)
\path(4200,-1000)(4400,-700)(4600,-1000)
\path(4400,-1000)(4400,-700)

\path(5200,-1000)(5400,-700)(5600,-1000)

\path(6300,-1000)(6300,-700)

\path(7000,-1000)(7200,-700)(7400,-1000)

\path(8000,-1000)(8200,-700)(8400,-1000)

\path(-2300,-2300)(-1700,-1700)(-1400,-2000)
\path(-1700,-2000)(-1700,-1700) \path(-2300,-2300)(-2300,-2600)
\path(-2600,-2600)(-2300,-2300)(-2000,-2600)

\path(0,-2000)(300,-1700)(600,-2000) \path(300,-2000)(300,-1700)

\put(-50,-2150){\makebox(0,0)[lb]{${\scriptscriptstyle 3}$}}
\put(250,-2150){\makebox(0,0)[lb]{${\scriptscriptstyle 2}$}}
\put(550,-2150){\makebox(0,0)[lb]{${\scriptscriptstyle 2}$}}

\path(1000,-2000)(1200,-1700)(1400,-2000)

\path(2000,-2000)(2400,-1700)(2800,-2000)
\path(2200,-2000)(2400,-1700)(2600,-2000)

\put(1950,-2150){\makebox(0,0)[lb]{${\scriptscriptstyle 2}$}}
\put(2150,-2150){\makebox(0,0)[lb]{${\scriptscriptstyle 3}$}}
\put(2550,-2150){\makebox(0,0)[lb]{${\scriptscriptstyle 2}$}}
\put(2750,-2150){\makebox(0,0)[lb]{${\scriptscriptstyle 2}$}}

\path(3200,-2000)(3400,-1700)(3600,-2000)
\path(3400,-2000)(3400,-1700)

\path(4000,-2000)(4400,-1700)(4800,-2000)
\path(4200,-2000)(4400,-1700)(4600,-2000)
\path(4400,-2000)(4400,-1700)

\path(5200,-2000)(5400,-1700)(5600,-2000)

\path(6300,-2000)(6300,-1700)

\path(7000,-2000)(7200,-1700)(7400,-2000)

\path(8000,-2000)(8200,-1700)(8400,-2000)

\path(-2300,-3300)(-1700,-2700)(-1400,-3000)
\path(-1700,-3600)(-1700,-2700) \path(-1750,-3300)(-1650,-3300)
\path(-2300,-3300)(-2300,-3600)
\path(-2600,-3600)(-2300,-3300)(-2000,-3600)

\path(0,-3000)(300,-2700)(600,-3000) \path(300,-3000)(300,-2700)

\put(-50,-3150){\makebox(0,0)[lb]{${\scriptscriptstyle 3}$}}
\put(250,-3150){\makebox(0,0)[lb]{${\scriptscriptstyle 3}$}}
\put(550,-3150){\makebox(0,0)[lb]{${\scriptscriptstyle 2}$}}

\path(1000,-3000)(1200,-2700)(1400,-3000)

\path(2000,-3000)(2400,-2700)(2800,-3000)
\path(2200,-3000)(2400,-2700)(2600,-3000)

\put(1950,-3150){\makebox(0,0)[lb]{${\scriptscriptstyle 2}$}}
\put(2150,-3150){\makebox(0,0)[lb]{${\scriptscriptstyle 3}$}}
\put(2550,-3150){\makebox(0,0)[lb]{${\scriptscriptstyle 2}$}}
\put(2750,-3150){\makebox(0,0)[lb]{${\scriptscriptstyle 2}$}}

\path(3200,-3000)(3400,-2700)(3600,-3000)
\path(3400,-3000)(3400,-2700)

\path(4000,-3000)(4400,-2700)(4800,-3000)
\path(4200,-3000)(4400,-2700)(4600,-3000)
\path(4400,-3000)(4400,-2700)

\path(5200,-3000)(5400,-2700)(5600,-3000)

\path(6300,-3000)(6300,-2700)

\path(7000,-3000)(7200,-2700)(7400,-3000)

\path(8000,-3000)(8200,-2700)(8400,-3000)
\put(7950,-3150){\makebox(0,0)[lb]{${\scriptscriptstyle 3}$}}
\put(8350,-3150){\makebox(0,0)[lb]{${\scriptscriptstyle 2}$}}

\path(-2300,-4300)(-1700,-3700)(-1100,-4300)
\path(-1700,-4600)(-1700,-3700) \path(-1750,-4300)(-1650,-4300)
\path(-2300,-4300)(-2300,-4600)
\path(-1300,-4600)(-1100,-4300)(-900,-4600)
\path(-2600,-4600)(-2300,-4300)(-2000,-4600)


\path(0,-4000)(300,-3700)(600,-4000) \path(300,-4000)(300,-3700)

\put(-50,-4150){\makebox(0,0)[lb]{${\scriptscriptstyle 3}$}}
\put(250,-4150){\makebox(0,0)[lb]{${\scriptscriptstyle 3}$}}
\put(550,-4150){\makebox(0,0)[lb]{${\scriptscriptstyle 3}$}}

\path(1000,-4000)(1200,-3700)(1400,-4000)

\path(2000,-4000)(2400,-3700)(2800,-4000)
\path(2200,-4000)(2400,-3700)(2600,-4000)

\put(1950,-4150){\makebox(0,0)[lb]{${\scriptscriptstyle 2}$}}
\put(2150,-4150){\makebox(0,0)[lb]{${\scriptscriptstyle 3}$}}
\put(2550,-4150){\makebox(0,0)[lb]{${\scriptscriptstyle 2}$}}
\put(2750,-4150){\makebox(0,0)[lb]{${\scriptscriptstyle 2}$}}

\path(3200,-4000)(3400,-3700)(3600,-4000)
\path(3400,-4000)(3400,-3700)
\put(3150,-4150){\makebox(0,0)[lb]{${\scriptscriptstyle 3}$}}
\put(3350,-4150){\makebox(0,0)[lb]{${\scriptscriptstyle 2}$}}
\put(3550,-4150){\makebox(0,0)[lb]{${\scriptscriptstyle 2}$}}

\path(4000,-4000)(4400,-3700)(4800,-4000)
\path(4200,-4000)(4400,-3700)(4600,-4000)
\path(4400,-4000)(4400,-3700)

\path(5200,-4000)(5400,-3700)(5600,-4000)

\path(6300,-4000)(6300,-3700)

\path(7000,-4000)(7200,-3700)(7400,-4000)

\path(8000,-4000)(8200,-3700)(8400,-4000)
\put(7950,-4150){\makebox(0,0)[lb]{${\scriptscriptstyle 3}$}}
\put(8350,-4150){\makebox(0,0)[lb]{${\scriptscriptstyle 2}$}}

\path(-2300,-5300)(-1700,-4700)(-1100,-5300)
\path(-1700,-5300)(-1700,-4700) \path(-2300,-5300)(-2100,-5600)
\path(-1300,-5600)(-1100,-5300)(-900,-5600)
\path(-2300,-5300)(-2300,-5600) \path(-2300,-5300)(-1700,-5300)

\path(0,-5000)(300,-4700)(600,-5000) \path(300,-5000)(300,-4700)

\put(-50,-5150){\makebox(0,0)[lb]{${\scriptscriptstyle 3}$}}
\put(250,-5150){\makebox(0,0)[lb]{${\scriptscriptstyle 3}$}}
\put(550,-5150){\makebox(0,0)[lb]{${\scriptscriptstyle 3}$}}

\path(1000,-5000)(1200,-4700)(1400,-5000)

\path(2000,-5000)(2400,-4700)(2800,-5000)
\path(2200,-5000)(2400,-4700)(2600,-5000)

\put(1950,-5150){\makebox(0,0)[lb]{${\scriptscriptstyle 3}$}}
\put(2150,-5150){\makebox(0,0)[lb]{${\scriptscriptstyle 3}$}}
\put(2550,-5150){\makebox(0,0)[lb]{${\scriptscriptstyle 2}$}}
\put(2750,-5150){\makebox(0,0)[lb]{${\scriptscriptstyle 2}$}}

\path(3200,-5000)(3400,-4700)(3600,-5000)
\path(3400,-5000)(3400,-4700)
\put(3150,-5150){\makebox(0,0)[lb]{${\scriptscriptstyle 3}$}}
\put(3350,-5150){\makebox(0,0)[lb]{${\scriptscriptstyle 2}$}}
\put(3550,-5150){\makebox(0,0)[lb]{${\scriptscriptstyle 2}$}}

\path(4000,-5000)(4400,-4700)(4800,-5000)
\path(4200,-5000)(4400,-4700)(4600,-5000)
\path(4400,-5000)(4400,-4700)

\path(5200,-5000)(5400,-4700)(5600,-5000)

\path(6300,-5000)(6300,-4700)

\path(7000,-5000)(7200,-4700)(7400,-5000)

\path(8000,-5000)(8200,-4700)(8400,-5000)
\put(7950,-5150){\makebox(0,0)[lb]{${\scriptscriptstyle 3}$}}
\put(8350,-5150){\makebox(0,0)[lb]{${\scriptscriptstyle 3}$}}

\end{picture}
}
\end{center}
\caption{Schematic drawing of the growth of the SPG from the node
1 with $N=9$ and the updating of the labels. The stubs without
labels have label 1. The first line shows the $N$ different
degrees. The growth process
starts by choosing the first stub of node 1 whose stubs are
labeled by 2 as illustrated in the second line, while all the
other stubs maintain the label 1. Next, we uniformly choose a stub
with label 1 or 2. In the example in line 3, this is the second
stub from node 3, whose stubs are labeled by 2 except for the second
stub which is labeled 3. The left hand side column visualizes growth of
the SPG by the attachment of stub 2 of node 3 to the first stub of
node 1. Once an edge is established the paired stubs are labeled
3. In the next step, the next stub of node one is again matched to a
uniform stub out of those
with label 1 or 2. In the example in line 4, it is the first stub
of the last node that will be attached to the second stub of node
1, the next in sequence to be paired. The last line exhibits the
result of creating a cycle when the first stub of node 3 is chosen
to be attached to the last stub of node 9 (the last node).
This process is continued until there are no more stubs with labels 1 or 2.
In this example, we have $Z_1^{\smallsup{1}}=3$ and $Z_2^{\smallsup{1}}=6$.}
\label{fig-1}
\end{figure}
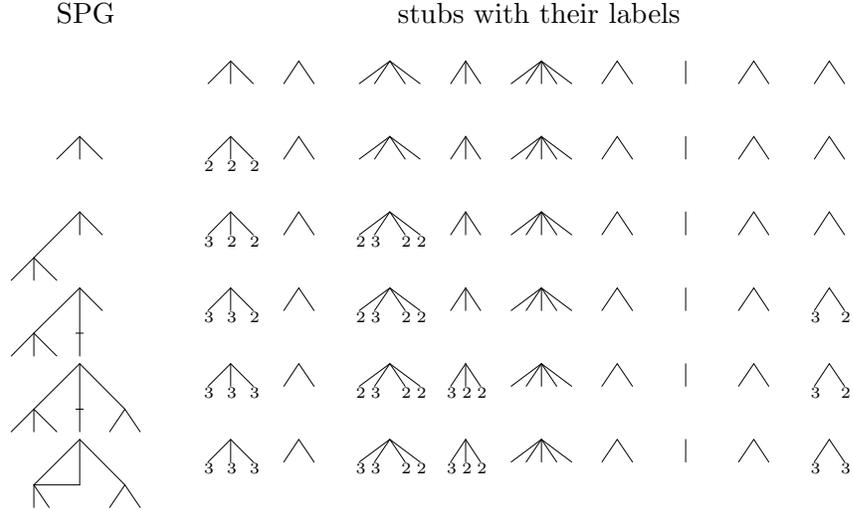

\begin{figure}[t]
\begin{center}
\setlength{\unitlength}{0.0007in}
{
\begin{picture}(4500,1700)(0,-500)
\path(0,0)(500,500)(1000,0) \path(300,-300)(0,0)
\path(-100,-300)(0,0) \path(100,-300)(0,0)
\dottedline{50}(-300,-300)(0,0) \dottedline{50}(800,-300)(1000,0)
\path(1200,-300)(1000,0)

{\small \put(-700,-100) {\makebox(0,0)[lb]{$D_i\!=\!5$}}
\put(300,-100) {\makebox(0,0)[lb]{$D_j\!=\!3$}} }

\path(2000,0)(2500,500)(3000,0)
\path(1900,-300)(2000,0) \path(2100,-300)(2000,0)
\path(2300,-300)(2000,0) \thicklines \path(2000,0)(3000,0)
\thinlines \path(3200,-300)(3000,0)

\put(2300,800) {\makebox(0,0)[lb]{SPG}}

\path(4000,0)(4500,500)(5000,0) \path(3500,-500)(4000,0)
\path(3900,-300)(4000,0) \path(4100,-300)(4000,0)
\path(4300,-300)(4000,0) \path(3400, -800)(3500,-500)
\path(3600,-800)(3500,-500) \put(3500, -500){\circle{150}}
\path(5200,-300)(5000,0) \path(4800,-300)(5000,0)

\put(4400,800) {\makebox(0,0)[lb]{BP}}

\end{picture}
}
\end{center}
\caption{Example of the coupling when a cycle occurs. Edges have
twice the length of stubs. In the SPG the two dotted stubs in the
left picture are to be connected. The middle picture gives the
result of creating the cycle in the SPG where the bold line is the
edge creating the cycle. The third figure draws the BP where the
cycle is removed and the degree of the circled node is 3.}
\label{fig-3}
\end{figure}
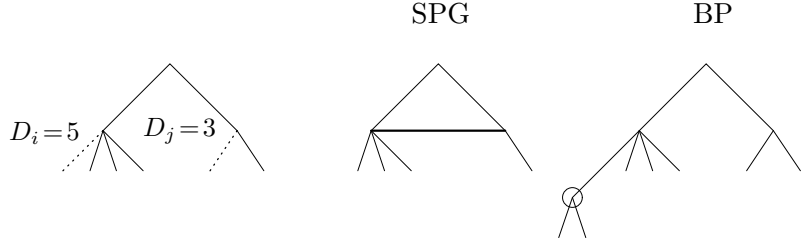

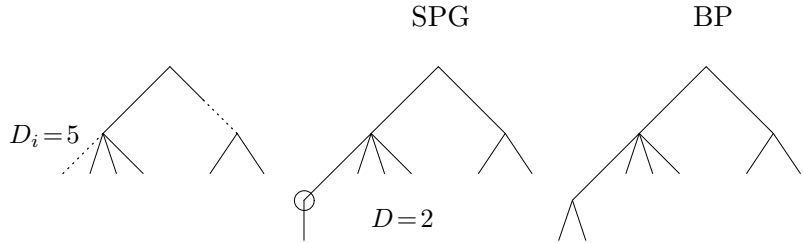
\begin{figure}[t]
\begin{center}
\setlength{\unitlength}{0.0007in}
{
\begin{picture}(4500,1700)(0,-500)
\path(0,0)(500,500) \path(500,500)(750,250)
\dottedline{50}(750,250)(1000,0) \path(300,-300)(0,0)
\path(-100,-300)(0,0) \path(100,-300)(0,0)
\dottedline{50}(-300,-300)(0,0)
\path(800,-300)(1000,0)
\path(1200,-300)(1000,0)

{\small \put(-700,-100) {\makebox(0,0)[lb]{$D_i\!=\!5$}} }

\path(2000,0)(2500,500)(3000,0)
\path(1500,-500)(2000,0) \path(1900,-300)(2000,0)
\path(2100,-300)(2000,0) \path(2300,-300)(2000,0)
\path(2800,-300)(3000,0)
\path(3200,-300)(3000,0)
\put(1500,-500){\circle{150}}
\path(1500,-800)(1500,-500)

\put(2300,800) {\makebox(0,0)[lb]{SPG}} {\small
\put(2000,-700){\makebox(0,0)[lb]{$D\!=\!2$}} }

\path(4000,0)(4500,500)(5000,0)
\path(3500,-500)(4000,0)
\path(3900,-300)(4000,0)
\path(4100,-300)(4000,0)
\path(4300,-300)(4000,0)
\path(3400, -800)(3500,-500)
\path(3600,-800)(3500,-500)
\path(5200,-300)(5000,0)
\path(4800,-300)(5000,0)

\put(4400,800) {\makebox(0,0)[lb]{BP}}

\end{picture}
}
\end{center}
\caption{An example of the coupling where we need to perform a
redraw. In the draw from
$g^{\smallsup{N}}$, we draw the dotted stub in the SPG with degree
3. In the BP, we keep this degree, while in the SPG we draw again
from the conditional distribution given that we do not draw a
stub with label 3. In this example, this redraw gives the value
$D=2$.} \label{fig-2}
\end{figure}

We draw repeatedly and independently from the distribution
$\{g_j^{\smallsup{N}}\}$. This is done conditionally given
$D_1,D_2,\ldots,D_{\sN}$, so that we draw from
the {\it random} distribution (\ref{gnN}).
After each draw we will update the realization of the SPG and the
BP, and classify the stubs according to three categories, which
will be labelled 1, 2 and 3. These labels will be updated as the growth
of the SPG proceeds. The labels have the following meaning:
\begin{enumerate}
\item[1.]
Stubs with label 1 are stubs belonging to a node that is not yet
attached to the SPG.
\item[2.]
Stubs with label 2 are attached to the SPG
(because the corresponding node has been chosen), but not yet
paired with another stub. These are called `free stubs'.
\item[3.]
Stubs with label 3 in the SPG are paired with another stub to form
an edge in the SPG.
\end{enumerate}

The growth process as depicted in Figure \ref{fig-1} starts by
giving all stubs label 1. Then, because we construct the SPG starting
from node $1$, we relabel the $D_1$ stubs of node $1$ with the label $2$.
We note that $Z_1^{\smallsup{1}}$ is equal to the number of stubs connected to
node 1, and thus $Z_1^{\smallsup{1}}=D_1$. We next identify $Z_j^{\smallsup{1}}$ for $j>1$.
$Z_j^{\smallsup{1}}$ is obtained by sequentially growing the SPG from the free stubs
in generation $Z_{j-1}^{\smallsup{1}}$. When all free stubs in generation $j-1$
have chosen their connecting stub, $Z_j^{\smallsup{1}}$ is equal to the number of stubs
labelled 2 (i.e., free stubs) attached to the SPG. Note that not necessarily
each stub of $Z_{j-1}^{\smallsup{1}}$ contributes to stubs of $Z_{j}^{\smallsup{1}}$, because a cycle may
`swallow' two free stubs in generation $j-1$. This is the case precisely
when a stub with label $2$ is chosen.

For the BP, we start with $\hat Z_1^{\smallsup{1}}=D_1$,
and grow from the free stubs available in the BP tree by sequentially
growing from the stubs (alike for the SPG). For the coupling, as long as
there are free stubs in {\it both} the BP and the SPG in a given generation,
we couple the BP and SPG in the following way.
At each step we will take an independent draw from all stubs,
according to the distribution (\ref{gnN}). Since the stubs are
specified by their label (1, 2 or 3), we can now present the
construction rules for the BP and the SPG.
\begin{enumerate}
\item[1.] If the chosen stub has label 1, then in both the BP and the SPG we will
connect the present stub to the chosen stub to form an edge and
attach the remaining stubs of the chosen node as children. We update the labels
as follows. The present and chosen stub melt together to form an edge
and both are assigned label 3. All `brother' stubs (except for the
chosen stub) belonging to the same node of the chosen stub receive
label 2.

\item[2.] In this case we choose a stub with label 2, which is
already connected to the SPG. For the BP, the chosen stub is
simply connected to the stub which is grown, and the number
of free stubs is the number of `brother stubs' of the chosen stub. For the SPG, a
self-loop is created when the chosen stub and present stub
are `brother' stubs which belong to the same
node. When they are not `brother' stubs, then a cycle is formed. Neither
a self-loop nor a cycle changes the distances in the SPG. Note that
for the SPG {\it two} free stubs are used, while for the BP only
{\it one} stub is used. This is illustrated in Figure \ref{fig-3}.

The updating of the labels solely consists of changing the label
of the present and the chosen stub from 2 to 3.

\item[3.]
A stub with label 3 is chosen. This case is illustrated in Figure
\ref{fig-2}. This possibility of choosing an already matched stub
with label 3 must be included for the BP which relies on the property
that all subsequent iterations in the process are i.i.d. Note that
this includes the case where we draw the present stub, which of
course is impossible for the SPG.

The rule now for the BP is that the corresponding node with the
prescribed number of stubs is simply attached.  Since for the SPG,
we sample {\it without} replacement, we have to resample from
distribution (\ref{gnN}), until we draw a stub with label 1 or 2.
This procedure is referred to as {\em a redraw}. Since we
sample uniformly from all stubs, the conditional sampling until we
hit a stub with label 1 or 2 is also uniform out of the set of
all stubs with labels 1 and 2, so that it has the correct
distribution. Obviously there are two cases: either we draw a stub
with label 1 or one with label 2. When we draw a stub with label
1 in the SPG then we update as under rule 1 above, while
when we draw a stub having label 2 in the SPG, we update as
under rule 2 above.
\end{enumerate}

Clearly, the redraws and the cycles cause possible
differences between the BP and the SPG: the degrees of
the chosen node are possibly different. We will need to
show that the above difference only leads to an error term.

The above process stops in the $j^{\rm th}$ generation when
there are no more free stubs in generation $j-1$ for either the BP or
for the SPG. When there are no more free stubs for the SPG, we complete the
$j^{\rm th}$ generation for the BP by drawing from distribution (\ref{gnN})
for all the remaining free stubs. The labels of the stubs remain unchanged.
When there are no more free stubs for the BP, we complete the
$j^{\rm th}$ generation for the SPG by drawing from distribution (\ref{gnN})
iteratively until we draw a stub with label 1 or 2. This is done
for all the remaining free stubs in the $j^{\rm th}$ generation
of the SPG. The labels are updated as under 1 and 2 above.

We continue the above process of drawing stubs until there are no
more stubs having label 1 or 2, so that all stubs have label 3.
Then, the construction is finalized, and we have generated the SPG
as seen from node 1. We have thus obtained the structure of
the SPG, and know how many nodes there are at a
given distance from node 1.

The above construction will be performed similarly
from node $2$. This construction is close to being independent as long as the
SPG's from the roots 1 and 2 do not share any nodes. More precisely,
the corresponding BP's are independent.
Thus, we have now constructed the SPG's and BP's
from both node 1 and node 2.

\subsection{Coupling with a BP with offspring distribution $\{g_j^{\smallsup{N}}\}$}
\label{sec-bdscoupling}
In the previous section, we have obtained a coupling of the
SPG and the BP with offspring distribution $\{g_j^{\smallsup{N}}\}$.
In this and the next section, we will summarize bounds on the couplings that
we need for the proof of Theorems \ref{thm-tau>3} and \ref{thm-ll}. These results
will be repeated in the appendix together with a full proof.
We start with the coupling of the number of stubs
$Z_j^{\smallsup{1}}$ in the SPG and the number of children $\hat Z_j^{\smallsup{1}}$
in the $j^{\rm th}$ generation of the BP with offspring distribution $\{g_j^{\smallsup{N}}\}$.
\begin{prop}[Coupling SPG with the BP with random offspring distribution]
\label{prop-bdscoupling1}
    There exist $\eta,\beta>0$, $\alpha>\frac 12+\eta$ and a constant $C$,
    such that for all $j\leq (\frac 12 +\eta) \log_\nu N$,
    \begin{equation}
    \prob\Big((1-N^{-\alpha} \nu^{j}) \hat Z_j^{\smallsup{1}}\leq Z_j^{\smallsup{1}}
    \leq (1+N^{-\alpha} \nu^{j}) \hat Z_j^{\smallsup{1}}\Big)\geq 1- C j N^{-\beta}.
    \end{equation}
\end{prop}
\subsection{Coupling with a BP with offspring distribution $\{g_j\}$}
\label{subsec-BP3}
We next describe the coupling with the
BP with offspring distribution $\{g_j\}$ and their bounds.
A classical coupling argument is used (see e.g. \cite{thor}). Let
$X^{\smallsup{N}}$ have law $\{g_j^{\smallsup{N}}\}$ and $X$ have
law $\{g_j\}$. We define $Y^{\smallsup{N}}$ by
    \begin{equation}
    \prob(Y^{\smallsup{N}}=n)=\min(g_n^{\smallsup{N}}, g_n),\qquad
    \prob(Y^{\smallsup{N}}=\infty)=1-\sum_{n=0}^{\infty}\min(g_n^{\smallsup{N}}, g_n)
    =\frac12 \sum_{n=0}^{\infty} |g_n^{\smallsup{N}}-g_n|.
    \end{equation}
Let $\hat X^{\smallsup{N}}=Y^{\smallsup{N}}$ when
$Y^{\smallsup{N}}<\infty$, and  $\prob(X^{\smallsup{N}}=n,
Y^{\smallsup{N}}=\infty)=g_n^{\smallsup{N}}-\min(g_n^{\smallsup{N}},
g_n)$, whereas $\hat X=X$ when
$Y^{\smallsup{N}}<\infty$, and  $\prob(X=n,
Y^{\smallsup{N}}=\infty)=g_n-\min(g_n^{\smallsup{N}}, g_n)$. Then
$\hat X^{\smallsup{N}}$ has law $g^{\smallsup{N}}$, and $\hat X$
has law $g$. Moreover, with large probability, $\hat
X^{\smallsup{N}}=\hat X$ due to Proposition \ref{prop-TV1} below.

This coupling argument is applied to each node in the BP
$\{{\hat Z}_i^{\smallsup{1}}\}_{i\geq 0}$ and $\{{\hat
Z}_i^{\smallsup{2}}\}_{i\geq 0}$. The BP's with
offspring distribution $\{g_j\}$ will be denoted by $\{{\cal
Z}_i^{\smallsup{1}}\}_{i\geq 0}$ and $\{{\cal
Z}_i^{\smallsup{2}}\}_{i\geq 0}$. We can interpret this coupling as
follows. Each node has an i.i.d.\ indicator variable which equals one with
probability
    \eq
    \label{pNdef}
    p_{\sN}=\frac12 \sum_{n=0}^{\infty} |g_n^{\smallsup{N}}-g_n|.
    \en
When at a certain node this indicator variable is 0, then the
offspring in $\{{\hat Z}_i^{\smallsup{1}}\}_{i\geq 0}$ or $\{{\hat
Z}_i^{\smallsup{2}}\}_{i\geq 0}$ equals the one in $\{{\cal
Z}_i^{\smallsup{1}}\}_{i\geq 0}$ or $\{{\cal
Z}_i^{\smallsup{2}}\}_{i\geq 0}$, and the node is successfully
coupled. When the indicator is 1, then an error has occurred, and
the coupling is not successful. In this case, the laws of the
offspring of $\{{\hat Z}_i^{\smallsup{1}}\}_{i\geq 0}$ or $\{{\hat
Z}_i^{\smallsup{2}}\}_{i\geq 0}$ is different from the one in
$\{{\cal Z}_i^{\smallsup{1}}\}_{i\geq 0}$ or $\{{\cal
Z}_i^{\smallsup{2}}\}_{i\geq 0}$, and we record an error.
Below we will use the notation $\prob_{\sN}$ to denote the conditional expectation
given $D_1,D_2,\ldots,D_{\sN}$ and $\expec_{\sN}$ to denotes the expectation with respect to
the probability measure $\prob_{\sN}$. Finally, we write
\begin{equation}
\label{nu-N}
\nu_{\sN}=\sum_{n=0}^\infty n g_n^{\smallsup{N}}.
\end{equation}

In the following proposition, we prove that at any fixed time, we
can couple the SPG to the delayed BP with law $\{g_j\}$:

\begin{prop}[Coupling at fixed time]
\label{prop-caft}
For any $m\in \mathbb{N}$ fixed,
there exist {\it independent} branching processes ${\cal Z}^{\smallsup{1}},
{\cal Z}^{\smallsup{2}}$, such that
    \eq
    \lim_{N\rightarrow \infty}
    \prob(Z_m^{\smallsup{i}}={\cal Z}_m^{\smallsup{i}})=
    1.
    \en
\end{prop}

In the course of the proof we will also rely on the following more
technical claims:

\begin{prop}[Convergence in total variation distance]
    \label{prop-TV1}
    There exist $\alpha_2,\beta_2>0$ such that
    \begin{equation}
    \prob\big(\sum_{n=0}^{\infty} (n+1)|g_n^{\smallsup{N}}-g_n| \geq N^{-\alpha_2}\big)\leq N^{-\beta_2}.
    \end{equation}
Consequently,
\begin{equation}
\label{propnu}
\prob(|\nu_{\sN}-\nu|>N^{-\alpha_2})\leq N^{-\beta_2},
\end{equation}
and
\begin{equation}
\label{proppN}
\prob(p_{\sN}>N^{-\alpha_2})\leq N^{-\beta_2}.
\end{equation}
\end{prop}

\begin{corr}[Coupling of sums]
There exist $\eps,\beta,\eta>0$ such that for all $j\le (1+2\eta)\log_{\nu} N$,
as $N\to \infty$,
    \label{corr-coup1}
   \begin{equation}
   \label{vgl-coup1}
     \prob\Big(\frac 1N\Big|\sum_{i=1}^{j}{\cal
    Z}_{\lceil i/2\rceil}^{\smallsup{1}}{\cal Z}_{\lfloor i/2\rfloor}^{\smallsup{2}}-
    \sum_{i=1}^{j}{\hat
Z}_{\lceil i/2\rceil}^{\smallsup{1}}{\hat Z}_{\lfloor i/2\rfloor}^{\smallsup{2}}\Big|>N^{-\eps}\Big)
    =O(N^{-\beta}).
   \end{equation}
\end{corr}

\section{Proof of Theorem \ref{thm-tau>3} and \ref{thm-ll}}
\label{sec-pftau>3} The proof consists of four steps.
\begin{enumerate}
\item[1.] We first express the survival probability $\prob(H_{\sN}>j)$
in the number of stubs $\{Z_i^{\smallsup{k}}\}, k=1,2$, of the
SPG's. For $j\leq (1+2\eta)\log_{\nu}N$,
where $\eta$ is specified in Proposition
\ref{prop-bdscoupling1}, we will show that
    \begin{equation}
    \label{compeven} \prob(H_{\sN}>j)= \expec\left[
    \exp\left\{
    \frac{-\sum_{i=2}^{j+1} Z_{\lceil i/2 \rceil}^{\smallsup{1}}
    Z_{\lfloor i/2\rfloor}^{\smallsup{2}}} {L_{\sN}} \right\}+RM_{\sN}(j)\right],
    \end{equation}
with
    $$
    RM_{\sN}(j)=
    O\left(\sum_{i=2}^{j+1}
    \frac{Z^{\smallsup{1}}_{\lceil i/2\rceil}Z_{\lfloor
    i/2\rfloor}^{\smallsup{2}}
        \sum_{k=1}^{\lceil i/2\rceil}(Z_k^{\smallsup{1} }+Z_k^{\smallsup{2}})}{L_{\sN}^2}
    \right).
    $$

\item[2.] We use Proposition \ref{prop-bdscoupling1} to show that in (\ref{compeven}) we can replace
$\{Z_k^{\smallsup{i}}\}, i=1,2$  by the BP $\{{\hat
Z}_k^{\smallsup{i}}\}, i=1,2$. The error term $\expec[|RM_{\sN}(j)|]$ and the error
involved in replacing the SPG by the BP is bounded by a constant times $N^{-\beta}$,
for some $\beta >0$, uniformly in $j\leq
(1+2\eta)\log_{\nu} N$.

\item[3.]
In this step we show that there exists $\beta>0$ such that for all
$j\leq (1+2\eta)\log_{\nu}N$, as $N\to \infty$,
    \begin{equation}
    \label{hhoednaarexp} \prob(H_{\sN}>j)= \expec\left[
    \exp\left\{
    \frac{-\sum_{i=2}^{j+1} {\cal Z}_{\lceil i/2
    \rceil}^{\smallsup{1}} {\cal Z}_{\lfloor
    i/2\rfloor}^{\smallsup{2}} } {\mu N} \right\}\right]
    +O(N^{-\beta}),
    \end{equation}
where ${\cal Z}_k^{\smallsup{i}},\, i=1,2,$ denotes the delayed BP with
offspring distribution (\ref{outgoing degree}).

\item[4.] We complete the proof of Theorem \ref{thm-tau>3} and
\ref{thm-ll}, using {\it step 3}, and the almost sure limit
in (\ref{martin}) applied to ${\cal Z}_n^{\smallsup{1}}$ and ${\cal Z}_n^{\smallsup{2}}$.
We finally use the speed of convergence of the above martingale
limit result to obtain (\ref{uitspraakrep}).
\end{enumerate}

\noindent{\bf Step 1: A formula for $\prob(H_{\sN}>j)$.}
The following lemma expresses $\prob(H_{\sN}>j)$
in terms of $\Q_Z^{\smallsup{k,l}}$, the conditional probabilities
given $\{Z_s^{\smallsup{1}}\}_{s=1}^k$ and $\{Z_s^{\smallsup{2}}\}_{s=1}^l$.
For $l=0$, we only condition on $\{Z_s^{\smallsup{1}}\}_{s=1}^k$.
\begin{lemma}
\label{lem-conditioneel}
For $j\geq 1$,
\begin{equation}
\label{multiplication rule}
\prob(H_{\sN}>j)=\expec\Big[\prod_{i=2}^{j+1}
    \Q_Z^{\smallsup{\lceil i/2\rceil,\lfloor i/2\rfloor}}
    (H_{\sN}>i-1|H_{\sN}>i-2)\Big].
\end{equation}
\end{lemma}

\proof
We first compute that
     \[
     \prob(H_{\sN}>j)
     =\expec\big[\Q_Z^{\smallsup{1,1}}(H_{\sN}>j)\big]
     =\expec\big[\Q_Z^{\smallsup{1,1}}(H_{\sN}>1)\Q_Z^{\smallsup{1,1}}(H_{\sN}>j|H_{\sN}>1)\big].
     \]
Continuing this further, and writing $\expec_Z^{\smallsup{k,l}}$ for the expectation with respect to $\Q_Z^{\smallsup{k,l}}$, \begin{eqnarray*}
     \Q_Z^{\smallsup{1,1}}(H_{\sN}>j|H_{\sN}>1)
     &=&\expec_Z^{\smallsup{1,1}}\big[\Q_Z^{\smallsup{2,1}}(H_{\sN}>j|H_{\sN}>1)\big]\\
     &=&\expec_Z^{\smallsup{1,1}}\big[\Q_Z^{\smallsup{2,1}}(H_{\sN}>2|H_{\sN}>1)
     \Q_Z^{\smallsup{2,1}}(H_{\sN}>j|H_{\sN}>2)\big].
\end{eqnarray*}
Therefore,
\begin{eqnarray*}
     \prob(H_{\sN}>j)
     &=&\expec\big[\Q_Z^{\smallsup{1,1}}(H_{\sN}>1)\expec_Z^{\smallsup{1,1}}
     \big[\Q_Z^{\smallsup{2,1}}(H_{\sN}>2|H_{\sN}>1)
     \Q_Z^{\smallsup{2,1}}(H_{\sN}>j|H_{\sN}>2)\big]\big]\\
     &=&\expec\big[\expec_Z^{\smallsup{1,1}}
     \big[\Q_Z^{\smallsup{1,1}}(H_{\sN}>1)\Q_Z^{\smallsup{2,1}}(H_{\sN}>2|H_{\sN}>1)
     \Q_Z^{\smallsup{2,1}}(H_{\sN}>j|H_{\sN}>2)\big]\big]\\
     &=&\expec\big[\Q_Z^{\smallsup{1,1}}(H_{\sN}>1)\Q_Z^{\smallsup{2,1}}(H_{\sN}>2|H_{\sN}>1)
     \Q_Z^{\smallsup{2,1}}(H_{\sN}>j|H_{\sN}>2)\big],
\end{eqnarray*}
where, in the second equality, we use that $\Q_Z^{\smallsup{1,1}}(H_{\sN}>1)$
is measurable with respect to the $\sigma$-algebra generated by
$Z_1^{\smallsup{1,N}}$.
This proves the claim for $j=2$.

More generally, we obtain that for $k,l$ such that $k+l\leq j-1$, \begin{eqnarray*}
    && \Q_Z^{\smallsup{k,l}}(H_{\sN}>j|H_{\sN}>k+l-1)
     =\expec_Z^{\smallsup{k,l}}\big[\Q_Z^{\smallsup{k,l+1}}(H_{\sN}>j|H_{\sN}>k+l-1)\big]\\
     &&\qquad=\expec_Z^{\smallsup{k,l}}\big[\Q_Z^{\smallsup{k,l+1}}(H_{\sN}>k+l|H_{\sN}>k+l-1)
     \Q_Z^{\smallsup{k,l+1}}(H_{\sN}>j|H_{\sN}>k+l)\big],
\end{eqnarray*}
and, similarly,
\[
\Q_Z^{\smallsup{k,l}}(H_{\sN}>j|H_{\sN}>k+l-1)=
\expec_Z^{\smallsup{k,l}}\big[\Q_Z^{\smallsup{k+1,l}}(H_{\sN}>k+l|H_{\sN}>k+l-1)
     \Q_Z^{\smallsup{k+1,l}}(H_{\sN}>j|H_{\sN}>k+l)\big].
\]
In the above formulas, we can choose to increase $k$ or $l$
by one depending on $\{Z_s^{\smallsup{1,N}}\}_{s=1}^k$ and
$\{Z_s^{\smallsup{2,N}}\}_{s=1}^l$.  We will iterate the above recursions,
until $k+l=j-1$, when the last term becomes 1. This yields
that
\begin{equation}
\label{multiplication rule2} \prob(H_{\sN}>j)=\expec\Big[\prod_{i=1}^{j}
     \Q_Z^{\smallsup{\lfloor i/2\rfloor+1,\lceil i/2\rceil}}
     (H_{\sN}>i|H_{\sN}>i-1)\Big].
\end{equation}
Renumbering gives the final result.
\qed


We will next prove (\ref{compeven}). In order to do so,
we start by proving upper and lower bounds on the probabilities of not connecting two
sets of stubs to each other. For this, suppose we have two disjoint sets of stubs
$A$ with $|A|=n$ and $B$ with $|B|=m$ out of a total of $L$ stubs. We match stubs
at random, in such a way that two stubs form one edge, as in the construction
of the SPG. In particular, loops are possible.

\medskip

Let $p(n,m,L)$
denote the probability  that none of the $n$ stubs in $A$ attaches to one of the $m$ stubs in $B$.
Then, by conditioning on whether we choose a stub in $A$ or not, we obtain the recursion
    \begin{equation}
    \label{recursion}
    p(n,m,L)=\frac{n-1}{L-1}p(n-2,m,L-2)+\left(1-\frac{m+n-1}{L-1}\right)p(n-1,m,L-2)
    \end{equation}
Since $p(n-2,m,L-2)\geq p(n-1,m,L-2)$, because we have to match one additional stub, we obtain
    \begin{equation}
    \label{lowerbound}
    p(n,m,L)\geq \left(1-\frac{m}{L-1}\right)p(n-1,m,L-2)\geq\prod_{i=0}^{n-1}\left(1-\frac{m}{L-2i-1}\right).
    \end{equation}
On the other hand, we can rewrite (\ref{recursion}) as
    \begin{equation}
    \label{rewrite}
    p(n,m,L)=\left(1-\frac{m}{L-1}\right)p(n-1,m,L-2)+\frac{n-1}{L-1}
    \left(p(n-2,m,L-2)-p(n-1,m,L-2)\right).
    \end{equation}
We claim that
    \begin{equation}
    \label{ineq}
    p(n-2,m,L-2)-p(n-1,m,L-2)
    =
    \frac{m}{L-3}p(n-2,m-1,L-2)\leq \frac{m}{L-3}.
    \end{equation}
Indeed, the difference $p(n-2,m,L-2)-p(n-1,m,L-2)$ is equal to the
probability of the event that the first $n-2$ stubs do not connect
to $B$, while the last one does. By exchangeability of the stubs, this probability
equals the probability that the first stub is attached to a stub in $B$,
and the remaining $n-2$ stubs are not. This latter probability is equal
to $\frac{m}{L-3}p(n-2,m-1,L-2)$.

The equations (\ref{rewrite}) and (\ref{ineq}) yield
    $$
    p(n,m,L)\leq \left(1-\frac{m}{L-1}\right)p(n-1,m,L-2)+\frac{n-1}{(L-1)}
    \frac{m}{(L-3)}.
    $$
Iteration gives the upper bound
    \begin{equation}
    \label{upperbound}
    p(n,m,L)\leq
    \left[ \prod_{i=0}^{n-1}\left(1-\frac{m}{L-2i-1}\right)\right]+\frac{n^2m}{(L-2n)^2} .
    \end{equation}

\medskip

Since the event $\{H_{\sN}>1\}$ holds if and
only if no stubs of root $1$ attaches to one of those of root $2$,
we obtain, using (\ref{lowerbound}) and (\ref{upperbound}), that
    \begin{equation}
    \label{low-bound}
    \prod_{i=0}^{Z_1^{\smallsup{1}}-1} \left(1- \frac
    {Z_1^{\smallsup{2}}} {L_{\sN}-2i-1} \right)
    \leq \Q_Z^{\smallsup{1,1}}(H_{\sN}>1)\leq \left[
    \prod_{i=0}^{Z_1^{\smallsup{1}}-1} \left(1- \frac
    {Z_1^{\smallsup{2}}} {L_{\sN}-2i-1} \right) \right]
    +\frac{(Z_1^{\smallsup{1}})^2Z_1^{\smallsup{2}}}{(L_{\sN}-2Z_1^{\smallsup{1}})^2}.
    \end{equation}
Similarly,
    \begin{equation}
    \Q_Z^{\smallsup{2,1}}(H_{\sN}>2|H_{\sN}>1)\geq \prod_{i=0}^{Z_1^{\smallsup{2}}-1} \left(1-
    \frac {Z_2^{\smallsup{1}}} {L_{\sN}-2Z_1^{\smallsup{1}}-2i-1} \right),
    \end{equation}
with a matching upper bound with an error term bounded by $\frac{(Z_1^{\smallsup{2}})^2Z_2^{\smallsup{1}}}{(L_{\sN}-2Z_1^{\smallsup{1}}-2Z_1^{\smallsup{2}})^2}.$

We use that, for natural numbers $n,m,M$ with $M+n+m=o(L)$,
    \begin{equation}
    \label{benadering} \prod_{i=0}^{n-1} \left(1- \frac {m} {L-M-2i-1}
    \right)=
    e^{-\frac{nm}{L}}
    +O\left(\frac{nm(M+n+m)}{L^2}\right), \quad L\to\infty.
    \end{equation}

\noindent
Using (\ref{benadering}), the bounds in (\ref{low-bound}) yield
    $$
    \Q_Z^{\smallsup{1,1}}(H_{\sN}>1)= \exp
    \left\{-\frac{Z^{\smallsup{1}}_{1}Z_{1}^{\smallsup{2}}}{L_{\sN}}\right\}
    +O\left(\frac{Z^{\smallsup{1}}_{1}Z_{1}^{\smallsup{2}}
        (Z_1^{\smallsup{1}}+Z_1^{\smallsup{2}})}{L_{\sN}^2}\right).
    $$
\noindent
Similarly, we can conclude that, as long as $\sum_{k=1}^{\lceil i/2 \rceil}(Z_k^{\smallsup{1}}+Z_k^{\smallsup{2}})=o(L_{\sN})$,
we have
    \begin{eqnarray}
    \label{hoofd-eerste}
    &&\Q_Z^{\smallsup{\lceil
    i/2\rceil,\lfloor i/2\rfloor}}(H_{\sN}>i-1|H_{\sN}>i-2)\nonumber\\
    &&\quad = \exp \left\{-\frac{Z^{\smallsup{1}}_{\lceil
    i/2\rceil}Z_{\lfloor i/2\rfloor}^{\smallsup{2}}}{L_{\sN}}\right\}
    +O\left(\frac{Z^{\smallsup{1}}_{\lceil i/2\rceil}Z_{\lfloor
    i/2\rfloor}^{\smallsup{2}}
        (\sum_{k=1}^{\lceil i/2 \rceil}(Z_k^{\smallsup{1}}+Z_k^{\smallsup{2}}))}{L_{\sN}^2}\right).
    \end{eqnarray}
From (\ref{multiplication rule}) and taking expectations, the main
term in (\ref{compeven}) is evident. For the error term, we obtain that,
as long as $\sum_{k=1}^{\lceil i/2 \rceil}(Z_k^{\smallsup{1}}+Z_k^{\smallsup{2}})=o(L_{\sN})$,
    $$
    RM_{\sN}(j)=
    \sum_{i=2}^{j+1}
    \frac{Z^{\smallsup{1}}_{\lceil i/2\rceil}Z_{\lfloor
    i/2\rfloor}^{\smallsup{2}}
    \sum_{k=1}^{\lceil i/2\rceil}(Z_k^{\smallsup{1} }+Z_k^{\smallsup{2}})}{L_{\sN}^2},
    $$
and we will show at the end of step 2 that for all $j<(1+2\eta)\log N$,
we have $\sum_{k=1}^{\lceil j/2 \rceil}(Z_k^{\smallsup{1}}+Z_k^{\smallsup{2}})=o(L_{\sN})$
and that there exists a $\beta>0$ such that
    \begin{equation}
    \label{ordefout}
    \expec[RM_{\sN}(j)]=O(N^{-\beta}).
    \end{equation}
\bigskip

\noindent{\bf Step 2:  Coupling of SPG  to the BP with offspring $\{g_j^{\smallsup{N}}\}$.}
We start by showing that for some $\beta>0$
and uniformly in $j\le (1+2\eta)\log_{\nu}N$, the main term in (\ref{compeven}) satisfies
    \begin{equation}
    \label{naarzhoed}
    \expec\left[
    \exp\left\{
    \frac{-\sum_{i=2}^{j+1} Z_{\lceil i/2 \rceil}^{\smallsup{1}}
    Z_{\lfloor i/2\rfloor}^{\smallsup{2}} } {L_{\sN}} \right\}\right]
    =\expec\left[
    \exp\left\{
    \frac{-\sum_{i=2}^{j+1} {\hat Z}_{\lceil i/2
    \rceil}^{\smallsup{1}} {\hat Z}_{\lfloor
    i/2\rfloor}^{\smallsup{2}} } {L_ N} \right\}\right]
    +O(N^{-\beta}).
    \end{equation}
We will deal with the error term (\ref{ordefout}) at the end of
this step. Bound
    \begin{eqnarray*}
    &&\big| \sum_{i=2}^{j+1} Z_{\lceil i/2
    \rceil}^{\smallsup{1}}Z_{\lfloor i/2\rfloor}^{\smallsup{2}} -{\hat
    Z}_{\lceil i/2 \rceil}^{\smallsup{1}}{\hat Z}_{\lfloor
    i/2\rfloor}^{\smallsup{2}}
    \big|\leq  \sum_{i=2}^{j+1} Z_{\lceil i/2
    \rceil}^{\smallsup{1}} \big| Z_{\lfloor
    i/2\rfloor}^{\smallsup{2}}-{\hat Z}_{\lfloor
    i/2\rfloor}^{\smallsup{2}} \big|+ \sum_{i=2}^{j+1} {\hat
    Z}_{\lfloor i/2\rfloor}^{\smallsup{2}} \big| Z_{\lceil i/2
    \rceil}^{\smallsup{1}}-{\hat Z}_{\lceil i/2 \rceil}^{\smallsup{1}}
    \big|.
    \end{eqnarray*}
By Proposition \ref{prop-bdscoupling1} and uniformly in $j \leq
(1+2\eta)\log_{\nu}N$, we have, with probability exceeding
$1-O(N^{-\beta}\log_{\nu} N)$, that
    $$
    \max \left( \sum_{i=2}^{j+1} Z_{\lceil i/2 \rceil}^{\smallsup{1}}
    \big| Z_{\lfloor i/2\rfloor}^{\smallsup{2}}-{\hat Z}_{\lfloor
    i/2\rfloor}^{\smallsup{2}} \big| ,\sum_{i=2}^{j+1} {\hat
    Z}_{\lfloor i/2\rfloor}^{\smallsup{2}} \big| Z_{\lceil i/2
    \rceil}^{\smallsup{1}}-{\hat Z}_{\lceil i/2 \rceil}^{\smallsup{1}}
    \big|
    \right)=O(\nu^{(\frac12+\eta)\log_{\nu}N}N^{-\alpha})\sum_{i=2}^{j+1}
    {\hat Z}_{\lceil i/2 \rceil}^{\smallsup{1}}{\hat Z}_{\lfloor
    i/2\rfloor}^{\smallsup{2}}.
    $$
Since $\alpha>\frac12+\eta$, we have
$\nu^{(\frac12+\eta)\log_{\nu}N}N^{-\alpha}
=N^{\frac12+\eta-\alpha}=N^{-\alpha_1},$ for some $\alpha_1>0$.
Hence, for any $\eps$ with $0<\eps<\alpha_1$, where as before $\prob_{\sN}$ denotes the
conditional probability given the degrees $D_1,D_2,\ldots,D_{\sN}$,
and $\expec_{\sN}$ the expectation with respect to $\prob_{\sN}$, we have
    \begin{eqnarray*}
    &&\prob_{\sN}\left( \frac1{N} \big| \sum_{i=2}^{j+1} Z_{\lceil i/2
    \rceil}^{\smallsup{1}}Z_{\lfloor i/2\rfloor}^{\smallsup{2}} -{\hat
    Z}_{\lceil i/2 \rceil}^{\smallsup{1}}{\hat Z}_{\lfloor
    i/2\rfloor}^{\smallsup{2}} \big|
    >N^{-\eps}\right)\\
    &&\qquad\leq O(N^{-\beta}\log_{\nu} N)+
    \prob_{\sN}\left(\frac1{N}\sum_{i=2}^{j+1} {\hat Z}_{\lfloor
    i/2\rfloor}^{\smallsup{1}}{\hat Z}_{\lceil i/2
    \rceil}^{\smallsup{2}}
    >O\left(N^{\alpha_1-\eps}\right)\right)\\
    &&\qquad\leq O(N^{-\beta}\log_{\nu} N)+
    O\left(N^{\eps-\alpha_1}\right) \sum_{i=2}^{j+1}\expec_{\sN}[ {\hat
    Z}_{\lceil i/2 \rceil}^{\smallsup{1}}{\hat Z}_{\lfloor
    i/2\rfloor}^{\smallsup{2}}].
    \end{eqnarray*}
The involved conditional expectation can be computed explicitly
and we obtain
    $$
    \sum_{i=2}^{j+1}\expec_{\sN}[{\hat Z}_{\lceil i/2 \rceil}^{\smallsup{1}}
    {\hat Z}_{\lfloor i/2\rfloor}^{\smallsup{2}}]
    =D_1D_2\sum_{i=2}^{j+1} \nu_{\sN}^{\lceil i/2 \rceil-1}\nu_{\sN}^{\lfloor
    i/2\rfloor-1} =D_1D_2\sum_{i=0}^{j-1} \nu_{\sN}^{i}\le c D_1D_2
    \nu_{\sN}^{j},
    $$
for some constant $c$. Proposition \ref{prop-TV1} implies that we
can bound $\nu_{\sN}^{j}$ by $\nu^j(1+N^{-\alpha_2})^j$, with
probability exceeding $1-N^{-\beta_2}$, for some
$\alpha_2,\beta_2>0$, whereas Lemma \ref{lem-mom} implies $L_{\sN}^{-1}$ can be
replaced by $(\mu N)^{-1}$ with probability exceeding
$1-N^{-\beta_3}$, for some $\beta_3>0$. Putting this together we
obtain after taking the expectation with respect to
$D_1,D_2,\ldots,D_{\sN}$,
    \begin{eqnarray*}
    &&\prob\left(\frac1{L_{\sN}} \big| \sum_{i=2}^{j+1} Z_{\lfloor
    i/2\rfloor}^{\smallsup{1}}Z_{\lceil i/2 \rceil}^{\smallsup{2}}
    -{\hat Z}_{\lfloor i/2\rfloor}^{\smallsup{1}}{\hat Z}_{\lceil i/2
    \rceil}^{\smallsup{2}} \big|
    >N^{-\eps}\right)\\
    &&\qquad\le O(N^{-\beta}\log_{\nu}
    N)+O(N^{-\beta_1})+O(N^{-\beta_2})+O(N^{-\beta_3})+
    O\left(\frac{\nu^{j}(1+O(\log_{\nu}N/N^{\alpha_2}))}{N^{1+\alpha_1-\eps}}
    \right).
    \end{eqnarray*}
Since $\nu^{j}\le N^{1+2\eta}$ for $j\leq (1+2\eta)\log_{\nu}N$,
we obtain
    \begin{equation}
    \label{naarzhoed1} \prob\left(\frac1{L_{\sN}} \big| \sum_{i=2}^{j+1}
    Z_{\lceil i/2 \rceil}^{\smallsup{1}}Z_{\lfloor
    i/2\rfloor}^{\smallsup{2}} -{\hat Z}_{\lceil i/2
    \rceil}^{\smallsup{1}}{\hat Z}_{\lfloor i/2\rfloor}^{\smallsup{2}}
    \big|
    >N^{-\eps}\right)=O(N^{-\beta}),
    \end{equation}
for some $\beta>0$ by taking $\beta, \beta_2,\beta_3,\eta$ and
$\eps$ sufficiently small. For $x-y$ small, and $x,y\geq 0$, we find
$e^{-y}=e^{-x}+O(x-y)$, so that
    $$
    \exp\left\{ -\frac{\sum_{i=2}^{j+1} Z_{\lceil i/2
    \rceil}^{\smallsup{1}}Z_{\lfloor
    i/2\rfloor}^{\smallsup{2}}}{L_{\sN}} \right\} - \exp\left\{
    -\frac{\sum_{i=2}^{j+1}{\hat Z}_{\lfloor
    i/2\rfloor}^{\smallsup{1}}{\hat Z}_{\lceil i/2
    \rceil}^{\smallsup{2}}}{L_{\sN}} \right\}=O(N^{-\eps}),
    $$
with probability exceeding $1-O(N^{-\beta})$. In combination with the inequality
$e^{-x}\le 1$ for $x\geq 0$, we obtain (\ref{naarzhoed}).

We turn to the proof of (\ref{ordefout}) and the assumption that
$\sum_{k=1}^{\lceil j/2 \rceil}(Z_k^{\smallsup{1}}+Z_k^{\smallsup{2}})=o(L_{\sN})$.
From Proposition \ref{prop-bdscoupling1} and, uniformly in $j \leq
(1+2\eta)\log_{\nu}N$, we have with probability exceeding
$1-O(N^{-\beta}\log_{\nu} N)$ that
    \eq
    \label{ZhatZbd}
    \sum_{k=1}^{\lceil j/2\rceil}(Z_k^{\smallsup{1}
    }+Z_k^{\smallsup{2}}) \leq (1+O(N^{\frac12+\eta-\alpha}))
    \sum_{k=1}^{\lceil j/2\rceil}({\hat Z}_k^{\smallsup{1} }+{\hat
    Z}_k^{\smallsup{2}}).
    \en
so that, for all $i\le j$,
    $$
    \prob_{\sN}\left( \frac {\sum_{k=1}^{\lceil
    i/2\rceil}(Z_k^{\smallsup{1} }+Z_k^{\smallsup{2}})}{L_{\sN}^{3/4}}>
    N^{-\eps} \right) \leq O(N^{-\beta}\log_{\nu}
    N)+(1+O(N^{\frac12+\eta-\alpha})) \expec_{\sN}\left[
    \frac{\sum_{k=1}^{\lceil j/2\rceil}({\hat Z}_k^{\smallsup{1}
    }+{\hat Z}_k^{\smallsup{2}})}{N^{-\eps}L_{\sN}^{3/4}} \right].
    $$
Thus, in particular, using
(\ref{ZhatZbd}), $\sum_{k=1}^{\lceil j/2 \rceil}
(Z_k^{\smallsup{1}}+Z_k^{\smallsup{2}})=o(L_{\sN})$
on the above event.
Bounding the expectation of ${\hat Z}_k^{\smallsup{i}}$, we find
for $0<\eps<1/4$ and for all $i\le j\le (1+2\eta)\log_{\nu} N$,
    $$
    \prob\left( \frac {\sum_{k=1}^{\lceil j/2\rceil}(Z_k^{\smallsup{1}
    }+Z_k^{\smallsup{2}})}{L_{\sN}^{3/4}}> N^{-\eps} \right) \leq
    N^{-\beta}+(1+O(N^{-\alpha_1}))\frac{N^{\frac12+\eta}}{N^{\frac34-\eps}}=O(N^{-\beta}),
    $$
for some $\beta>0$. Hence, for $\eps_1>0$,
    $$
    \prob\left( \sum_{i=2}^{j+1} \frac{Z^{\smallsup{1}}_{\lceil
    i/2\rceil}Z_{\lfloor i/2\rfloor}^{\smallsup{2}}
        \sum_{k=1}^{\lceil i/2\rceil}(Z_k^{\smallsup{1} }+Z_k^{\smallsup{2}})}{L_{\sN}^2}
    >N^{-\eps_1}\right)\leq
    O(N^{-\beta})+ \prob\left( \sum_{i=2}^{j+1}
    \frac{Z^{\smallsup{1}}_{\lceil i/2\rceil}Z_{\lfloor
    i/2\rfloor}^{\smallsup{2}}}{L_{\sN}^{5/4}}
    >N^{\eps-\eps_1}\right).
    $$
By Proposition \ref{prop-bdscoupling1}, the product
$Z^{\smallsup{1}}_{\lceil i/2\rceil}Z_{\lfloor
i/2\rfloor}^{\smallsup{2}}$ can be bounded by
$(1+O(N^{\frac12+\eta-\alpha})) {\hat Z}^{\smallsup{1}}_{\lceil
i/2\rceil}{\hat Z}_{\lfloor i/2\rfloor}^{\smallsup{2}}$ and
$\expec[\sum_{i=2}^{j+1} \hat Z^{\smallsup{1}}_{\lceil
i/2\rceil}\hat Z_{\lfloor i/2\rfloor}^{\smallsup{2}}]\leq N^{1+2\eta}$,
while $L_{\sN}^{5/4}$ is of order $N^{5/4}$. Therefore, we obtain from the
Markov inequality that
    $$
    \prob\left( \sum_{i=2}^{j+1} \frac{Z^{\smallsup{1}}_{\lceil
    i/2\rceil}Z_{\lfloor i/2\rfloor}^{\smallsup{2}}
        \sum_{k=1}^{\lceil i/2\rceil}(Z_k^{\smallsup{1} }+Z_k^{\smallsup{2}})}{L_{\sN}^2}
    >N^{-\eps_1}\right)\leq
    O(N^{-\beta}),
    $$
for some $\beta>0$. Since $RM_{\sN}(j)$ is the difference of two numbers
between $0$ and $1$ and hence $|RM_{\sN}(j)|\le 1$, we obtain that, when $\eps_1\geq \beta$,
    \begin{equation}
    \expec[RM_{\sN}(j)]\leq N^{-\eps_1}+\prob\left( \frac{1}{L_{\sN}^2}\sum_{i=2}^{j+1}
    Z^{\smallsup{1}}_{\lceil i/2\rceil}Z_{\lfloor i/2\rfloor}^{\smallsup{2}}
    \sum_{k=1}^{\lceil i/2\rceil}(Z_k^{\smallsup{1}}+Z_k^{\smallsup{2}})
    >N^{-\eps_1}\right)\leq O(N^{-\beta}).
    \end{equation}
This proves (\ref{ordefout}).

\bigskip

\noindent {\bf Step 3: Coupling to the BP with offspring $\{g_j\}$.}
Corollary \ref{corr-coup1} combined with Lemma \ref{lem-mom} yields
    \begin{eqnarray*}
        \prob\Big(\frac{1}{L_{\sN}}\Big|\sum_{i=2}^{j+1}{\cal
    Z}_{\lceil i/2 \rceil}^{\smallsup{1}} {\cal Z}_{\lfloor i/2
    \rfloor}^{\smallsup{1}}-
            \sum_{i=2}^{j+1}{\hat
    Z}_{\lceil i/2 \rceil}^{\smallsup{1}} \hat Z_{\lfloor i/2
    \rfloor}^{\smallsup{2}}\Big|>N^{-\eps}\Big)
        =O(N^{-\beta}).
    \end{eqnarray*}
From this result we obtain, as in the first half of {\it step 2},
    $$\expec\left[\exp\left\{ -\frac{\sum_{i=2}^{j+1}{\hat Z}_{\lfloor
    i/2\rfloor}^{\smallsup{1}}{\hat Z}_{\lceil /2
    \rceil}^{\smallsup{2}}}{L_{\sN}} \right\}\right]=
    \expec\left[\exp\left\{ -\frac{\sum_{i=2}^{j+1} {\cal Z}_{\lceil
    i/2 \rceil}^{\smallsup{1}}{\cal Z}_{\lfloor
    i/2\rfloor}^{\smallsup{2}}}{L_{\sN}} \right\}\right]+O(N^{-\beta}),
    $$
where, as before, $\beta$ is a generic small positive number.
Using (\ref{compeven}) and the result of {\it step 2}, it follows that
    $$
    \prob(H_{\sN}>j)= \expec\left[
    \exp\left\{-
    \frac{\sum_{i=2}^{j+1} {\cal Z}_{\lceil i/2
    \rceil}^{\smallsup{1}} {\cal Z}_{\lfloor
    i/2\rfloor}^{\smallsup{2}} } {L_{\sN}} \right\}\right] +O(N^{-\beta}).
    $$
To obtain (\ref{hhoednaarexp}),
we finally replace, again at the cost of an additional term
$O(N^{-\beta})$, the random number $L_{\sN}$ by $\mu N(1+O(N^{-a}))$.

\bigskip

\noindent {\bf Step 4: Evaluation of the limit points.}
We start from (\ref{hhoednaarexp}) with
$j=k+\sigma_{\sN}\le (1+2\eta)\log_{\nu} N$, where
$\sigma_{\sN}=\lfloor\log_{\nu} N\rfloor$, to obtain
    \begin{equation}
    \label{biggerthanlogN}
    \prob(H_{\sN}>\sigma_{\sN}+k)= \expec\left[
    \exp\left\{
    \frac{-\sum_{i=2}^{\sigma_{\sN}+k+1} {\cal Z}_{\lceil i/2
    \rceil}^{\smallsup{1}} {\cal Z}_{\lfloor
    i/2\rfloor}^{\smallsup{2}} } {\mu N} \right\}\right]
    +O(N^{-\beta}).
    \end{equation}
We write $N=\nu^{\log_{\nu}N}=\nu^{\sigma_{\sN}-a_{\sN}},$
where we recall that $a_{\sN}=\lfloor
\log_{\nu}N\rfloor-\log_{\nu}N$. Then
$$
\frac{\sum_{i=2}^{\sigma_{\sN}+k+1} {\cal Z}_{\lceil i/2
\rceil}^{\smallsup{1}} {\cal Z}_{\lfloor
i/2\rfloor}^{\smallsup{2}}}{\mu N} = \mu \nu^{a_{\sN}+k}
\frac{\sum_{i=2}^{\sigma_{\sN}+k+1} {\cal Z}_{\lceil i/2
\rceil}^{\smallsup{1}} {\cal Z}_{\lfloor
i/2\rfloor}^{\smallsup{2}}}{\mu^2 \nu^{\sigma_{\sN}+k}}.
$$
In the above expression, the factor $\nu^{a_{\sN}}$ prevents proper
convergence. Without the factor $\mu \nu^{a_{\sN}+k}$, we obtain from
(\ref{martin}), with probability $1$,
    \begin{eqnarray*}
    \lim_{N\to \infty} \frac{\sum_{i=2}^{\sigma_{\sN}+k+1} {\cal
    Z}_{\lceil i/2 \rceil}^{\smallsup{1}} {\cal Z}_{\lfloor
    i/2\rfloor}^{\smallsup{2}}}{\mu^2 \nu^{\sigma_{\sN}+k}} =\frac{{\cal
    W}^ {\smallsup{1}}{\cal W}^{\smallsup{2}}}{\nu-1}.
    \end{eqnarray*}
Using (\ref{big-O}) we conclude that for each $\alpha >0$, there
is a $\beta>0$ such that
    \begin{eqnarray*}
    \prob\left( \left| \frac{\sum_{i=2}^{\sigma_{\sN}+k+1} {\cal Z}_{\lceil
    i/2 \rceil}^{\smallsup{1}} {\cal Z}_{\lfloor
    i/2\rfloor}^{\smallsup{2}}}{\mu^2 \nu^{\sigma_{\sN}+k}} -\frac{{\cal
    W}^ {\smallsup{1}}{\cal W}^{\smallsup{2}}}{\nu-1}\right|>O((\log
    N)^{-\alpha})\right)=O(N^{-\beta}).
    \end{eqnarray*}
Hence, for $k\leq 2\eta\log_{\nu} N$ and each $\alpha >0$,
    \begin{equation}
    \label{uitspraak}
    \prob(H_{\sN}>\sigma_{\sN}+k)
    =\expec\big(\exp\{-\kappa \nu^{a_{\sN}+k}{\cal W}^\smallsup{1}
    {\cal W}^\smallsup{2}\}
    \big)+O((\log N)^{-\alpha}),
    \end{equation}
where $\kappa=\mu/(\nu-1)$. This proves (\ref{uitspraakrep}).

We proceed by proving (\ref{hop}), with $R_{a}$ given in (\ref{Zdistr2}).
For this, we need to condition on node 1 and node 2 being connected.
Node 1 and node 2 are connected if and only if $H_{\sN}<\infty$.
Using (\ref{uitspraak}), for (\ref{hop}), it suffices to prove
that
    \begin{equation}
    \label{connclaim}
    \prob(H_{\sN}<\infty)=
    q^2+o(1),\quad
    \text{where}
    \quad
    q=\prob({\cal W}^\smallsup{1}>0).
    \end{equation}
We prove (\ref{connclaim}) using upper and lower bounds.
We note that, with $k=\eta \log_\nu N$,
    \begin{equation}
    \label{low-boundSP}
    \prob(H_{\sN}<\infty)\geq
    \prob(H_{\sN}\leq \sigma_{\sN}+k)=\expec\big(1-\exp\{-\kappa \nu^{a_{\sN}+k}{\cal W}^\smallsup{1}
    {\cal W}^\smallsup{2}\}
    \big)+O((\log N)^{-\alpha}).
    \end{equation}
Therefore,
    \begin{equation}
    \label{low-boundSP2}
    \prob(H_{\sN}<\infty)\geq
    q^2 \expec\big(1-\exp\{-\kappa \nu^{a_{\sN}+k}{\cal W}^\smallsup{1}
    {\cal W}^\smallsup{2}\}\big| {\cal W}^\smallsup{1}
    {\cal W}^\smallsup{2}>0\big)+O((\log N)^{-\alpha}).
    \end{equation}
By dominated convergence, for $k=2\eta \log_\nu N$, the
conditional expectation converges to 1, so that indeed
$\prob(H_{\sN}<\infty)\geq q^2+o(1)$.
For the upper bound, we rewrite, for any $m$,
    \begin{equation}
    \prob(H_{\sN}<\infty) = \prob(H_{\sN}<\infty, Z_m^{\smallsup{1}}Z_m^{\smallsup{2}}=0)
    +\prob(H_{\sN}<\infty, Z_m^{\smallsup{1}}Z_m^{\smallsup{2}}>0).
    \end{equation}
The second term is bounded from above by
    \begin{equation}
    \prob(H_{\sN}<\infty, Z_m^{\smallsup{1}}Z_m^{\smallsup{2}}>0)
    \leq \prob(Z_m^{\smallsup{1}}Z_m^{\smallsup{2}}>0)
    =\prob({\cal Z}_m^{\smallsup{1}}{\cal Z}_m^{\smallsup{2}}>0)+o(1)
    =q^2_m +o(1),
    \end{equation}
where we use Proposition \ref{prop-caft}, and we write
$q_m=\prob({\cal Z}_m^{\smallsup{1}}>0)$. When $m\rightarrow \infty$, we have that
$q_m\rightarrow q$, so that we are done when we can show that
for any $m$ fixed, $\prob(H_{\sN}<\infty, Z_m^{\smallsup{1}}Z_m^{\smallsup{2}}=0)=o(1).$
We note that if $Z_m^{\smallsup{1}}Z_m^{\smallsup{2}}=0$, then $H_{\sN}\leq m-1$.
Therefore, using (\ref{uitspraak}) with $k=m-\sigma_{\sN}-1$, we conclude
    \begin{equation}
    \prob(H_{\sN}<\infty, Z_m^{\smallsup{1}}Z_m^{\smallsup{2}}=0)\leq \prob(H_{\sN} \leq m-1)
    =\expec\big(1-\exp\{-\kappa \nu^{a_{\sN}+k}{\cal W}^\smallsup{1}
    {\cal W}^\smallsup{2}\}
    \big)+o(1)=o(1).
    \end{equation}
This completes the proof of (\ref{connclaim}).
We finally complete the proof of Theorems \ref{thm-tau>3} and \ref{thm-ll}
using (\ref{connclaim}),
which, together with (\ref{uitspraak}), implies that, for $k\leq 2\eta \log_\nu N$,
    \begin{equation}
    \prob(H_{\sN}\leq \sigma_{\sN}+k|H_{\sN}<\infty)
    =
    \expec\big(1-\exp\{-\kappa \nu^{a_{\sN}+k}{\cal W}^\smallsup{1}
    {\cal W}^\smallsup{2}\}| {\cal W}^\smallsup{1}
    {\cal W}^\smallsup{2}>0\big)+o(1).
    \end{equation}
\qed



\section{On the connected components}
\label{sec-clusters}

In this section, we will investigate the sizes of the connected components and prove
Theorem \ref{thm-cluster}.

\noindent {\bf Proof of Theorem \ref{thm-cluster}.} In the proof,
we will make essential use of the results in \cite{MR95,MR98},
where the statement in Theorem \ref{thm-cluster} is
proved for certain degree sequences. Indeed, denote by
    \eq
    d_i(N) = \sum_{j=1}^N I[D_j=i], \qquad i=0,1,\ldots,
    \label{din}
    \en
the degree sequence of our random graph $G$, where
$D_1,D_2,\ldots,D_{\sN}$ is the i.i.d.\  sequence with distribution $F$
introduced in (\ref{kansen}) and satisfying (\ref{distribution}).
In \cite{MR95}, the bounds on the connected components in Theorem
\ref{thm-cluster} are proved with only a lower bound on the
largest connected component size, while in \cite{MR98}, the asymptotic size of
the largest connected component is determined. Both papers assume a number of
hypotheses on the degree sequence $\{d_i(N)\}_{i\geq 0}$. Thus,
Theorem \ref{thm-cluster} follows when we can
show that the probability that our degree sequences in (\ref{din})
satisfy the restrictions is at least $1-o(1)$. In fact, we need to
alter the random graph $G$ in a certain way to meet the conditions
of Molloy and Reed, and subsequently need to prove that the
alteration does not affect the results. We now go over their
conditions and definitions.

Firstly, the degree sequence needs to be {\bf feasible}, meaning
that there exists at least one graph with the degree sequence.
This is true, since $L_{\sN}$ is even and we have that
$$
\sum_{i=1}^\infty i d_i(N) =\sum_{i=1}^\infty i \sum_{j=1}^N
I[D_j=i] = \sum_{j=1}^N\sum_{i=1}^\infty  i I[D_j=i] =\sum_{j=1}^N
D_j=L_{\sN}.
$$

Secondly, the degree sequence needs to be {\bf smooth}, meaning
that for some sequence $\lambda_i$, we have
    \[\lim_{N\rightarrow \infty}  \frac{d_i(N)}{N}=\lambda_i.\]
In our setting, this follows almost surely from the law of large
numbers, with $\lambda_i=f_i=\prob(D=i)$.

Thirdly, and this is the most serious condition, the degree
sequence needs to be {\bf well-behaved}, meaning that it is
smooth, feasible, and that for every $\epsilon'$, there exists
$N'=N'(\epsilon')$, such
that for all $N>N'$, we have that\\
1.
    \eq
    \label{cond1}
    \sup_i \big| i(i-2) \frac{d_i(N)}{N}-i(i-2) \lambda_i\big |<\epsilon';
    \en
2. there exists $i^*$
    \eq
    \label{cond2}
    \Big|\sum_{i=1}^{i^*} i(i-2) \frac{d_i(N)}{N} -\sum_{i=1}^{\infty} i(i-2) \lambda_i\Big|
    \leq \epsilon';
    \en
3. there exists an $\epsilon>0$ such that
$d_i(N)=0$ for all $i\geq \lceil N^{\frac 14 -\epsilon}\rceil$.

We start with the last assumption, which is not satisfied by our
graph. Indeed, the last restriction means that all nodes have
degree at most $\lceil N^{\frac 14 -\epsilon}\rceil-1$. We will
first alter the graph, and thus the degree sequences, in the
following way. Fix $\epsilon>0$ small. For nodes $j$ with $D_j\geq
\lceil N^{\frac 14 -\epsilon}\rceil$, we remove $D_j-\lceil
N^{\frac 14 -\epsilon}\rceil+1$ edges. We do this
by first removing in a uniform way edges between pairs $i,j$ where
the degrees  of $D_i$ and $D_j$ {\it both} exceed $\lceil N^{\frac
14 -\epsilon} \rceil-1$. When there are no more edges between
nodes with degrees exceeding $\lceil N^{\frac 14 -\epsilon}
\rceil-1$, we remove edges uniformly from the nodes with degrees
exceeding $\lceil N^{\frac 14 -\epsilon} \rceil-1$. Thus, we end
up with a graph $G'$ such that all degrees are at most $\lceil
N^{\frac 14 -\epsilon}\rceil-1$. Moreover each node $j$ for which
$D_j\geq \lceil N^{\frac 14 -\epsilon}\rceil$ has degree equal to
$\lceil N^{\frac 14 -\epsilon}\rceil-1$ in the altered graph $G'$.
This will be the graph to which we apply the results of Molloy
and Reed. Let $D_j'$ be the degree of the node $j$ in $G'$, and
write $d_i'(N)$ for the number of nodes with degree equal to $i$
in $G'$. Then $d_i'(N)=0$ for $i\geq \lceil N^{\frac 14 -\epsilon}
\rceil,$ as required.

We first compute the number of removed edges, which we denote by
$R_{\sN}$. Its expectation is bounded above by
    \begin{eqnarray*}
    \expec[R_{\sN}]&\leq&\expec\big[\sum_{j=1}^N (D_j+1- \lceil N^{\frac 14 -\epsilon}\rceil)^+
    I[D_j\geq \lceil N^{\frac 14 -\epsilon}\rceil]\big]
    \leq N\sum_{l\geq \lceil N^{\frac 14 -\epsilon}\rceil-1}
    \prob(D_1>l)\\
    &\leq&c N \sum_{l\geq \lceil N^{\frac 14 -\epsilon}\rceil-1 }
    l^{-\tau+1} =C N^{1-(\tau-2)(\frac 14 -\epsilon)}<N^{\frac34},
    \end{eqnarray*}
for $\tau>3$ and $\epsilon$ sufficiently small. We are hence removing only a fraction of the $L_{\sN}$
available edges and all degrees go down (see Lemma \ref{lem-mom}
that $L_{\sN}$ is close to $\mu N$). Moreover, with probability
converging to one, we have that $R_{\sN}\leq 2 N^{\frac34}$, since by
a computation analogous to the one given above for $\expec[R_{\sN}]$, we have
$\mbox{Var}(R_{\sN})\leq CN^{1-(\tau-3)(\frac 14 -\epsilon)}$, so that by
the Chebychev inequality,
\begin{equation}
\prob(R_{\sN}> 2 N^{\frac34}) \leq \prob(|R_{\sN}-\expec[R_{\sN}]|> N^{\frac34})
\leq N^{-\frac32}\mbox{Var}(R_{\sN})
\leq C N^{-\frac12-(\tau-3)(\frac14-\epsilon)}\leq CN^{-\frac12}.
\end{equation}

We start by checking (\ref{cond1}) for the graph $G'$, with
$\lambda_i=f_i$ in (\ref{kansen}). For this, we will use the
following bound from \cite[Corollary 1.4(i)]{bol01}, which
states that if $S_{\sN}$ is binomial with
parameters $N$ and $p$, and if $x =(Np(1-p))^{1/2} \geq 1$, then
    \eq
    \prob\big( |S_{\sN} -Np|\geq x (Np(1-p))^{1/2}\big) \leq \frac 1x e^{-x^2/2}.
    \label{binbd}
    \en
We first check condition (\ref{cond1}) for $i=\lceil N^{\frac 14
-\epsilon}\rceil-1 $. By construction, we have that for $i=\lceil
N^{\frac 14 -\epsilon}\rceil-1$,
    \eq
    \label{bdd's}
    d_i'(N)= \sum_{j\geq i} d_j(N).
    \en
Hence, $d_i'(N)$ is a binomial
random variable with parameters $N$ and $p=1-F(\lceil N^{\frac 14
-\epsilon}\rceil-2)$. Thus, by (\ref{binbd}), with $x=C\sqrt{\log
N}$, we have that
    \eq
    \label{bdd's2}
    \prob\big( |d'_i(N) -Np|\geq C ((\log{N}) Np(1-p))^{1/2}\big) \leq \frac 1C N^{-C^2/2}.
    \en
Thus, we have that
for $i=\lceil N^{\frac 14 -\epsilon}\rceil-1$, $\lambda_i=f_i$ and $p=1-F(\lceil N^{\frac 14
-\epsilon}\rceil-2)$,
    \begin{eqnarray*}
    i(i-2)\big|\frac{d_i'(N)}{N} - \lambda_i\big|
    &\leq& i^2\big|\frac{d_i'(N)}{N}-p+p - \lambda_i\big|
    \leq i^2|\frac{d_i'(N)}{N}-p|+i^2|p - \lambda_i|\\
    &\leq& i^2 \frac{C\sqrt{\log{N}}}{\sqrt{N}} [1-F(\lceil N^{\frac 14 -\epsilon}\rceil-2)]^{1/2}
    +i^2 f_i +i^2 [1-F(\lceil N^{\frac 14 -\epsilon}\rceil-2)]    \\
    &\leq& CN^{(\frac 12-2\epsilon)}\cdot \frac{\log N}{\sqrt{N}}\cdot N^{\frac12(1-\tau)(\frac14-\epsilon)}
    +2N^{\frac12-2\epsilon} \cdot c\lceil N^{\frac 14 -\epsilon}\rceil^{1-\tau}<\epsilon' ,
    \end{eqnarray*}
for $\tau>3$. This proves (\ref{cond1}) for $i=\lceil N^{\frac 14
-\epsilon}\rceil-1$.

We next prove (\ref{cond1}) for $i<\lceil N^{\frac 14
-\epsilon}\rceil-1$. For this, we use the triangle inequality
    \eq
    \label{split}
    i(i-2)\big|\frac{d_i'(N)}{N} - \lambda_i\big| \leq i^2\big|\frac{d_i'(N)}{N}-\frac{d_i(N)}{N}\Big|+
    i^2\big|\frac{d_i(N)}{N} - \lambda_i\big|,
    \en
and we bound these two terms separately.

We start with the second term, and use (\ref{binbd}), which gives
that
    \eq
    \prob\big( |d_i(N)-Nf_i| \geq C (f_i N\log{N})^{1/2}\big)\leq N^{-C^2/2}.
    \en
We will take $C>2$, so that
    \eq
    \label{bdbin}
    \prob\big(\exists i< \lceil N^{\frac 14 -\epsilon}\rceil -1: |d_i(N)-Nf_i| \geq C (f_i N\log{N})^{1/2}\big)\leq
    \sum_{i=1}^N N^{-C^2/2}= N^{1-C^2/2}.
    \en
On the complementary event, we have that
    \eq
    \sup_{i< \lceil N^{\frac 14 -\epsilon}\rceil-1} \big| i(i-2) \frac{d_i(N)}{N}-i(i-2) \lambda_i\big |
    \leq C \sup_{i< \lceil N^{\frac 14 -\epsilon}\rceil-1} i^2 \big(\frac{f_i\log{N}}{N}\big)^{1/2} =o(1).
    \en
Thus, we have bounded the second term in (\ref{split}). We next
turn to the first term in (\ref{split}). First, we clearly have that
$|d'_i(N)-d_i(N)|\leq R_{\sN}$. Thus, since $R_{\sN}\leq 2 N^{\frac34}$,
    \[
    i^2\big|\frac{d_i'(N)}{N}-\frac{d_i(N)}{N}\Big|
    \leq i^2 \frac{R_{\sN}}{N} \leq 2i^2 N^{-\frac14}
    \leq 2N^{-\frac14+2(\frac18-\epsilon)}\leq \epsilon',
    \]
for $i\leq N^{\frac18-\epsilon}$. For $i> N^{\frac18-\epsilon}$,
we bound  $d_i'(N)\leq \sum_{j\geq i} d_j(N)$, so that, again using
(\ref{bdd's}--\ref{bdd's2}),
\[
 i^2\big|\frac{d_i'(N)}{N}-\frac{d_i(N)}{N}\Big|
 \le
 \frac{2i^2}{N}\sum_{j\geq i} d_j(N)= 2i^2(1-F(i-1))(1+o(1))\leq 2c N^{(\frac18-\eps)(3-\tau)}\to 0.
\]

To check (\ref{cond2}), we first take $i^*$ fixed so that
    \eq
    \sum_{i=i^*+1}^{\infty} i(i-2) \lambda_i
    \leq \epsilon'/2.
    \en
This is possible, since $\expec[D^2]<\infty$. Thus, we are left to show that
    \eq
    \sum_{i=1}^{i^*} i(i-2) \big|\frac{d_i(N)}{N}-\lambda_i\big|
    \leq \epsilon'/2.
    \en
In order to do so, we use the bound in (\ref{bdbin}) to
obtain that
    \eq
    \sum_{i=1}^{i^*} i(i-2) \big|\frac{d_i(N)}{N}-\lambda_i\big|
    \leq C\sum_{i=1}^{i^*} i^2\big(\frac{f_i\log{N}}{N}\big)^{1/2}\leq C(i^*)^3 \big(\frac{\log{N}}{N}\big)^{1/2}
    \leq \epsilon'/2,
    \en
whenever $N$ is sufficiently large. The same result applies to
$d_i'(N)$, since $|d_i'(N)-d_i(N)|\leq R_{\sN}$, and $R_{\sN}=o(N)$,
so that
\[
\sum_{i=1}^{i^*} i(i-2) \big|\frac{d_i(N)}{N}-\frac{d'_i(N)}{N}|
\leq (i^*)^3 \frac{R_{\sN}}{N}=o(1).
\]

Therefore, we have proved all conditions for the
graph $G'$, and thus obtain the result in Theorem
\ref{thm-cluster} for $G'$. To complete the proof, we need to show
that the result for $G'$ implies the result for $G$.

This implication is proved in several small steps. First, denote
the largest connected components of $G$ and $G'$ by $LC_G$ and $LC_{G'}$.
Since $G$ can be obtained from $G'$ by adding the removed edges
back, we obtain that (since we put back at most $R_{\sN}$ connected components of
size at most $\gamma \log N$),
    \eq
    |LC_{G'}|\leq |LC_G| \leq |LC_{G'}| + R_{\sN} \cdot \gamma \log{N}.
    \en
Thus, since $|LC_{G'}|= qN(1+o(1))$ and $R_{\sN}\leq 2N^{\frac34}$ with
probability $1+o(1)$, we
obtain that
    \eq
    qN(1+o(1)) \leq |LC_G| \leq qN(1+o(1)) + O(N^{\frac34}\log{N})=qN(1+o(1)),
    \en
so that the largest connected component has size $qN(1+o(1))$ with
probability $1+o(1)$, as claimed.

To see that all other connected components in $G$ have size at most $\gamma
\log{N}$, we note that in $G'$ the removed edges are all connected
to nodes with degree $\lceil N^{\frac 14 -\epsilon}\rceil$. We
first show that with overwhelming probability these nodes are
already in the largest connected component in $G'$. Since in $G'$ only the
largest connected component has at least $N^{\delta}$ nodes
for any $\delta>0$ and since $\gamma \log{N}=o(N^{\delta})$, it suffices
to check that nodes in $G'$ with degree $\lceil N^{\frac 14
-\epsilon}\rceil$ are connected to at least $N^{\delta}$ other
nodes. Since the probability of picking a node different from the
ones already connected to the node under observation is bounded
from below by $1-N^{2(\frac 14 -\epsilon)-1}$ (since all degrees
in $G'$ are bounded above by $\lceil N^{\frac 14
-\epsilon}\rceil$), the probability that at most $N^{\delta}$
different nodes are chosen is bounded by the probability that a
binomial random variable, with parameters $p=1-N^{2(\frac 14 -\epsilon)-1}$ and
$n=\lceil N^{\frac 14 -\epsilon}\rceil$, is bounded from above by
$N^{\delta}$. By (\ref{binbd}), this probability is negligible
whenever $\delta <\frac 14 -\epsilon$. Thus, we may assume that
all nodes with degree $\lceil N^{\frac 14 -\epsilon}\rceil$ are in
the largest connected component in $G'$. Therefore, we obtain that the nodes
that must be added to $G'$ to form $G$ are attached to the largest
connected component of $G'$. Thus, the size of the second largest connected component of
$G$ is bounded from above by the size of the second largest
connected component of $G'$, which is bounded from above by $\gamma \log N$.
\qed

\bigskip

\noindent{\bf Acknowledgement.}
The work of RvdH was supported in part by Netherlands Organisation for
Scientific Research (NWO). We thank Dmitri Znamenski for the Figures 1 and 2
and for useful comments on a previous version. We thank the two referees for
many suggestions that improved on the readability of the paper.


\renewcommand{\thesection}{\Alph{section}}
\setcounter{section}{0}

\numberwithin{equation}{subsection}
\numberwithin{theorem}{subsection}
\section{Appendix.}

\subsection{Proof of Proposition \ref{prop-TV1} } \label{appendix}
In this part of the appendix, we prove Proposition \ref{prop-TV1},
which we restate here for convenience as Proposition
\ref{prop-TV2}. At the end of this section, we restate and prove
Corollary \ref{corr-coup1}.

\begin{prop}
    \label{prop-TV2}
    There exist $\alpha_2,\beta_2>0$ such that
    \begin{equation}
    \prob\Big(\sum_{n=0}^{\infty} (n+1)|g_n^{\smallsup{N}}-g_n| \geq N^{-\alpha_2}\Big)\leq N^{-\beta_2}.
    \end{equation}
\end{prop}
\noindent
In the proof, we need the following lemma.
\begin{lemma}
    \label{lem-boundsf}
    Fix $\tau>1$. For each non-negative integer $s$, there exists
    a constant $C>0$, such that
        \begin{eqnarray}
        \sum_{j=m}^{n} (j+1)^s f_{j+1} &\leq & C m^{-(\tau-1-s)}+Ch(n).
        \label{boundf}
        \end{eqnarray}
    where
    $$
    h(n)=
    \left\{
    \begin{array}{ll}
    0,& s<\tau-1,\\
    \log(n+1),&s=\tau-1,\\
    (n+1)^{s-\tau+1},& s>\tau-1.
    \end{array}
    \right.
    $$
    \end{lemma}
\noindent
We defer the proof of Lemma \ref{lem-boundsf} to the end of this
section.

\noindent {\bf Proof of Proposition \ref{prop-TV2}.} Fix $a,b,
\alpha>0$. Define
    \begin{eqnarray}
    F&=&\left\{|\frac{L_{\sN}}{\mu N}-1|\leq N^{-\alpha}\right\}
    \cap \left \{\frac{1}{N}\sum_{i=1}^N (D_i+1)^2 I[D_i \geq N^a] \leq N^{-b}\right\}
    \label{evFdef}\\
    &&\qquad \cap \left\{ \frac{1}{N}\sum_{n=0}^{N^a} (n+1)^2 \big|\sum_{i=1}^N\big(I[D_i=n+1]-f_{n+1}\big)\big| \leq N^{-b}
    \right\}.
    \nonumber
    \end{eqnarray}
The constants $a,b$ and $\alpha$ will be chosen appropriately in the
proof. The strategy of the proof is as follows. We will prove that
    \eq
    \label{probFbd}
    \prob(F^c)\leq N^{-\beta_2},
    \en for some $\beta_2>0$, and that on
$F$,
    \eq
    \label{onFbd}
    \sum_{n=0}^{\infty} (n+1)|g_n^{\smallsup{N}}-g_n|\leq
    N^{-\alpha_2},
    \en
for some $\alpha_2$. This proves Proposition
\ref{prop-TV2}. We start by showing (\ref{onFbd}).

We bound
    \eq
    \label{split2}
    \sum_{n=0}^{\infty} (n+1)|g_n^{\smallsup{N}}-g_n|
    \leq \sum_{n=0}^{\infty} (n+1)|g_n^{\smallsup{N}}-\frac{N\mu}{L_{\sN}}g_n|
    +(\nu+1) \big|\frac{N\mu}{L_{\sN}}-1\big|.
    \en
The second term is bounded by $(\nu+1) N^{-\alpha}$ by the first event
in $F$. The first term in (\ref{split2}) can be bounded, for $N$
sufficiently large, as, again using the first event in $F$,
    \begin{eqnarray}
    \sum_{n=0}^{\infty} (n+1)|g_n^{\smallsup{N}}-\frac{N\mu}{L_{\sN}}g_n|
    &=&\frac{1}{L_{\sN}}\sum_{n=0}^{\infty} (n+1)^2\big|\sum_{i=1}^N \big(I[D_i=n+1]-f_{n+1}\big)\big|\\
    &\leq& \frac{2}{\mu N} \sum_{n=0}^{\infty} (n+1)^2\big|\sum_{i=1}^N \big(I[D_i=n+1]-f_{n+1}\big)\big|\nonumber.
    \end{eqnarray}
We next split the sum over $n$ into $n> N^{a}$ and $n\leq N^{a}$
for some appropriately chosen $a\in (0,1]$. On $F$, the contribution from
$n\leq N^{a}$ is at most $\frac{1}{\mu}N^{-b}$, whereas we can bound the contribution from
$n>N^a$ by
    \begin{eqnarray}
    \frac{2}{\mu N}\sum_{n=N^a}^{\infty} (n+1)^2\sum_{i=1}^N (I[D_i=n+1]+f_{n+1})
    =\frac{2}{\mu N}\sum_{i=1}^N (D_i+1)^2 I[D_i \geq N^a] +\frac2{\mu}\sum_{n=N^a}^{\infty}(n+1)^2 f_{n+1}.
    \nonumber
    \end{eqnarray}
For $\tau>3$, the second term is bounded by $CN^{-a(\tau-3)}$ by Lemma
\ref{lem-boundsf}. The first term is bounded by $\frac{\mu}{2}N^{-b}$ by the
second event in $F$. Thus, we obtain (\ref{onFbd}) with
$\alpha_2=\min\{b,a(\tau-3)\}$.

We now prove (\ref{probFbd}). For this, we use
that $F$ is an intersection of three events which we will write as
$F_1, F_2$ and $F_3$, so that
    \eq
    \prob(F^c)\leq \prob(F^c_1)+\prob(F^c_2)+\prob(F^c_3).
    \en
The first probability is bounded by $\prob(F^c_1)\leq
c\cdot  N^{2\alpha-1}$, by Lemma \ref{lem-mom}. For $\prob(F^c_2)$, we
use the Markov inequality, to obtain that
    \eq
    \prob(F^c_2)\leq N^b \mathbb{E}\big[(D_1+1)^2 I[D_1 \geq N^a]\big]
     \leq N^{b-a(\tau-3)},
     \en
by Lemma \ref{lem-boundsf}.
For $\prob(F^c_3)$, we use in turn the Markov
inequality, Cauchy-Schwarz in the form
$\sum_{n=0}^{N^a} b_n \le (\sum_{n=0}^{N^a} 1^2 \sum_{n=0}^{N^a} b_n^2)^{\frac12}$,
and the Jensen inequality applied to $x \mapsto \sqrt{x}$ (a concave function), to obtain
    \begin{eqnarray}
    \prob(F^c_3)
    &\leq& N^{b-1}\mathbb{E}\big[\sum_{n=0}^{N^a} (n+1)^2
    \big|\sum_{i=1}^N\big(I[D_i=n+1]-f_{n+1}\big)\big|\big]\\
    &\leq & N^{b-1  } (N^a+1)^{\frac12}\mathbb{E}
    \Big(\sum_{n=0}^{N^a} (n+1)^4\big(\sum_{i=1}^N \big(I[D_i=n+1]-f_{n+1}\big)\big)^2\Big)^{1/2}
    \nonumber\\
    &\leq & 2N^{b+a/2-1}\sum_{n=0}^{N^a} (n+1)^4
    \mathbb{E}\big(\sum_{i=1}^N \big(I[D_i=n+1]-f_{n+1}\big)\big)^2\big]\Big)^{1/2}
    \nonumber\\
    &\leq & 2N^{b+a/2-1}\Big(\sum_{n=0}^{N^a} (n+1)^4 Nf_{n+1} \Big)^{1/2}
    \leq 2N^{b+a/2-1/2} N^{a\max\{0,5-\tau\}/2},\nonumber
    \end{eqnarray}
where in the last inequalities, we have used Lemma \ref{lem-boundsf} and
    $$
    \expec\big[\big(\sum_{i=1}^N[I[D_i=n+1]-f_{n+1}]\big)^2\big]
    =\mbox{Var}\left(\sum_{i=1}^N I[D_i=n+1]\right)=Nf_{n+1}(1-f_{n+1})\leq Nf_{n+1}.
    $$
Thus, we obtain
the statement in Proposition \ref{prop-TV2} with
    $$
    \beta_2=\min\{1/2-b-a\max\{1,6-\tau\}/2, a(\tau-3)-b,
    (2\alpha-1)\}.
    $$
By picking first $b$ small, and then $a$ small, we
see that $\alpha_2, \beta_2>0$.
\qed

\begin{remark}
\label{exttau2}
When (\ref{distribution}) holds for some $\tau>2$ (rather than $\tau>3$),
then the above proof can be repeated to show that
    \begin{equation}
    \prob\big(\sum_{n=0}^{\infty}|g_n^{\smallsup{N}}-g_n| \geq N^{-\alpha_2}\big)\leq N^{-\beta_2}.
    \end{equation}
Indeed, in the definition of the event $F$ in (\ref{evFdef}), we can replace
$(D_i+1)^2$ by $(D_i+1)$ in the second event, and $(n+1)^2$ by
$(n+1)$ in the third event. Then, by adapting the above argument,
the event $F$ implies that
$\sum_{n=0}^{\infty}|g_n^{\smallsup{N}}-g_n| \leq N^{-\alpha_2}$.
The proof that $\prob(F^c)\leq N^{-\beta_2}$ can be adapted accordingly.

\end{remark}

\medskip

\noindent {\bf Proof of Lemma \ref{lem-boundsf}.} Define a density
$f(x)=\sum_{j=0}^{\infty} f_j I[j\le x < j+1]$, and the
corresponding distribution function ${\tilde F}(x)=\int_0^x
f(u)\,du$. Then for integer-valued $j>0$,
    $$
    {\tilde F}(j)=f_0+\ldots+f_{j-1}=F(j-1), \quad F(j-1)\leq {\tilde
    F}(x)\le F(j),\quad x \in (j,j+1).
    $$
Moreover
    $$
    \sum_{j=m}^{n} (j+1)^s f_{j+1}\leq \int_{m+1}^{n+2}
    x^sf(x)\,dx=-\int_{m+1}^{n+2} x^s \, d(1-{\tilde F}(x)).
    $$
Using partial integration and the upper bound
    $$
    1-{\tilde F}(x)\le 1-F(j-1)\leq c (j-1)^{1-\tau},
    $$
for $x\in (j,j+1)$, we conclude that
    \begin{eqnarray*}
    \sum_{j=m}^{n} (j+1)^s f_{j+1}
    &\leq&(m+1)^s(1-{\tilde F}(m+1))-(n+2)^s(1-{\tilde F}(n+2))
    +\int_{m+1}^{n+2} (1-{\tilde F}(x))\, dx^s\\
    &\leq& c\left[
    m^{1+s-\tau}+\int_m^{n+1} y^{s-\tau}\,dy \right].
    \end{eqnarray*}
This yields the upper bound.
\qed

\medskip

We finally prove Corollary \ref{corr-coup1}. In order to do so, we first formulate
and prove an intermediate result. This result will be followed by the
reformulation of Corollary \ref{corr-coup1}, which now becomes
Corollary \ref{corr-coup2}, and its proof.
\begin{prop}
There exist $\eps,\beta,\eta>0$ such that for all $j\le (\frac12+\eta)\log_{\nu} N$,
as $N\to \infty$,
    \label{prop-coup}
   \begin{equation}
   \label{verg-coup}
     \prob\Big(\frac{1}{\sqrt{N}}\Big|\sum_{i=1}^{j}{\cal
    Z}_i^{\smallsup{1}}-
        \sum_{i=1}^{j}{\hat
Z}_i^{\smallsup{1}}\Big|>N^{-\eps}\Big)
    =O(N^{-\beta}).
   \end{equation}
\end{prop}

\proof
Let
    \begin{equation}
        \label{ENdef}
        F_{\sN}=\{ \sum_{n=0}^{\infty} n|g_n^{\smallsup{N}}-g_n|< N^{-\alpha_2}\},
    \end{equation}
then according to Proposition \ref{prop-TV2} we have
$\prob(F^c_{\sN})\le N^{-\beta_2}$. We claim that for all $i\geq 1$,
\begin{equation}
    \label{diffexpbd}
    \expec_{\sN} |{\cal Z}_i^{\smallsup{1}}
    -{\hat Z}_i^{\smallsup{1}}| \leq
    \max\{\nu-\alpha_{\sN}, \nu_{\sN}-\alpha_{\sN}\}
    \sum_{m=1}^i \expec_{\sN}
    [{\hat Z}_{m}^{\smallsup{1}}]
    (\max\{\nu,\nu_{\sN}\})^{i-m},
 \end{equation}
where
    \eq
    \alpha_{\sN} = \sum_{n=0}^\infty n \min\{g_n, g_n^{\smallsup{N}}\}=\nu -
    \sum_{n=0}^\infty n \big(g_n-\min\{g_n, g_n^{\smallsup{N}}\}\big)
    =\nu_{\sN} -
    \sum_{n=0}^\infty n \big(g_n^{\smallsup{N}}-\min\{g_n, g_n^{\smallsup{N}}\}\big).
    \en
We first prove
(\ref{diffexpbd}). For ${\cal Z}_i^{\smallsup{1}}\neq
\hat Z_i^{\smallsup{1}}$, the coupling is not
successful in at least one of the generations $m, 1\le m \le i$.
Let $m$ be the first generation for which the coupling
is unsuccessful. There are at most $\hat Z_m^{\smallsup{1}}$ nodes
for which the coupling can fail. If the coupling
fails for a node, the expected difference between the offspring of that node is bounded
above by $\max\{\nu-\alpha_{\sN}, \nu_{\sN}-\alpha_{\sN}\}$. Finally, from
generation $m+1$ on, we again have two BP's with laws $g$
and $g^{\smallsup{N}}$, so that the expected offspring is bounded
by $(\max\{\nu,\nu_{\sN}\})^{i-m}$. This demonstrates the claim (\ref{diffexpbd}).

Furthermore, since $\expec_{\sN}[{\hat
Z}_{m}^{\smallsup{1}}]=D_1 \nu_{\sN}^{m-1}$, we end up with
   \eq
   \label{fournine}
    \expec_{\sN}|{\cal Z}_i^{\smallsup{1}}-{\hat Z}_i^{\smallsup{1}}| \leq
    \max\{\nu-\alpha_{\sN}, \nu_{\sN}-\alpha_{\sN}\}
    i D_1   (\max\{\nu, \nu_{\sN}\})^{i-1}.
    \en
By (\ref{ENdef}), on $F_{\sN}$ we have that
    \begin{eqnarray*}
    &&\max\{\nu-\alpha_{\sN}, \nu_{\sN}-\alpha_{\sN}\}\leq \sum_n n
        |g_n-g_n^{\smallsup{N}}|
        < N^{-\alpha_2},\\
    &&\frac{\max\{\nu, \nu_{\sN}\}}{\nu}
    =1+\nu^{-1}\max\{0,\sum_n n(g_n-g_n^{\smallsup{N}})\}=1+O(N^{-\alpha_2}).
    \end{eqnarray*}
Hence, for $j\le (\frac12+\eta)\log_{\nu} N$, using the abbreviation
    $$
    T_{\sN}=
    \frac{1}{\sqrt{N}}\Big|\sum_{i=1}^{j}{\cal Z}_i^{\smallsup{1}}-
        \sum_{i=1}^{j}\hat Z_i^{\smallsup{1}}\Big|,
    $$
we have
    \begin{eqnarray*}
    \prob\Big(T_{\sN}>N^{-\eps}\Big)
    &\leq & \prob(F_{\sN}^c)+\prob(T_{\sN}> N^{-\eps},F_{\sN})
    \leq N^{-\beta_2}+\expec\left[\prob_{\sN}\left(T_{\sN}I_{F_{\sN}}> N^{-\eps}\right)\right]\\
    & \leq &    N^{-\beta_2}+\expec\left[N^{\eps}\expec_{\sN}\left[T_{\sN}I_{F_{\sN}}\right]\right].
    \end{eqnarray*}
From (\ref{fournine}) and the estimates on $F_{\sN}$, we obtain
    \begin{eqnarray*}
    \expec\left[N^{\eps}\expec_{\sN}\left[T_{\sN}I_{F_{\sN}}\right]\right]
    &\leq& N^{\eps-\frac12}
    \expec_{\sN}\left[\sum_{i=1}^{j}
    |{\cal Z}_i^{\smallsup{1}}- \hat Z_i^{\smallsup{1}}|\cdot I_{F_{\sN}}\right]\\
    & \leq& \nu^j N^{\eps-\frac12}\expec[D_1\sum_{i=1}^j N^{-\alpha_2} i(1+O(N^{-\alpha_2}))^{i-1}]\\
    &\leq & \mu N^{\eps+\eta-\alpha_2}\sum_{i=1}^{\lfloor(\frac12+\eta)\log_{\nu}N \rfloor} i(1+O(N^{-\alpha_2}))^{i-1}\\
    &\leq & \mu N^{\eps+\eta-\alpha_2}\cdot (\log_{\nu} N)^2 \cdot N^{(\frac12+\eta)\log_{\nu}(1+O(N^{-\alpha_2}))},
    \end{eqnarray*}
using that for $x=1+O(N^{-\alpha_2})>1$, we have $\sum_{i=1}^n ix^{i-1}\leq n^2x^n$.
This proves the proposition since $\log_{\nu} N \cdot N^{(\frac12+\eta)\log_{\nu}(1+O(N^{-\alpha_2}))}$
can be bounded by any small power of $N$, and $\eps$ and $\eta$ can both be taken arbitrarily small,
whereas $\alpha_2>0$.

We finally restate and prove Corollary \ref{corr-coup1}.
\qed
\begin{corr}
There exist $\eps,\beta,\eta>0$ such that for all $j\le (1+2\eta)\log_{\nu} N$,
as $N\to \infty$,
    \label{corr-coup2}
   \begin{equation}
   \label{vgl-coup2}
   \prob\Big(\frac 1N\Big|\sum_{i=1}^{j}{\cal
    Z}_{\lceil i/2\rceil}^{\smallsup{1}}{\cal Z}_{\lfloor i/2\rfloor}^{\smallsup{2}}-
    \sum_{i=1}^{j}{\hat
Z}_{\lceil i/2\rceil}^{\smallsup{1}}{\hat Z}_{\lfloor i/2\rfloor}^{\smallsup{2}}\Big|>N^{-\eps}\Big)
    =O(N^{-\beta}),
   \end{equation}
\end{corr}

\proof  Bound
\begin{eqnarray}
&&\Big|\frac{\sum_{i=1}^{j}{\cal
    Z}_{\lceil i/2\rceil}^{\smallsup{1}}{\cal Z}_{\lfloor i/2\rfloor}^{\smallsup{2}}}{N}-
    \frac{\sum_{i=1}^{j}{\hat
Z}_{\lceil i/2\rceil}^{\smallsup{1}}{\hat Z}_{\lfloor i/2\rfloor}^{\smallsup{2}}}{N}\Big|\\
&&\qquad \leq
\left|\frac{
    \sum_{i=1}^{j}
    {\cal Z}_{\lfloor i/2\rfloor}^{\smallsup{2}}
    ({\cal Z}_{\lceil i/2\rceil}^{\smallsup{1}}-{\hat Z}_{\lceil i/2\rceil}^{\smallsup{1}})} {\sqrt{N}\sqrt{N}}\right|+
\left|\frac{
    \sum_{i=1}^{j}
    {\hat Z}_{\lceil i/2\rceil}^{\smallsup{1}}
    ({\cal Z}_{\lfloor i/2\rfloor}^{\smallsup{2}}-{\hat Z}_{\lfloor i/2\rfloor}^{\smallsup{2}})} {\sqrt{N}\sqrt{N}}\right|.\nonumber
\label{verg-phony}
\end{eqnarray}
Both terms on the right hand side of (\ref{verg-phony}) can be treated as in the
proof of Proposition \ref{prop-coup}, because the processes with sources (1) and (2)
are independent and uniformly in $i\leq (\frac12+\eta)\log_{\nu} N$,
    $$
    \max\left\{ \frac{\expec[{\cal Z}_i^{\smallsup{2}}],\expec[{\hat Z}_i^{\smallsup{1}}]}{\sqrt{N}}
    \right)=
    \max
    \{
    N^{\eta},N^{\eta}\cdot(1+O(N^{-\alpha_2}))
    ^{(\frac12+\eta)\log_{\nu} N}
    \},
    $$
on $F_{\sN}$. The right-hand side can again be bounded by any small power of $N$ by taking $\eta$
arbitrarily small. We omit further details. \qed

\subsection{Proof of Proposition \ref{prop-bdscoupling1}}
In this second part of the appendix, we restate our main
result on the coupling between the SPG and the BP with
offspring distribution $\{g^{\smallsup{N}}_n\}$ once more and give a full proof.

\begin{prop}
    \label{prop-bdscoupling2}
    There exist $\eta,\beta>0$, $\alpha>\frac 12+\eta$ and a constant $C$,
    such that for all\\ $j\leq (\frac 12 +\eta) \log_\nu N$,
    \begin{equation}
    \prob\Big((1-N^{-\alpha} \nu^{j}) \hat Z_j^{\smallsup{1}}\leq Z_j^{\smallsup{1}}
    \leq (1+N^{-\alpha} \nu^{j}) \hat Z_j^{\smallsup{1}}\Big)\geq 1- C j N^{-\beta}.
    \end{equation}
\end{prop}
\label{proof Prop4.1} This proof is divided into several lemmas.
It is rather involved, and we may think of Proposition
\ref{prop-bdscoupling2} as one of the key estimates of the paper. We start with an
explanation of the different steps in this proof.

The proof of Proposition \ref{prop-bdscoupling2} proceeds by
induction with respect to $j$. Note that for all $j\leq (\frac 12
+\eta) \log_\nu N$, we have $N^{-\alpha} \nu^{j}\leq
N^{(\frac12+\eta)-\alpha} \to 0$, as $N \to \infty$ and
when $\alpha>\eta$. When at level
$j-1$, the event in the statement of the proposition holds, we have
\[
|\hat Z_{j-1}^{\smallsup{1}}- Z_{j-1}^{\smallsup{1}}| \leq
\frac{\nu^{j-1}}{N^{\alpha}} \hat Z_{j-1}^{\smallsup{1}},
\]
so that we control the difference between the number of stubs $
Z_{j-1}^{\smallsup{1}}$ and the number of children $\hat
Z_{j-1}^{\smallsup{1}}$. The absolute value of this difference is
bounded by $ \hat Z_{j-1}^{\smallsup{1}}$ times a fraction that
converges to $0$. For generation $j$ we have to control the
difference $\hat Z_{j}^{\smallsup{1}}- Z_{j}^{\smallsup{1}}$.
Differences in generation $j$ arise from differences in
generation $j-1$ and from drawing stubs with label 2 or label 3.
If a label 2 stub is chosen, then the SPG will contain a loop or cycle
and hence no free stubs in level $j$
are created, whereas in the BP a non-negative number of offspring is attached.
If a label 3 stub is chosen, then the corresponding node with described number of children is attached in the BP,
whereas for the SPG we have to resample until we draw a stub labeled 1 or 2.
Hence, if $Z_{j}^{\smallsup{1}}\geq \hat Z_{j}^{\smallsup{1}}$,
so that the number of free stubs attached to nodes at distance $j-1$ of the
SPG exceeds the number of children in generation $j$ of the BP, then this
overshoot can only be caused by drawing label 3 stubs.
The number of stubs with label 3 is bounded by the total
number drawn in the SPG, i.e., by
$$
\sum_{i=1}^{j-1} Z_{i}^{\smallsup{1}} \le \sum_{i=1}^{j-1}
(1+N^{-\alpha}\nu^i) \hat Z_{i}^{\smallsup{1}} \le
2\sum_{i=1}^{j-1} \hat Z_{i}^{\smallsup{1}}.
$$
For $Z_{j}^{\smallsup{1}}\leq \hat Z_{j}^{\smallsup{1}}$, the number
of stubs with level $2$ or $3$  both matter and their total amount is bounded by
$$
\sum_{i=1}^j Z_{i}^{\smallsup{1}} = \sum_{i=1}^{j-1}
Z_{i}^{\smallsup{1}} +Z_{j}^{\smallsup{1}} \le \sum_{i=1}^{j-1}
(1+N^{-\alpha}\nu^i) \hat Z_{i}^{\smallsup{1}} +\hat
Z_{j}^{\smallsup{1}} \le 2\sum_{i=1}^j \hat Z_{i}^{\smallsup{1}}.
$$
In both cases the probability of drawing a label $2$ or $3$ stub
is bounded by
\eq
\frac{ 2\sum_{i=1}^j \hat Z_{i}^{\smallsup{1}}}{L_{\sN}} \leq
\frac{2N^{\frac12+\delta}}{L_{\sN}},
\en
on the event where $\sum_{i=1}^j \hat Z_{i}^{\smallsup{1}}\le
N^{\frac12+\delta}$. Using that $L_{\sN}$ is of order $\expec[L_{\sN}]=\mu
N$ (see Lemma \ref{lem-mom}), this probability is sufficiently
small to allow us to use Chebychev's inequality.

The main lemmas in this section are Lemma \ref{lem-<} and Lemma
\ref{lem->}. Together, they prove the induction step described
above. Lemmas \ref{lem-firststubs} up to \ref{lem-conti}
are preparations, the most important one being Lemma
\ref{lem-conti}. This lemma shows that if the {\it total} progeny
up to and including generation $j$
of $\{{\hat Z}_i^{\smallsup{1}}\}$ is larger than
$N^{\frac12-\delta}$, for some $\delta>0$, then with overwhelming
probability also each of the sizes of the last two generations, i.e.,
${\hat Z}_{j-1}^{\smallsup{1}}$ and ${\hat
Z}_j^{\smallsup{1}}$, exceed $N^{\frac12-2\delta}$.

As before, we will abbreviate the conditional probability and
expectation given $D_1, \ldots, D_{\sN}$ by $\mathbb P_{\sN}$ and $\mathbb
E_{\sN}$.
\begin{lemma} For $0<\eta<\frac12$ and all $j\geq 1$,
\label{lem-firststubs}
    \begin{equation}
    \prob_{\sN}\big(Z_j^{\smallsup{1}}\neq \hat Z_j^{\smallsup{1}},
    \sum_{i=1}^j \hat Z_i^{\smallsup{1}}\leq N^{\frac 12-\eta}\big)\leq
    \frac{N^{-2\eta}}{(L_{\sN}/N)}, \quad a.s.
    \end{equation}
\end{lemma}
Lemma \ref{lem-firststubs} together with Lemma \ref{lem-mom}
prove Proposition \ref{prop-bdscoupling2} for all $j$ such that
the total size of the BP is at most $N^{\frac
12-\eta}$.

\proof We denote by $l$ the first stub which is grown differently
in the SPG and in the BP. Assume that this $l^{\rm th}$ stub is in
the $j^{\rm th}$ generation or earlier.

Before the growth of the $l^{\rm th}$ stub, the BP and the SPG are
identical. Thus, we must have that $l\leq \sum_{i=1}^j \hat
Z_i^{\smallsup{1}}$. Hence, as we reach to the $l^{\rm th}$ stub,
the number of stubs having either label 2 or 3 is bounded above by
$\sum_{i=1}^j \hat Z_i^{\smallsup{1}}\leq N^{\frac 12-\eta}$. A
difference in the SPG and the BP can only arise when we draw a
stub for the BP having label 2 or 3. Thus, the probability that
the $l^{\rm th}$ stub is the first to create a difference between
the SPG and the BP is bounded above by $N^{\frac 12-\eta}/L_{\sN}$.
Therefore,
    \begin{eqnarray*}
    \prob_{\sN}\big(Z_j^{\smallsup{1}}\neq \hat Z_j^{\smallsup{1}},
    \sum_{i=1}^j \hat Z_i^{\smallsup{1}}\leq N^{\frac 12-\eta}\big)
    \leq \sum_{l=1}^{N^{\frac 12-\eta}}  \frac{N^{\frac 12-\eta}}{L_{\sN}}
    = \frac{N^{-2\eta}}{(L_{\sN}/N)}.
    \end{eqnarray*}
\qed

Recall that $\nu_{\sN}=\sum_{n=0}^{\infty} n g_n^{\smallsup{N}}$
is the expected offspring of the BP $\{\hat Z^{\smallsup{1}}\}_j$
under $\mathbb P_{\sN}$. Note from Proposition \ref{prop-TV2} that
$\nu_{\sN}$ is close to $\nu$ with probability close to one.
In the statement of the next lemma, we write
    \eq
    D_{\sN}^{\smallsup{N}}=\max_{1\leq i\leq N} D_i.
    \en

\begin{lemma}
\label{lem-maxbd} For every $\gamma>0$,
    \begin{equation}
    \prob\big(D_{\sN}^{\smallsup{N}} \geq N^{\gamma}\big)\leq cN^{1-(\tau-1)\gamma}.
    \end{equation}
\end{lemma}

\proof We use Boole's inequality to obtain from
(\ref{distribution}) that
    \begin{equation}
    \prob\big(D_{\sN}^{\smallsup{N}} \geq N^{\gamma}\big)\leq \sum_{i=1}^N \prob(D_i \geq N^{\gamma})
    \leq cN^{1-(\tau-1)\gamma}.
    \end{equation}
\qed

\begin{lemma}
\label{lem-MI} For $\eta, \delta\in (-\frac12,\frac12)$, and all
$j\leq (\frac 12 +\eta) \log_\nu N$, there exists $\beta_2>0$ such
that
    \begin{equation}
    \prob \big(\sum_{i=1}^j \hat Z_i^{\smallsup{1}}\geq N^{\frac 12+\delta}\big)\leq C N^{\eta-\delta}+N^{-\beta_2}.
    \end{equation}
\end{lemma}

\proof By Proposition \ref{prop-TV1}, we can include the indicator
that $|\nu_{\sN}-\nu| \leq N^{-\alpha_2}$; this explains the
additional error term $N^{-\beta_2}$. By the Markov inequality, we
obtain for $j\leq (\frac 12 +\eta) \log_\nu N$,
    \begin{eqnarray*}
    \prob \big(\sum_{i=1}^j \hat Z_i^{\smallsup{1}}\geq N^{\frac 12+\delta}, |\nu_{\sN}-\nu| \leq
    N^{-\alpha_2}\big)
    \leq N^{-\frac 12-\delta}{\mathbb E} \big(\sum_{i=1}^j \hat Z_i^{\smallsup{1}}I[|\nu_{\sN}-\nu| \leq
    N^{-\alpha_2}]\big).\\
    \end{eqnarray*}
The expectation on the right-hand side can be computed by
conditioning:
\begin{eqnarray*}
&&{\mathbb E} \big[\hat Z_i^{\smallsup{1}}I[|\nu_{\sN}-\nu| \leq
N^{-\alpha_2}]\big]=
\expec[\expec_{\sN} \big[\hat Z_i^{\smallsup{1}}I[|\nu_{\sN}-\nu| \leq N^{-\alpha_2}]\big]\\
&&\qquad=\expec[I[|\nu_{\sN}-\nu| \leq N^{-\alpha_2}]\cdot D_1
\nu_{\sN}^{i-1}]\leq (\nu+N^{-\alpha_2})^{i-1}\expec[D_1].
\end{eqnarray*}
Hence,
\begin{eqnarray*}
\prob \big(\sum_{i=1}^j \hat Z_i^{\smallsup{1}}\geq N^{\frac
12+\delta}\big)&\leq&
 N^{-\beta_2}+\mu N^{-\frac12-\delta}\sum_{i=1}^{j}(\nu+N^{-\alpha_2})^{i-1}\\
&\leq& N^{-\beta_2}+\mu
N^{-\frac12-\delta}\frac{(\nu+N^{-\alpha_2})^{j}-1}{(\nu+N^{-\alpha_2})-1}
\leq N^{-\beta_2}+CN^{\eta-\delta}.
\end{eqnarray*}
\qed


In the lemma below, we write $d$ for a random variable with
discrete distribution $\{g_n^{\smallsup{N}}\}$ given    in (\ref{gnN}), and
$\mbox{Var}_{\sN}(d)$ for the variance of $d$ under $\prob_{\sN}$.
Furthermore, we let, for any $0<a <\frac12$,
    $$
    A_{\sN}=A_{\sN}(a, \gamma,\alpha_2)=\left\{\left|\frac{L_{\sN}}{\mu N}-1\right|\leq
    N^{-a}\right\}\cap \{D_{\sN}^{\smallsup{N}} \leq N^{\gamma}\}
    \cap \{|\nu_{\sN}-\nu|\leq N^{-\alpha_2}\},
    $$
then, according to Proposition \ref{prop-TV1}, Lemmas \ref{lem-mom} and \ref{lem-maxbd}, we have
    \eq
    \label{Acbd}
    \prob(A_{\sN}^c)=O(N^{-\epsilon}),
    \en
where $\epsilon=b\wedge ((\tau-1)\gamma-1)\wedge \beta_2>0$ whenever $\gamma>1/(\tau-1)$.
On $A_{\sN}$, we have
    \begin{equation}
    \label{quotient LNandN} \frac{1}{\mu(1+N^{-a})}\leq
    \frac{N}{L_{\sN}}\leq \frac{1}{\mu(1-N^{-a})}.
    \end{equation}
This will be used in the following lemma.
\begin{lemma}
\label{lem-varbd} For every $\gamma>0$,
    \begin{equation}
    \expec \big({\rm Var}_{\sN}(d) I[A_{\sN}]\big)
    \leq C N^{(4-\tau)^+\gamma},
    \end{equation}
where $x^+=\max(0,x)$.
\end{lemma}

\proof Since the variance of a random variable is bounded by its
second moment,
    $$
    \mbox{Var}_{\sN}(d)\leq \sum_{n=0}^\infty n^2 g_n^{\smallsup{N}}=
    \sum_{n=0}^\infty\sum_{j=1}^N \frac{n^2(n+1)}{L_{\sN}}I[D_j=n+1]
    \leq \frac{1}{L_{\sN}}\sum_{j=1}^N D_j^3,
    $$
and so, for $\tau \in (3,4]$,
    \[
    \expec \big(\mbox{Var}_{\sN}(d)
    I[A_{\sN}]\big)\leq \sum_{j=1}^N\expec\big[\frac{1}{L_{\sN}} D_j^3 I[A_{\sN}]\big]
    \leq \frac{N}{\mu N}\expec\big[D^3 I[D \leq
    N^{\gamma}] \big] \leq C\sum_{i=1}^{\lceil N^{\gamma}\rceil }i^3 f_i
    \leq N^{\gamma(4-\tau)},
    \]
by Lemma \ref{lem-boundsf}. For $\tau>4$, the third moment of $D$ is finite, and
the result is also true even without the indicator $I[D_{\sN}^{\smallsup{N}} \leq N^{\gamma}]$. \qed

\begin{lemma}
\label{lem-conti} For all $(\frac 12 -2\eta) \log_\nu N\leq j\leq
(\frac 12 +2\eta) \log_\nu N$, there exists $\delta,\beta>0$ such
that
    \begin{eqnarray}
    &&\prob \big(\sum_{i=1}^{j} \hat Z_i^{\smallsup{1}}\geq N^{\frac 12-\delta},
    \hat Z_{j-1}^{\smallsup{1}}\leq N^{\frac 12 -2\delta}\big)\leq C N^{-\beta},
    \label{conti}\\
    &&\prob \big(\sum_{i=1}^{j} \hat Z_i^{\smallsup{1}}\geq N^{\frac 12-\delta},
    \hat Z_{j}^{\smallsup{1}}\leq N^{\frac 12 -2\delta}\big)\leq C N^{-\beta}.
    \end{eqnarray}
\end{lemma}
{\bf Remark:} The statements of the lemma are almost identical,
the difference being that the index of $\hat
Z_{j-1}^{\smallsup{1}}$ in the first statement is replaced by the
index $j$ in the second statement. We will be satisfied with a
proof for the first statement only, the proof with index $j$ is a
straightforward extension.

\proof Since $\sum_{i=1}^{j} \hat Z_i^{\smallsup{1}}\geq N^{\frac
12-\delta}$, there must be an $i\leq j\leq (\frac 12+2\eta)\log_\nu N$
such that for $N$ large
enough
    \[\hat Z_i^{\smallsup{1}}\geq N^{\frac 12-\delta}/j
    \geq
    \frac{N^{\frac 12-\delta}}{(\frac 12 +2\eta) \log_\nu N}
    \geq N^{\frac 12-\frac 32 \delta}.
    \]
We write $I$ for the first $i\leq j$ such that $\hat
Z_i^{\smallsup{1}}\geq N^{\frac 12-\frac 32 \delta}$. It suffices to
bound
    \eq
    \sum_{i=1}^j
    \prob\big(\sum_{k=1}^{j} \hat Z_k^{\smallsup{1}}\geq N^{\frac 12-\delta}, I=i,
    \hat Z_{j-1}^{\smallsup{1}}\leq N^{\frac 12 - 2\delta}\big).
    \en
The contribution from $I=j-1$ is 0. When $I=j$, then $\hat
Z_j^{\smallsup{1}}\geq N^{\frac 12-\frac 32\delta}$, but $\hat
Z_{j-1}^{\smallsup{1}}\leq N^{\frac 12-2\delta}$ so that from
the Markov inequality
    \begin{eqnarray}
    &&\prob\big(\sum_{k=1}^{j} \hat Z_k^{\smallsup{1}}\geq N^{\frac 12-\delta}, I=j,
    \hat Z_{j-1}^{\smallsup{1}}\leq N^{\frac 12 -2\delta}\big)\\
    &&\qquad \leq \expec\Big[I[\hat Z_{j-1}^{\smallsup{1}}\leq N^{\frac 12-2\delta}]
    \prob_{\sN}\big(\hat Z_j^{\smallsup{1}}\geq N^{\frac 12-\frac 32\delta}\big|
    \hat Z_{j-1}^{\smallsup{1}}\big)\Big]\nonumber\\
    &&\qquad \leq N^{-\frac 12+\frac 32\delta}\expec
    \Big[I[\hat Z_{j-1}^{\smallsup{1}}\leq N^{\frac 12-2\delta}]\expec_{\sN}\big[\hat Z_j^{\smallsup{1}}\big|
    \hat Z_{j-1}^{\smallsup{1}}\big]\Big]\nonumber\\
    &&\qquad = N^{-\frac 12+\frac 32\delta}\expec
    \Big[I[\hat Z_{j-1}^{\smallsup{1}}\leq N^{\frac 12-2\delta}]\nu_{\sN} \hat Z_{j-1}^{\smallsup{1}}\Big]
    \leq C N^{-\frac 12+\frac 32\delta} N^{\frac 12-2\delta}=C N^{-\delta/2}.\nonumber
    \end{eqnarray}
Thus, we are left to deal with the cases where $I<j-1$. Then,
there exists an $i<j-1$ such that $\hat Z_i^{\smallsup{1}}\geq
N^{\frac 12-\frac 32\delta}$, but $\hat Z_{j-1}^{\smallsup{1}}\leq
N^{\frac 12-2\delta}$. Thus, there must be a first $s\geq i$ such
that $\hat Z_{s+1}^{\smallsup{1}}\leq \hat Z_{s}^{\smallsup{1}}$.
Consequently, $\hat Z_s^{\smallsup{1}} \geq \hat
Z_i^{\smallsup{1}}\geq N^{\frac 12-\frac 32\delta}$. We will bound,
uniformly in $s$,
    \eq
    \prob(\hat Z_{s+1}^{\smallsup{1}}\leq \hat Z_{s}^{\smallsup{1}}, \hat Z_s^{\smallsup{1}}\geq
    N^{\frac 12-\frac 32\delta}\big)\leq N^{-\beta},
    \en
for some $\beta>0$. This proves (\ref{conti}), since the total number of possible $i$ and $s$
with $i\leq s\leq j$ is bounded by $(\log_{\nu}{N})^2$.

We use Lemma \ref{lem-maxbd} to see that we may include the
indicator on $A_{\sN}$ for any
$\gamma>1/(\tau-1)$. We will use the Chebychev inequality and Lemma
\ref{lem-varbd} to obtain that
    \begin{eqnarray}
    &&\prob(\hat Z_{s+1}^{\smallsup{1}}\leq \hat Z_{s}^{\smallsup{1}}, \hat Z_s^{\smallsup{1}}\geq
    N^{\frac 12-\frac 32\delta}, A_{\sN}\big)\\
    &&\qquad =\expec\Big[I[\hat Z_s^{\smallsup{1}}\geq  N^{\frac 12-\frac 32\delta},
    A_{\sN}]
    \prob_{\sN}\big(\hat Z_{s+1}^{\smallsup{1}}\leq \hat Z_{s}^{\smallsup{1}}\big|
    \hat Z_s^{\smallsup{1}}\big)\Big]\nonumber\\
    &&\qquad \leq \expec\Big[I[\hat Z_s^{\smallsup{1}}\geq  N^{\frac 12-\frac 32\delta},
    A_{\sN}]
    \prob_{\sN}\big(\big|\hat Z_{s+1}^{\smallsup{1}}-\nu_{\sN}\hat Z_{s}^{\smallsup{1}}\big|\geq (\nu_{\sN}-1)
    \hat Z_{s}^{\smallsup{1}}\big|
    \hat Z_s^{\smallsup{1}}\big)\Big]\nonumber\\
    &&\qquad \leq \expec\Big[I[\hat Z_s^{\smallsup{1}}\geq  N^{\frac 12-\frac 32\delta},
    A_{\sN}]
    (\nu_{\sN}-1)^{-2} \frac{\mbox{Var}_{\sN}(d_1)}{\hat Z_s^{\smallsup{1}}}\Big]\nonumber\\
    &&\qquad \leq C N^{(4-\tau)^+ \gamma -\frac 12+\frac 32\delta}\leq
    N^{-\beta},\nonumber
    \end{eqnarray}
with $C=2(\nu-1)^{-2}$, and since $(4-\tau)^+\gamma<1/2$ and
$\delta>0$ can be taken arbitrarily small. \qed

\medskip
\noindent
We are now ready to give the proof of Proposition \ref{prop-bdscoupling2}.

\noindent{\bf Proof of Proposition \ref{prop-bdscoupling2}}.

\noindent We first set the stage for the proof by induction in $j$. Fix
$\eta<\delta<2\eta$, and $\alpha>\frac12+\eta$, and define
    \eq
    E_j=\big\{\forall i\leq j: (1-N^{-\alpha} \nu^{i}) \hat Z_i^{\smallsup{1}}\leq Z_i^{\smallsup{1}}
    \leq (1+N^{-\alpha} \nu^{i}) \hat Z_i^{\smallsup{1}}\big\}.
    \en
We will prove by induction  that for all $j \leq
(\frac12+\eta)\log_{\nu} N$,
    \begin{equation}
    \prob\big(E_j^c\big)\leq C j N^{-\beta},
    \label{todo}
    \end{equation}
which implies  Proposition \ref{prop-bdscoupling2} by taking the
complementary event. First, by Lemma \ref{lem-firststubs} and \ref{lem-MI} and since
$\eta<\delta$ we see that it is sufficient to prove for $j \leq
(\frac12+\eta)\log_{\nu} N$,
$$
\prob\big(E_j^c, N^{\frac 12-\delta} \leq \sum_{i=1}^j \hat
Z_i^{\smallsup{1}}\leq N^{\frac 12+\delta}\big)\leq C j
N^{-\beta}.
$$

For $j<(\frac 12-2\eta)\log_\nu N$, we bound
    \begin{eqnarray}
    \prob\Big(E_j^c, N^{\frac 12-\delta} \leq \sum_{i=1}^j
    \hat  Z_i^{\smallsup{1}}\leq N^{\frac 12+\delta}\Big)
    \leq \prob\Big(\sum_{i=1}^j \hat  Z_i^{\smallsup{1}}\ge N^{\frac 12-\delta}\Big)
    \leq        N^{-\beta}+CN^{-2\eta+\delta}\nonumber,
    \end{eqnarray}
by the Markov inequality and using Proposition \ref{prop-TV1} in a similar
way as in Lemma \ref{lem-MI}. Hence, the statement in (\ref{todo})
follows for $j<(\frac12-2\eta)\log_\nu N$. This initializes the
induction in $j$.

To advance the induction, we bound
    \begin{eqnarray*}
    \prob\big(E_j^c, N^{\frac 12-\delta} \leq
    \sum_{i=1}^j \hat Z_i^{\smallsup{1}}\leq N^{\frac 12+\delta}\big)
    &\leq& \prob(E_{j-1}^c) + \prob\big(E_j^c\cap E_{j-1}, N^{\frac 12-\delta} \leq
    \sum_{i=1}^j \hat Z_i^{\smallsup{1}}\leq N^{\frac 12+\delta}\big)\\
    &\leq& C (j-1) N^{-\beta}+\prob\big(E_j^c\cap E_{j-1}, N^{\frac 12-\delta} \leq
    \sum_{i=1}^j \hat Z_i^{\smallsup{1}}\leq N^{\frac 12+\delta}\big),
    \end{eqnarray*}
where the last inequality follows by the induction hypothesis.
Thus, it suffices to prove that
    \eq
    \prob\big(E_j^c\cap E_{j-1}'\big)\leq CN^{-\beta},
    \label{todo2}
    \en
where
    \[E_{j-1}'=E_{j-1}\cap \{N^{\frac 12-\delta} \leq
    \sum_{i=1}^j \hat Z_i^{\smallsup{1}}\leq N^{\frac 12+\delta}\}.\]
Note that
    \eq
    \label{union}
    E_j^c\cap E_{j-1}'= \Big(\big\{Z_j^{\smallsup{1}}<(1-N^{-\alpha} \nu^{j})\hat Z_j^{\smallsup{1}}\big\}
    \cap E_{j-1}'\Big)\bigcup
    \Big(\big\{Z_j^{\smallsup{1}}>(1+N^{-\alpha} \nu^{j})\hat Z_j^{\smallsup{1}}\big\}
    \cap E_{j-1}'\Big).
    \en
We write the disjoint events on the right-hand side of
(\ref{union}) as $E_{j,<}^c$ and $E_{j,>}^c$  and bound the probability
of these
events separately. We will start with $E_{j,<}^c$. This result is
stated in the following lemma:

\begin{lemma}
\label{lem-<} There exists $\beta>0$ such that for all $(\frac
12 -2\eta) \log_\nu N< j\leq (\frac 12 +\eta) \log_\nu N$,
    \begin{equation}
     \prob(E_{j,<}^c)\leq C N^{-\beta}.
    \label{<bd}
    \end{equation}
\end{lemma}

\proof We note that on $E_{j,<}^c$, we have that
    \[
    \sum_{i=1}^j Z_i^{\smallsup{1}}
    \leq \sum_{i=1}^j (1+\nu^i N^{-\alpha})\hat Z_i^{\smallsup{1}}\leq
    (1+N^{\frac12+\eta}N^{-\alpha})\sum_{i=1}^j\hat Z_i^{\smallsup{1}}
    \leq 2N^{\frac 12+\delta},
    \]
because $\alpha > \frac12+\eta$. Thus, for every stub which is
grown simultaneously for the BP and the SPG, there is a
probability bounded from above by $2N^{\frac 12+\delta}/L_{\sN} $ that
a difference is created between the BP and the SPG (such a
difference is called a miscoupling). Denote by $U$ the number of
stubs where such a difference occurs. Then, $U$ is bounded from
above by a binomial random variable with $n=N^{\frac 12+\delta}$
and $p=2N^{\frac 12+\delta}/L_{\sN}$. Thus, by the Markov inequality,
we have,
    $$
    \prob_{\sN}(U\geq N^a)\leq \frac{2 N^{-a+1+2\delta}}{L_{\sN}}.
    $$
Using (\ref{quotient LNandN}), we obtain, for $2\delta<a$,
    \eq
    \label{U-bd}
    \prob(U\geq N^a)\leq C N^{-a+2\delta}+N^{-b}\le N^{-\beta}.
    \en

Observe that differences between $Z_j^{\smallsup{1}}$ and $\hat
Z_j^{\smallsup{1}}$ can only arise through (i) different numbers
of stubs in the $(j-1)^{\rm st}$ generation, and (ii) differences
created in the $j^{\rm th}$ generation which we previously called
miscouplings. In the first case, the difference in the number of
stubs is bounded from below by an independent draw from
$g^{\smallsup{N}}$. A miscoupling occurs if we draw a stub with
label $2$ or $3$. Hence,
    \eq
    Z_j^{\smallsup{1}}-\hat Z_j^{\smallsup{1}}
    \geq -\sum_{i=1}^{\big(\hat Z_{j-1}^{\smallsup{1}}-Z_{j-1}^{\smallsup{1}}\big)^+} d_i
    -\sum_{i=1}^U \tilde d_i,
    \en
where $\{d_i\}_{i\geq 1}$ are independent draws from
$g^{\smallsup{N}}$ and $\{\tilde d_i\}_{i\geq 1}$ are draws
conditionally on drawing a stub labeled 2 or 3.  On $E_{j-1}'$,
we have that
    \eq
    \big(\hat Z_{j-1}^{\smallsup{1}}-Z_{j-1}^{\smallsup{1}}\big)^+
    \leq N^{-\alpha} \nu^{j-1} \hat Z_{j-1}^{\smallsup{1}},
    \en
so that on $E_{j<}^c$, introducing the notation
$\alpha_{\sN,j}=N^{-\alpha} \nu^{j-1} \hat Z_{j-1}^{\smallsup{1}}$,
        \eq
    \sum_{i=1}^{\alpha_{\sN,j}} d_i
    +\sum_{i=1}^{U} \tilde d_i
    \geq \sum_{i=1}^{\big(\hat Z_{j-1}^{\smallsup{1}}-Z_{j-1}^{\smallsup{1}}\big)^+} d_i
    +\sum_{i=1}^U \tilde d_i> N^{-\alpha} \nu^{j}\hat Z_j^{\smallsup{1}}.
    \en

Combining this with (\ref{U-bd}) and using the definition of
$\alpha_{\sN,j}$, we see that in order to prove (\ref{<bd}) it
suffices to show that

    \eq \label{replacant} \prob\Big(\Big\{ \sum_{i=1}^{\alpha_{\sN,j}}
    d_i
    +\sum_{i=1}^{\lceil N^a\rceil } \tilde d_i >N^{-\alpha}\nu^j\hat Z_j^{\smallsup{1}}\Big\}\cap E_{j-1}'
    \Big)\leq CN^{-\beta}.
\en We will first show that on $E_{j-1}'$ the term
$\sum_{i=1}^{\lceil N^a \rceil} \tilde d_i$ is small compared to
$N^{-\alpha}\nu^j\hat Z_j^{\smallsup{1}}$, if we choose $a$
sufficiently small. On $E_{j-1}'$, we have $\sum_{i=1}^j
Z_{i}^{\smallsup{1}}\geq N^{\frac12-\delta}$, and so, with
probability larger than $1-CN^{-\beta}$, according to Lemma
\ref{lem-conti}, we have that also $Z_{j}^{\smallsup{1}}\geq
N^{\frac12-2\delta}$. Hence,
\begin{eqnarray*}
\prob\Big(\Big\{ \sum_{i=1}^{\lceil N^a\rceil} \tilde d_i >\frac 12
N^{-\alpha}\nu^j\hat Z_j^{\smallsup{1}}\Big\}\cap E_{j-1}'
    \Big)&\leq& CN^{-\beta}+\prob\Big(
\sum_{i=1}^{\lceil N^a\rceil} \tilde d_i >\frac 12 N^{-\alpha}\nu^j N^{\frac12-2\delta}\Big)\\
&\leq& CN^{-\beta}+\prob\Big( \sum_{i=1}^{\lceil N^a\rceil} \tilde d_i >\frac 12
N^{\frac12-2\eta-\alpha} N^{\frac12-2\delta}
    \Big)\\
    &\leq& CN^{-\beta}+\prob\Big(
N^a D_{\sN}^{\smallsup{N}} >\frac 12 N^{1-2\eta-\alpha-2\delta}\Big)\\
 &\leq&
 CN^{-\beta} +cN^{1-(\tau-1)\gamma},
 \end{eqnarray*}
where $\gamma=1-2\eta-\alpha-2\delta-a<\frac12$, but can be taken
arbitrary close to $\frac12$. Since $\tau>3$, we then have that
$cN^{1-(\tau-1)\gamma}<N^{-\beta}$.

Hence it suffices to prove the statement in (\ref{replacant})
without the term $\sum_{i=1}^{\lceil N^a \rceil} \tilde d_i $, that is, it
suffices to prove
    \eq \label{replacant2} \prob\Big(
    \sum_{i=1}^{\alpha_{\sN,j}} d_i
    >\frac 12 N^{-\alpha}\nu^j\hat Z_j^{\smallsup{1}}, E_{j-1}'
    \Big)\leq CN^{-\beta}.
    \en
Since we can write $\hat
Z_{j}^{\smallsup{1}}=\sum_{i=1}^{\hat Z_{j-1}^{\smallsup{1}}}
d_i$, and, using again Lemma  \ref{lem-conti}, we have that
$E_{j-1}'$ implies $\hat Z_{j-1}^{\smallsup{1}}\geq
N^{\frac12-2\delta}$, with probability larger than
$1-CN^{-\beta}$, it is sufficient to prove that
    \eq
    \label{replacant3} \prob\Big( (1-N^{-\alpha}
    \nu^{j})\sum_{i=1}^{\alpha_{N,j}} d_i
    >  \frac 12 N^{-\alpha} \nu^{j}\sum_{i=\alpha_{N,j}+1}^{\hat Z_{j-1}^{\smallsup{1}}} d_i,
    Z_{j-1}^{\smallsup{1}}\geq N^{\frac12-2\delta}\Big)\le
    N^{-\beta}.
    \en
Now $\expec_{\sN}[d]=\nu_{\sN}$ and, given $\hat
Z_{j-1}^{\smallsup{1}}$, the variance of $\sum_{i=1}^{\hat
Z_{j-1}^{\smallsup{1}}} (d_i-\nu_{\sN})$ equals $\hat
Z_{j-1}^{\smallsup{1}}\mbox{Var}_{\sN}(d)$. Therefore, by the Chebychev
inequality,
    \begin{eqnarray*}
    && \prob_{\sN} \Big( (1-N^{-\alpha} \nu^{j})\sum_{i=1}^{\alpha_{\sN,j}}
    d_i
    >  \frac 12 N^{-\alpha} \nu^{j}\sum_{i=\alpha_{\sN,j}+1}^{\hat Z_{j-1}^{\smallsup{1}}} d_i
    \Big|\hat Z_{j-1}^{\smallsup{1}}
        \Big)\\
    &&\qquad \leq \prob_{\sN} \Big( (1-N^{-\alpha}
    \nu^{j})\sum_{i=1}^{\alpha_{N,j}} (d_i-\nu_{\sN}) -  N^{-\alpha}
    \nu^{j}\sum_{i=\alpha_{N,j}+1}^{\hat Z_{j-1}^{\smallsup{1}}}
    (d_i-\nu_{\sN})
    >\frac 12 \nu_{\sN} \alpha_{N,j}(\nu-1)
    \Big|\hat Z_{j-1}^{\smallsup{1}}
        \Big)\\
    &&\qquad
        \leq  \frac{4\hat Z_{j-1}^{\smallsup{1}}\mbox{Var}_{\sN}(d_1)}
        {(\nu_{\sN} \alpha_{\sN,j}(\nu-1))^2}
        = \frac{4\mbox{Var}_{\sN}(d_1)}
        {\hat Z_{j-1}^{\smallsup{1}}
        N^{-2\alpha}
        \nu^{2j}
        (1-\nu^{-1})^2\nu_{\sN}^2}.
    \end{eqnarray*}
We use Lemma \ref{lem-varbd}.
Hence, by intersecting with the event $I[A_{\sN}]$ and its complement,
and using (\ref{Acbd}), we obtain for $j\geq (\frac12-2\eta)\log_\nu {N}$,
        \begin{eqnarray*}
        &&\prob\Big(
        (1-N^{-\alpha} \nu^{j})\sum_{i=1}^{\alpha_{N,j}} d_i
        >  N^{-\alpha} \nu^{j}\sum_{i=\alpha_{N,j}+1}^{\hat Z_{j-1}^{\smallsup{1}}} d_i,
        Z_{j-1}^{\smallsup{1}}\geq N^{\frac12-2\delta}\Big)\\
        &&\quad\leq  c_1 N^{-\epsilon}
        +c_2\frac{\expec\Big[\mbox{Var}_{\sN}(d_1)I[A_{\sN}]\Big]}
        {N^{\frac12-2\delta}N^{-2\alpha}\nu^{2j}}\\
        &&\quad
        \leq  c_1 N^{-\epsilon}+c_2N^{2\alpha+2\delta-\frac12-1+4\eta } N^{(4-\tau)^+\gamma}
        \leq c_1 N^{-\epsilon}+ c_2N^{-\beta}\leq CN^{-\beta},
        \end{eqnarray*}
by fixing $\alpha>\frac 12+\eta$ so that the exponent is negative
(using that $\gamma<\frac 12$ and $(4-\tau)^+\leq 1$), and writing
$\beta=\frac32-2\alpha-2\delta-4\eta-(4-\tau)^+\gamma>0$.
This proves (\ref{replacant3}) and completes the proof of Lemma
\ref{lem-<}.    \qed

\medskip

Before turning to the proof of the bound on $\prob(E_{j,>}^c)$ in
Lemma \ref{lem->} below, we start with a preparatory lemma and
some definitions.
Suppose we have $L$ objects divided into $N$ groups of sizes $d_1,
\ldots, d_{\sN}$, so that $L=\sum_{i=1}^N d_i$. Suppose we draw an
object at random, and we define a random variable by $d_I-1$ when
the object is taken from the $I^{\rm th}$ group. This gives a
distribution $g^{\smallsup{\vec{d}}}$, i.e.,
    \eq
    g^{\smallsup{\vec{d}}}_n = \frac{1}{L} \sum_{i=1}^N d_iI[d_i=n+1].
    \en
Clearly, $g^{\smallsup{N}}=g^{\smallsup{\vec{D}}}$, where
$\vec{D}=(D_1, \ldots, D_{\sN})$.

We next label $M$ of the $L$ objects, and suppose that the
distribution $g^{\smallsup{\vec{d}}}(M)$ is obtained in a similar
way from drawing conditionally on drawing an unlabelled object.
More precisely, we remove the labelled objects from all objects
thus creating new $d_1', \ldots, d_{\sN}',\, \sum d_i'=L-M$, and we
let $g^{\smallsup{\vec{d}}}(M)=g^{\smallsup{\vec{d}'}}.$ Even
though this is not indicated, the law $g^{\smallsup{\vec{d}}}(M)$
depends on what objects have been labelled.

Lemma \ref{lem-stochbds} below shows that the law
$g^{\smallsup{\vec{d}}}(M)$ can be bounded above and below by two
specific ways of labeling the $M$  objects. Before we can state
the lemma, we need to describe those specific labellings.

For a vector $\vec{d}$, we let $d_{\smallsup{1}}, \ldots,
d_{\smallsup{N}}$ be the ordered vector, so that
$d_{\smallsup{1}}=\min_{i=1,\ldots, N} d_i$ and
$d_{\smallsup{N}}=\max_{i=1,\ldots, N} d_i$. Then the laws
$f^{\smallsup{\vec{d}}}(M)$ and $h^{\smallsup{\vec{d}}}(M)$,
respectively, are defined by successively decreasing
$d_{\smallsup{N}}$ and $d_{\smallsup{1}}$ respectively, by one.
Thus,
    \begin{eqnarray}
    f^{\smallsup{\vec{d}}}_n(1)&=&\frac{1}{L-1} \sum_{i=1}^{N-1}
    d_{\smallsup{i}}I[d_{\smallsup{i}}=n+1] +
    \frac{(d_{\smallsup{N}}-1)I[d_{\smallsup{N}}-1=n+1]}{L-1}\label{f1}\\
    h^{\smallsup{\vec{d}}}_n(1)&=&\frac{1}{L-1} \sum_{i=2}^{N}
    d_{\smallsup{i}}I[d_{\smallsup{i}}=n+1] +
    \frac{(d_{\smallsup{1}}-1)I[d_{\smallsup{1}}-1=n+1]}{L-1}.\label{h1}
    \end{eqnarray}
For $f^{\smallsup{\vec{d}}}(M)$ and $h^{\smallsup{\vec{d}}}(M)$,
respectively, we repeat the above change $M$ times. Here we note
that when $d_{\smallsup{1}}=1$, and for
$h^{\smallsup{\vec{d}}}(1)$ we decrease it by one, that we only
keep the $d_i\geq 1$. Thus, in this case, the number of groups of
objects is decreased by 1.


Finally, we write that $f \preceq g$ when the distribution $f$ is
stochastically dominated by $g$, i.e., when $\sum_{i=0}^n f_i
\geq \sum_{i=0}^n g_i$ for
all $ n\geq 0$. Similarly, we write that $X \preceq Y$ when for the probability
mass functions $f_X, f_Y$ we have that $f_X \preceq f_Y$.

We next prove stochastic bounds on the distribution
$g^{\smallsup{\vec{d}}}(M)$ that are uniform in the choice of the
$M$ labelled objects.

\begin{lemma}
\label{lem-stochbds} For all choices of $M$ labelled objects
    \begin{equation}
    f^{\smallsup{\vec{d}}}(M)\preceq g^{\smallsup{\vec{d}}}(M)
    \preceq h^{\smallsup{\vec{d}}}(M).
    \label{bdsdistr}
    \end{equation}
Thus, the expectation and variance of the random variable $X(M)$
with probability mass function $g^{\smallsup{\vec{d}}}(M)$ are
bounded by
    \eq
    \expec[X(M)]\leq \expec[\overline{X}(M)],
    \qquad {\rm Var}[X(M)]\leq \expec[\overline{X}(M)^2],
    \label{boundsmoments}
    \en
where $\overline{X}(M)$ has probability mass function
$h^{\smallsup{\vec{d}}}(M)$.

Moreover, when $X_1, \ldots, X_l$ are draws from
$g^{\smallsup{\vec{d}}}(M_1), \ldots,
g^{\smallsup{\vec{d}}}(M_l)$, where the only dependence between
the $X_i$ resides in the labelled objects, then
    \eq
    \sum_{i=1}^l \underline{X}_i \preceq \sum_{i=1}^l X_i \preceq \sum_{i=1}^l \overline{X}_i,
    \label{bdssums}
    \en
where $\{\underline{X}_i\}_{i=1}^l$ and
$\{\overline{X}_i\}_{i=1}^l$, respectively, are i.i.d.\ copies of
$\underline{X}$ and $\overline{X}$ with laws
$f^{\smallsup{\vec{d}}}(M)$ and $h^{\smallsup{\vec{d}}}(M)$ for
$M=\max_{i=1}^l M_i$, respectively.
\end{lemma}
In the proof of Proposition \ref{prop-bdscoupling2}, we will only
use the upper bounds in Lemma \ref{lem-stochbds}.

\proof In order to prove (\ref{bdsdistr}), we will use induction
in $M$. We note that $f^{\smallsup{\vec{d}}}(0)=
g^{\smallsup{\vec{d}}}(0)=
h^{\smallsup{\vec{d}}}(0)=g^{\smallsup{\vec{d}}}$, and this
initializes the induction. To advance the induction, we note that
we need to investigate the effect of labelling one extra object.
For $f^{\smallsup{\vec{d}}}(M)$, we need to maximize the
cumulative distribution function, whereas for
$h^{\smallsup{\vec{d}}}(M)$, we need to minimize it. Clearly,
(\ref{f1}-\ref{h1}) are optimal. This advances the induction. The
statement in (\ref{boundsmoments}) follows from  (\ref{bdsdistr})

To prove (\ref{bdssums}), we see that for every $j$, conditionally
on the `past' $(X_1, \ldots, X_{j-1})$, the random variable $X_j$
is stochastically bounded by $\underline{X}_j$ and
$\overline{X}_j$, respectively. This completes the proof
of Lemma \ref{lem-stochbds}.\qed

\begin{lemma}
\label{lem->} There exists $\beta>0$ such that for all $j\leq
(\frac 12 +\eta) \log_\nu N$,
    \begin{equation}
    \prob(E_{j,>}^c)\leq C N^{-\beta}.
    \end{equation}
\end{lemma}

\proof The proof of Lemma \ref{lem->} follows the proof of Lemma
\ref{lem-<}, and we focus on the differences only.

Let $V$ denote the number of stubs out of the $\hat
Z_{j-1}^{\smallsup{1}}$ stubs that are attached to stubs with
label 3 in the BP. Since for each stub in the $(j-1)^{\rm st}$
generation, on $E'_{j-1}$, we have that there are at most
$2\sum_{i=1}^{j-1} Z_i^{\smallsup{1}}\leq 2 N^{\frac 12+\delta}$
stubs with label 3, we have that $V$ is bounded from above by a
binomial random variable with $n=N^{\frac 12+\delta}$ and $p=2
N^{\frac 12+\delta}/L_{\sN}$. Thus, by the Markov inequality, we have
that for any $a>2\delta$,
    \eq
    \prob(V\geq N^a)\leq C N^{-\beta},\qquad \mbox{with}\qquad \beta=a-2\delta>0,
    \en
where we can take $a$ arbitrarily small by choosing $\delta>0$
small.

We thus assume that $V\leq N^a$. We next proceed by investigating
$\prob(E_{j,>}^c).$ Now, on $E_{j,>}^c\cap E_{j-1}$, we have that
    \eq
    Z_j^{\smallsup{1}}> (1+N^{-\alpha} \nu^{j})\hat Z_j^{\smallsup{1}}.
    \en
Thus,  $Z_j^{\smallsup{1}}$ is larger than $\hat Z_j^{\smallsup{1}}$.
We note that $Z_j^{\smallsup{1}}$ can only become larger than
$\hat Z_j^{\smallsup{1}}$ from (a) a redraw and the redraw exceeds
the original draw from
$g^{\smallsup{N}}$; and (b) stubs in $Z_{j-1}^{\smallsup{1}}$ that
are not in $\hat Z_{j-1}^{\smallsup{1}}$ which give rise to new stubs.
On $E_{j-1}$, we thus have that (recalling that
$\alpha_{N,j}=N^{-\alpha} \nu^{j-1}\hat Z_{j-1}^{\smallsup{1}}$)
    \eq
    Z_j^{\smallsup{1}}-\hat Z_j^{\smallsup{1}}\leq \sum_{i=1}^{\alpha_{N,j}}
    d_i' +\sum_{i=1}^V d_i'',
    \en
where $d_i', d_i''$ are drawn from the appropriate conditional
distributions given that we pick a stub with label unequal to 3.

We note that each of the $d_i', d_i''$ is obtained by drawing from
stubs conditionally on labels not being 3. Since the total number
of stubs labeled 3 is throughout the growth process bounded above
by $2\sum_{i=1}^{j-1} Z_i^{\smallsup{1}}\leq 2 N^{\frac
12+\delta}$, on $V\leq N^a$, we obtain that by Lemma
\ref{lem-stochbds}, $\{d'_i\}_{i=1}^{\alpha_{\sN,j}}$ and $\{d''_i\}_{i=1}^V$
are bounded above by $\alpha_{N,j}+\lceil N^a\rceil$ independent copies of
$\overline{X}_i(2 N^{\frac 12+\delta})$, where for any $M$,
$\overline{X}_i(M)$ has probability distribution
$h^{\smallsup{\vec D}}(M)$.

We note that by (\ref{boundsmoments}) and Proposition
\ref{prop-TV1}, the expectation of $\overline{X}_i(2 N^{\frac
12+\delta})$ is bounded above by $\nu+N^{-\alpha_2}$ for some
$\alpha_2>0$, and the variance of $\overline{X}_i(2 N^{\frac
12+\delta})$ obeys the same bound as $\mbox{Var}_{\sN}(d)$ in Lemma
\ref{lem-varbd}. Thus, we can copy the remaining part of the proof
from the proof of Lemma \ref{lem-<}. \qed

\subsection{Proof of Proposition \ref{prop-caft}}
\label{sec-pf33}
In this section, we prove Proposition \ref{prop-caft}.
In fact, we will prove a slightly different result, as
formulated in the next proposition. This proposition
summarizes the coupling results, and will be instrumental both in
this paper, as well as in \cite{HHZ04b}, in which we investigate the case where
$\tau \in (2,3)$.

\begin{prop}
\label{prop-caftrep}
Fix $\tau>2$, and assume that (\ref{distribution}) holds. For any $m$ such that,
for any $\eta>0$ small enough,
    \begin{equation}
    \label{asstau}
    \prob(\sum_{j=1}^m \hat Z_j^{\smallsup{i}} \geq N^{\eta})=o(1),
    \end{equation}
there exist {\it independent} branching processes ${\cal Z}^{\smallsup{1}},
{\cal Z}^{\smallsup{2}}$, such that
    \eq
    \label{conccaft}
    \lim_{N\rightarrow \infty}
    \prob(Z_m^{\smallsup{i}}={\cal Z}_m^{\smallsup{i}})=
    1.
    \en
\end{prop}
{\bf Remark:} For fixed $m$, by the Markov inequality, (\ref{asstau}) indeed holds.
Therefore, Proposition \ref{prop-caft} follows from (\ref{conccaft}). We are left
to prove Proposition \ref{prop-caftrep}.

\proof
By (\ref{asstau}), it suffices to show that
$\prob(Z_m^{\smallsup{i}}={\cal Z}_m^{\smallsup{i}},
    \sum_{j=1}^m \hat Z_j^{\smallsup{i}} < N^{\eta})=1+o(1)$.
For this, we use Lemma \ref{lem-firststubs} to conclude that,
for $\eta<1/2$,
    \eq\prob(Z_m^{\smallsup{i}}={\cal Z}_m^{\smallsup{i}},
    \sum_{j=1}^m \hat Z_j^{\smallsup{i}} < N^{\eta})
    =\prob(\hat Z_m^{\smallsup{i}}={\cal Z}_m^{\smallsup{i}},
    \sum_{j=1}^m \hat Z_j^{\smallsup{i}} < N^{\eta})
    +o(1).
    \en
By the coupling between $\hat Z_m^{\smallsup{i}}$ and ${\cal Z}_m^{\smallsup{i}}$,
a miscoupling occurs with probability equal to $p_{\sN}$ defined in (\ref{pNdef}).
Therefore, by Remark \ref{exttau2}, the probability of a miscoupling
for the offspring of a given individual is bounded from above by $N^{-\alpha_2}$
with probability $1+O(N^{-\beta_2})$.
On the event that $\sum_{j=1}^m \hat Z_j^{\smallsup{i}} < N^{\eta}$, the
number of individuals that need to be coupled is bounded from above by
$N^{\eta}$. We thus obtain that for any $\eta<\alpha_2$,
    \eq
    \prob(\hat Z_m^{\smallsup{i}}\neq {\cal Z}_m^{\smallsup{i}},
    \sum_{j=1}^m \hat Z_j^{\smallsup{i}} < N^{\eta}, p_{\sN}\leq N^{-\alpha_2})
    \leq N^{\eta} N^{-\alpha_2}=o(1),
    \en
which completes the proof.
\qed
\end{document}